\documentclass[10pt, leqno]{amsart}

\usepackage{color}
\usepackage{amsmath,amssymb,amsfonts,amsthm,bm}
\usepackage{graphicx}
\usepackage{float}
\usepackage{needspace}
\usepackage{enumerate}
\usepackage{appendix}
\usepackage[shortlabels]{enumitem}
\usepackage{tikz}
\usepackage{xcolor}
\usepackage{hyperref}
\hypersetup{pdfborder={0 0 0}}
\usepackage[numbers, sort&compress]{natbib}
 \usepackage[left=1.35in, right=1.35in]{geometry}

\usepackage[normalem]{ulem}

\usepackage{tikz}
\usetikzlibrary{arrows.meta,positioning}
 
\def\R{\mathbb R}
\def\Z{\mathbb Z}

\def\q{\mathbb Q}

\def\N{\mathbb N}
\def\F{\mathcal F}
\def\cB{\dot{\mathcal{B}}}

\def\cR{\mathcal R}
\def\cC{\mathcal C}

\def\d{\partial}

\def\dD{\dot \Delta}
\def\dS{\dot S}
\renewcommand\div{{\rm div}\,}
\renewcommand\exp{{\rm exp}\,}
\renewcommand\log{{\rm log}\,}
\renewcommand\det{{\rm det}\,}

\newcommand{\black}{}

\newcommand{\andf}{\quad\!\hbox{and}\!\quad}
\newcommand{\ifa}{\quad\!\hbox{if}\!\quad}

\newtheorem{theorem}{Theorem}[section]
 \newtheorem{corollary}[theorem]{Corollary}
 \newtheorem{lemma}[theorem]{Lemma}
 \newtheorem{proposition}[theorem]{Proposition}
 \theoremstyle{definition}
 \newtheorem{definition}[theorem]{Definition}
 \theoremstyle{remark}
 \newtheorem{remark}{Remark}
 
 \numberwithin{equation}{section}

\newcommand{\curl}{\mbox{curl }}

\newcommand{\abs}[1]{\left\vert#1\right\vert}

\newcommand{\norm}[1]{\Vert#1\Vert}

\begin{document}

\title[Critical stability of stratified IPM flows]
{Sharp asymptotic stability of the IPM equation in $\mathbb R^2$ \\via isopycnal coordinates}

\author{Roberta Bianchini}
\address[R. Bianchini]{Consiglio Nazionale Delle Ricerche, 00185, Rome, Italy.}
\email{roberta.bianchini@cnr.it}

\author{Timothée Crin-Barat}
\address[T.~Crin-Barat]{Universit\'e de Toulouse,  Institut de Math\'ematiques de Toulouse, Route de Narbonne 118, 31062 Toulouse Cedex 9, France.}
\email{timothee.crin-barat@math.univ-toulouse.fr}

\author{Ling-Yun Shou}
\address[L-Y. Shou]{School of Mathematical Sciences, Ministry of Education Key Laboratory of NSLSCS and Key Laboratory of Jiangsu Provincial Universities of FDMTA, Nanjing Normal University, Nanjing, 210023,
P. R. China.}
\email{shoulingyun11@gmail.com}

\author{Shan Wang}
\address[S. Wang]{Instituto de Ciencias Matem\'aticas CSIC-UAM-UCM-UC3M}
\email{shan.wang@icmat.es}

\begin{abstract}
We study the incompressible porous media (IPM) equation, a model describing the motion of a fluid with variable density through a porous medium under the action of gravity. We prove that the stably stratified state, in which heavier fluid lies below lighter fluid, is nonlinearly stable: small perturbations decay over time and the fluid returns to equilibrium. Our result holds in the sharp Sobolev setting $H^s(\mathbb R^2)$ for every $s>2$: this is the optimal threshold, since the equation is ill-posed in $H^2(\mathbb R^2)$. To achieve that, we first prove global well-posedness at critical Besov regularity.  We also provide optimal decay rates.

The main difficulty is that the stabilizing effect of stratification acts only in some directions and is degenerate, so that reaching the sharp regularity threshold requires new ideas. First, we rewrite the equation in \emph{isopycnal coordinates}, which follow the level sets of the density: in these coordinates the transport becomes one-dimensional as only the horizontal term survives. Second, we establish new \emph{anisotropic commutator estimates} which allow us to control the remaining vertical transport. Finally, we show that the estimates obtained in isopycnal coordinates can be transferred back to the original variables.
\end{abstract}

\date{}
\maketitle


\section{Introduction}

Active scalar equations, in which a transported scalar determines the velocity field through a singular integral operator, are a fundamental class of nonlinear PDEs in fluid dynamics. A key distinction between these models lies in the order of the operator relating the scalar to the velocity. For the incompressible porous media (IPM) equation, as for the surface quasi-geostrophic (SQG) equation, this operator is of order zero, so that the velocity has the same regularity as the transported scalar, while for the two-dimensional Euler equations in vorticity form, the Biot-Savart law gains one derivative. The lack of such a smoothing effect makes IPM (and SQG) more singular and their well-posedness theory more delicate. We refer to \cite{HPG,CGO} and the references therein for further details.

In this work, we focus on the two-dimensional incompressible porous media equation, which models the evolution of a variable-density fluid through a porous medium under the action of gravity, in the \emph{whole plane} $\mathbb R^2$ (see \cite{bear1} for a derivation):
\begin{equation}\label{eq:ipm}\tag{IPM}
\begin{cases}
\partial_t \rho + u\cdot\nabla \rho = 0, \\[4pt]
u = -\kappa\,\nabla P + \mathbf{g}\,\rho, \qquad \mathbf{g} = (0,-g)^T \qquad \text{(Darcy's law)}, \\
\div u = 0.
\end{cases}
\end{equation}
Here $\rho$, $u$ and $P$ denote the density, velocity and pressure of the fluid, $\kappa$ is the permeability and $g$ the gravitational acceleration. Hereafter, we normalize $\kappa=g=1$.

Equation \eqref{eq:ipm} contains no dissipation: the density is purely transported by the flow. Nevertheless, an effective damping mechanism emerges when one considers perturbations of a \emph{stable stratification}, namely a background density $\bar\rho(x_2)$ satisfying $\bar\rho'(x_2)<0$, so that heavier fluid lies below lighter fluid. We focus on the simplest stable profile $\bar\rho(x_2)=1-x_2$ and write $\rho=\bar\rho+b.$
Then the perturbation $b$ satisfies the \emph{stable IPM equation}
\begin{equation}\label{eq:ipm3}
\begin{cases}
\partial_t b+u\cdot\nabla b=u_2,\\
u=\mathcal R^\perp\mathcal R_1 b,\qquad \operatorname{div}u=0,
\end{cases}
\end{equation}
where $\mathcal R_j:=(-\Delta)^{-\frac12}\partial_{x_j}$, $j=1,2$, denote the Riesz transforms and $\mathcal R^\perp:=(-\mathcal R_2,\mathcal R_1)$.
The term $u_2=\mathcal R_1^2 b$ provides the damping mechanism. Indeed, the linearized equation
\[
\partial_t b-\mathcal R_1^2 b=0
\]
has the Fourier representation
\begin{equation}\label{eq:semigroup}
\widehat b(t,\xi)
=
e^{-t\,\frac{\xi_1^2}{|\xi|^2}}\widehat b_0(\xi).
\end{equation}
The damping rate $\frac{\xi_1^2}{|\xi|^2}$ is anisotropic and degenerates near the vertical frequency axis $\{\xi_1=0\}$. This angular degeneracy governs both the decay properties of the linearized flow and the anisotropic estimates required to control the nonlinear dynamics.

The asymptotic behavior depends crucially on the spatial domain. On $\mathbb T\times\mathbb R$, the set $\{\xi_1=0\}$ corresponds to a genuine undamped horizontal Fourier mode, and solutions converge to a stratified state that is close to, but in general different from, the reference profile $\bar\rho$, see \cite{RW}. On $\mathbb R^2$, by contrast, $\{\xi_1=0\}$ has zero Lebesgue measure in Fourier space, and no nontrivial undamped $L^2$ mode persists. Accordingly, the perturbation decays to zero in suitable norms, so that $\rho(t)$ converges back to $\bar\rho$. Elgindi \cite{TME} first exploited this damping mechanism to establish the stability of small perturbations in sufficiently high Sobolev spaces.

Note that the IPM system \eqref{eq:ipm3} is invariant under the scaling
$$
 b_\lambda(t,x)=\lambda^{-1}b(t,\lambda x),
 \qquad \lambda>0.
$$
A functional space is called critical if its norm is invariant under this scaling.
In two dimensions, both $\dot H^2(\mathbb R^2)$ and
$\dot B^2_{2,1}(\mathbb R^2)$ are critical in this sense.

\smallbreak
\textbf{The goal of this work.}
 Taking advantage of the \emph{anisotropic damping} due to the background stratification, we prove \emph{ asymptotic stability in $H^{2+}(\mathbb{R}^2)$}. This threshold is sharp in the framework of Sobolev spaces, since the IPM equation is ill-posed in $H^2(\mathbb{R}^2)$ even near a stable profile \cite{RDL}. Our contribution closes the regularity gap $(2,3]$ left open for the stable IPM equation on $\mathbb{R}^2$. This result follows from a proof of global well-posedness in the critical Besov framework $\dot B^1_{2,1}\cap\dot B^2_{2,1}$, with smallness required only in the critical norm $\dot B^2_{2,1}$.
\smallbreak
\textbf{Difficulties.}
In the stable regime, the damping mechanism in IPM is strongly anisotropic: at the energy level, it controls $\mathcal R_1 b$ rather than $b$ itself. The velocity inherits this horizontal structure through $u=\mathcal R^\perp\mathcal R_1 b$, whereas the nonlinear transport $u\cdot\nabla b$ involves a full derivative of $b$. As a result, standard energy estimates do not close at the sharp Sobolev regularity $H^{2+}(\mathbb R^2)$. The best previous result proves asymptotic stability in $H^{3+}(\mathbb R^2)$, see \cite{RTM}.

A second difficulty is the decay mechanism in the whole space. The linear damping rate $\frac{\xi_1^2}{|\xi|^2}$ is homogeneous of degree zero and degenerates near $\{\xi_1=0\}$. Hence increasing the isotropic regularity does not improve the decay rates. This sharply contrasts with the setting $\mathbb T\times\mathbb R$, where decay rates improve with regularity \cite{RW}. Thus, on $\mathbb R^2$, \emph{extra smoothness does not yield extra decay}: both the sharp regularity threshold and the optimal decay rates require a genuinely anisotropic analysis of the equation.
\smallbreak
\textbf{New ideas and strategy.}
To obtain the sharp $H^{2+}(\mathbb R^2)$ result, we first prove global-in-time well-posedness at the critical regularity $\dot B^1_{2,1}\cap\dot B^2_{2,1}$. Since $H^{2+}(\mathbb R^2)\hookrightarrow \dot B^1_{2,1}\cap\dot B^2_{2,1}$, the latter Besov space provides a low-regularity framework to obtain an $L^1$-in-time control of the Lipschitz norm of the vertical velocity, by exploiting the anisotropic damping.  More precisely, we use the chain of embeddings
\begin{equation}\label{eq:embeddings}
    H^{2+} \hookrightarrow \dot B^1_{2,1}\cap\dot B^2_{2,1}
    \hookrightarrow \cB^{\frac32,\frac12}\cap\cB^{\frac12,\frac12},
\end{equation}
where $\cB^{\frac32,\frac12}$ and $\cB^{\frac12,\frac12}$ are the anisotropic Besov-type spaces introduced in \cite{RTM} and adapted to the linearized dynamics. 

Compared with \cite{RTM}, two additional ingredients are needed to cover the sharp regularity threshold $H^{2+}$. First, we introduce \emph{isopycnal coordinates} \cite{BD2024}, yielding an equation with only \emph{horizontal transport} (the vertical transport is eliminated by the change of variables). Second, we prove a new anisotropic commutator estimate that controls the remaining horizontal transport. Together, these ingredients allow us to close the argument. 

\vspace{3mm}
\textbf{Main result in Sobolev spaces.} We prove asymptotic stability in $H^s(\mathbb{R}^2)$ for every $s>2$.

\begin{theorem}
\label{thm:ipm-higher-sobolev-intro}
Let $s>2$ and $\bar\rho(x_2)=1-x_2$. There exist $\delta_s>0$ and $C_s>0$, depending only on $s$, such that the following holds. If $\rho_0=\bar\rho+b_0$ with $b_0\in H^s(\mathbb{R}^2)$ and
\[
\|b_0\|_{H^s(\mathbb R^2)}\leq\delta_s,
\]
then \eqref{eq:ipm} has a unique global-in-time solution $\rho=\bar\rho+b$ with $b\in C(\mathbb{R}_+;H^s)$, and
\begin{equation}\label{ipm-uniform-sobolev-intro}
\sup_{t\geq0}\|b(t)\|_{H^s}^2+\int_0^\infty\|\mathcal R_1 b(t)\|_{H^s}^2\,dt\leq C_s\|b_0\|_{H^s}^2 .
\end{equation}
Moreover, if
$b_0\in W^{4,1}(\mathbb R^2)$, then, for every $t\geq0$,
\[
\begin{aligned}
 \|b(t)\|_{H^2}\lesssim (1+t)^{-\frac14},
 &\quad 
 \|u_1(t)\|_{H^2}&\lesssim (1+t)^{-\frac{3}{4}}
 \quad \text{and}\quad
 \|u_2(t)\|_{H^2}&\lesssim (1+t)^{-\frac{5}{4}}.
 \end{aligned}
\]
\end{theorem}

Theorem~\ref{thm:ipm-higher-sobolev-intro} gives a Sobolev formulation of our main result, just above the critical regularity $H^2(\mathbb R^2)$, where \eqref{eq:ipm} is ill-posed. We present it in this form to align with Elgindi's original framework~\cite{TME}. Our analysis in fact yields a stronger statement: global-in-time well-posedness holds at critical Besov regularity, see Theorem~\ref{them1} below.

\section{Previous results and critical obstruction}

We recall the stability and instability results for the IPM equation.

\subsection{Local well-posedness and singular behavior of IPM}

C\'ordoba, Gancedo and Orive \cite{CGO} established the local well-posedness of
\eqref{eq:ipm} in $C^{1+\alpha}$ and $\dot H^{2+\alpha}$, $0<\alpha<1$, together
with a blow-up criterion in $\dot H^{2+}$. C\'ordoba and Gancedo \cite{CG}
treated the associated interface problem, and Friedlander, Gancedo, Sun and
Vicol \cite{FGSV} studied a singular variant of IPM.

Beyond these well-posedness results, several works have revealed the singular
character of the IPM dynamics. Weak solutions need not be unique \cite{CFG,S},
and Kiselev and Yao \cite{KY} constructed solutions exhibiting unbounded
derivative growth as $t\to\infty$. Collot, Prange and Tan~\cite{CPT2026} proved the stability
of self-similar blow-up for a class of
infinite-energy solutions. For the forced equation, finite-time
singularities were constructed for spatially smooth forcing in \cite{CMZ24},
and subsequently extended to forcing that is smooth in both space and time up
to the blow-up time in \cite{ABC2026}. Very recently, the forcing
was removed in \cite{Dai2026} and \cite{Zheng2026}. 

\subsection{Large-time stability near a stratified density}

For the stable profile $\bar \rho(x_2)=1-x_2$, Elgindi \cite{TME} proved
nonlinear asymptotic stability in $\mathbb R^2$ for small, highly
regular perturbations. Castro, C\'ordoba and Lear \cite{CCL}
obtained asymptotic stability in a flat strip. Bianchini, Crin-Barat
and Paicu \cite{RTM} subsequently proved global existence in the
whole space for data in $\dot H^{1-\tau}\cap\dot H^s$, where
$0<\tau<1$ and $s>3$. Jo and Kim~\cite{JK1} and Park~\cite{JP} established
quantitative decay estimates in Sobolev and $L^2$ norms
on $\mathbb R^2$ and in periodic channels, respectively. More recently, Bianchini, Jo, Park and Wang
\cite{RW} proved the global stability in the sharp Sobolev regime $H^{2+}$ on
$\mathbb T\times\mathbb R$ for general
stratified steady states $\rho_s(x_2)$ with $\rho'(x_2)<0$. Very recently, \cite{BMN} established a stability-instability
dichotomy for stratified steady states in the periodic
channel.

The same stable stratification appears in related fluid equations.
For two-dimensional Boussinesq flows with velocity damping, global
stability and large-time behavior have been studied in
\cite{castro3,wan,BianchiniNatalini2DB,KL}. In particular,
\cite{RTM} proves that the stable IPM equation can be derived as a relaxation limit of
the damped Boussinesq system. Further results on stratified fluids are, for instance, \cite{DLS,BBCD,JW}.

\subsection{The critical endpoint}

Lipschitz control is central to the propagation of regularity for transport equations. In two dimensions, $H^2$ carries the critical number of
derivatives, yet it fails to embed into $W^{1,\infty}$. This gap is
not a technicality: Bianchini, C\'ordoba and Mart\'inez-Zoroa
\cite{RDL} proved nonexistence and strong ill-posedness for
arbitrarily small $H^2$ perturbations of $-x_2$, even though the
perturbation remains bounded in $C^1$. A similar ill-posedness result holds for
SQG \cite{JK}. In particular, the anisotropic damping generated by the
stable IPM profile is not enough to have well-posedness in the critical space $H^2$. The
question then is whether the asymptotic stability can be proved in (another) critical space that controls the Lipschitz norm. This work gives a positive answer to this question. Furthermore, it proves asymptotic stability for small data in $H^s(\mathbb R^2)$ with $s \in (2, 3]$ which remained open until now.

\subsection{Plan of the paper and notation}

Section~\ref{section2} states the well-posedness and decay results. In Section~\ref{sec:strategy}, we present the strategy of the proof.
Section~\ref{sectionabv} introduces the anisotropic functional
framework, Section~\ref{dcsec} develops the a priori estimates in
isopycnal coordinates, and Section~\ref{IPM} proves the well-posedness results.
Section~\ref{sec:linear-decay} proves the linear decay estimates, and
Section~\ref{sec:decay} treats the nonlinear problem.
The appendices contain the technical estimates.

Throughout, $A\lesssim B$ means $A\leq CB$, where $C$ is independent
of the relevant quantities and may change from line to line. We write
$x=(x_1,x_2)\in\mathbb R_h\times\mathbb R_v=\mathbb R^2$ and use
$L^p=L^p(\mathbb R^2)$, $L_h^p=L^p(\mathbb R_h)$ and
$L_v^p=L^p(\mathbb R_v)$, with analogous notation for Sobolev
spaces. For a Banach space $X$ and $T>0$, let
$L^p_T(X)=L^p((0,T);X)$.

\section{Main results}\label{section2}
We first state the critical well-posedness theorem, then the decay
theorems and the mechanism behind their proofs.

\subsection{Local and global well-posedness results}

\begin{theorem}[Local well-posedness]
\label{them1:local}
Assume that $b_0\in\dot B^1_{2,1}\cap\dot B^2_{2,1}$. Then there
exists $T>0$ such that the stratified IPM system \eqref{eq:ipm3}, with
initial datum $b_0$, admits a unique solution satisfying
\[
 b\in C([0,T];\dot B^1_{2,1}\cap\dot B^2_{2,1}).
\]
\end{theorem}

Here $\dot B^s_{2,1}$ denotes the usual homogeneous Besov space. Its definition is recalled in Section \ref{sectionabv}. Next, we state our global existence result, which requires smallness only at the critical Lipschitz level, while the lower-order norm may be arbitrarily large.

\begin{theorem}[Global {{stability}}]
\label{them1}
Assume that $b_0\in\dot B^1_{2,1}\cap\dot B^2_{2,1}$. There exists
$\delta_0>0$ such that, if
\begin{align}\label{a1}
\|b_0\|_{\dot B^2_{2,1}}\leq\delta_0,
\end{align}
then system \eqref{eq:ipm3}, with initial datum $b_0$, admits a unique
global solution satisfying
\begin{align*}
 &b\in C(\R_+;\dot B^1_{2,1}\cap\dot B^2_{2,1}),\qquad
 \cR_1b\in L^2(\R_+;\dot B^1_{2,1}\cap\dot B^2_{2,1}),\qquad u_2,\ \nabla u_2\in L^1(\R_+;L^\infty).
\end{align*}
For all $T>0$ and $s=1,2$, we have
\begin{equation}\label{r1}
 \|b\|_{L^\infty_T(\dot B^s_{2,1})}
 +\|\cR_1b\|_{L^2_T(\dot B^s_{2,1})}
 \leq C\|b_0\|_{\dot B^s_{2,1}}.
\end{equation}
The vertical velocity also satisfies
\begin{equation}\label{vel1}
 \int_0^T\|u_2(t)\|_{L^\infty}\,dt\leq C\|b_0\|_{\dot B^1_{2,1}}\quad \text{and}\quad 
 \int_0^T\|\nabla u_2(t)\|_{L^\infty}\,dt\leq C\|b_0\|_{\dot B^2_{2,1}}.
\end{equation}
\end{theorem}

The asymptotic stability in sharp Sobolev regularity $H^{s}(\mathbb R^2)$ with $s>2$ is stated below.

\begin{theorem}[Propagation of higher Sobolev regularity]
\label{thm:ipm-higher-sobolev}
Let $s>2$. There exists $0<\delta_s\leq\delta_0$, depending only on $s$,
such that the following holds. If $b_0\in H^s(\R^2)$ and
\begin{equation*}
 \|b_0\|_{\dot B^2_{2,1}}\leq\delta_s,
\end{equation*}
then the global solution $b$ provided by Theorem~\ref{them1} satisfies
$b\in C(\R_+;H^s)$, $\cR_1 b\in L^2(\R_+;H^s)$ and
\begin{equation}\label{ipm-uniform-sobolev}
 \sup_{t\geq0}\|b(t)\|_{H^s}^2
 +\int_0^\infty\|\cR_1b(t)\|_{H^s}^2\,dt
 \leq C_s\|b_0\|_{H^s}^2,
\end{equation}
where $C_s$ depends only on $s$.
\end{theorem}

\begin{remark}
Theorem~\ref{them1} should be compared with the strong ill-posedness in
$H^2$ proved in \cite{RDL}. Although $H^2$ and $\dot B^2_{2,1}$ have the
same derivative count in dimension two, the following embedding is crucial:
\[\dot B^{\frac{d}{2}}_{2,1}(\mathbb{R}^d) \hookrightarrow L^\infty(\mathbb R^d).\]
Moreover,
Theorem~\ref{thm:ipm-higher-sobolev} shows that any additional Sobolev
regularity is propagated globally in time.
\end{remark}

\begin{remark}[On the intersection $\dot B_{2,1}^1\cap \dot B^2_{2,1}$]
In Theorem \ref{them1}, the critical
$\dot B^2_{2,1}$ norm controls the Lipschitz regularity and is the only norm
required to be small. The $\dot B^1_{2,1}$ norm controls the low frequencies,
makes the intersection $\dot B_{2,1}^1\cap \dot B^2_{2,1}$ a Banach space and is used to prove the uniqueness.
\end{remark}

\subsection{Time-decay estimates}

We now state the optimal decay estimates under an additional
condition on the initial data. We define the mixed
norm by
\begin{equation}\label{def:mixed-X}
 \|f\|_{\mathcal X^{s,r}}
 :=\sum_{j\in\mathbb Z}2^{sj}
       \sup_{k\in\mathbb Z}
       \bigl(2^{rk}\|\dot\Delta_j\dot\Delta_k^h f\|_{L^2}\bigr),
 \qquad s,r\in\mathbb R.
\end{equation}
{Here $\dot\Delta_j$ and $\dot\Delta_k^h$ localize to
$|\xi|\sim2^j$ and $|\xi_1|\sim2^k$, respectively, see
Section~\ref{sectionabv} for their definitions.}

\begin{theorem}[Nonlinear decay]\label{thm:decay}
Under the assumptions of Theorem~\ref{them1}, let $(b,u)$ be the global
solution of \eqref{eq:ipm3}. Set $h_0:=b_0\circ X(0,\cdot),$
where $X(0,\cdot)$ denotes the inverse of the map $x\longmapsto (x_1,x_2-b_0(x)).$
Assume that
\[
h_0\in \mathcal X^{\frac 12,-\frac 12}\cap \mathcal X^{\frac 52,-\frac 12}.
\]
Then, for every $t\geq0$, {{there holds}}
\[
\|b(t)\|_{H^2}
 \leq C(1+t)^{-\frac14},
\quad
\|u_1(t)\|_{H^2}
 \leq C(1+t)^{-\frac{3}{4}},
\quad
\|u_2(t)\|_{H^2}
 \leq C(1+t)^{-\frac{5}{4}}
\]
as well as
\[
\|\nabla b(t)\|_{L^\infty}
 \leq C(1+t)^{-\frac12},
\quad
\|\nabla u_1(t)\|_{L^\infty}
 \leq C(1+t)^{-1},
\quad
\|\nabla u_2(t)\|_{L^\infty}
 \leq C(1+t)^{-\frac{3}{2}}.
\]
In all these estimates, the constant $C$ is independent of $t$.
\end{theorem}

\begin{remark}[Sharpness of the decay rates]\label{rem:decay-elgindi}
The decay rates obtained here are sharp, as they coincide with the optimal rates for the linear equation
\[
\partial_t b-\mathcal R_1^2 b=0;
\]
see Theorem~\ref{thm:linear-decay}. Elgindi~\cite{TME} obtained the same rates, for small initial data in $W^{4,1}\cap H^s$, with $s\geq20$. Note that the following embedding holds: \[
W^{4,1}(\mathbb R^2)
\hookrightarrow
\mathcal X^{\frac12,-\frac12}
\cap
\mathcal X^{\frac52,-\frac12},
\]
see Lemma~\ref{lem:mixed-embeddings}. 
\end{remark}

\begin{remark}
The spaces $\mathcal X^{s+\sigma,-\sigma}$ may be viewed as anisotropic substitutes of $W^{4,1}$: they quantify the concentration of the Fourier transform near the degenerate region $\{\xi_1=0\}$. The parameter $\sigma$ strengthens only the angular control and does not increase regularity. This is adapted to the dissipative symbol
\[
-\frac{\xi_1^2}{|\xi|^2},
\]
which is homogeneous of degree zero and depends only on the angle of $\xi$. Roughly, the parameter $\sigma$ provides additional decay without requiring additional regularity.
\end{remark}

\section{Proof Strategy}\label{sec:strategy}
\subsection{Overall Strategy}
{The overall strategy to prove global well-posedness is as follows.
\begin{itemize}

    \item We first establish a priori estimates in the critical isotropic space
$\dot B^1_{2,1}\cap \dot B^2_{2,1}.$
    Proposition~\ref{dc0} shows that solutions lie in such a space provided that the vertical velocity satisfies
    \[
    u_2\in L^1_T(W^{1,\infty}).
    \]
    The Sobolev result follows by
    \[
    H^{2+\varepsilon}
    \hookrightarrow
    \dot B^1_{2,1}\cap \dot B^2_{2,1}.
    \]

  \item The main issue is to recover the Lipschitz control of \(u_2\) at this low regularity. Following the approach of \cite{RTM}, we exploit the anisotropic structure of the dissipation. In particular, for the linear equation $\d_t b-\cR_1^2b=0$, one has
\[
\|\nabla u_2\|_{L^1_T(L^\infty)}
\lesssim
\|b\|_{L^1_T(\cB^{-\frac12,\frac52})}
\lesssim
\|b_0\|_{\cB^{\frac32,\frac12}}.
\]
In particular, the Lipschitz control of \(u_2\) can be recovered at a regularity level strictly below what would be required by a purely isotropic argument.
Moreover, the embedding
\[
\dot B^1_{2,1}\cap\dot B^2_{2,1}
\hookrightarrow
\cB^{\frac32,\frac12}\cap\cB^{\frac12,\frac12}
\]
allows us to go back to the isotropic framework. 

    \item The anisotropic estimate cannot, however, be closed directly in Eulerian coordinates. The difficulty comes from the transport term: at the critical regularity
    $\cB^{\frac32,\frac12}\cap\cB^{\frac12,\frac12},$
    the anisotropic commutator estimates behave badly with the vertical transport $u_2\partial_2 b$ (see Lemma \ref{lem:2.97}).
    In particular, the vertical derivative produces an angular loss that cannot be absorbed by the anisotropic dissipation.

    \item The key step is to pass to \emph{isopycnal coordinates}. Setting
    \[
    h:=b\circ X,
    \]
    with \(X\) given by Lemma~\ref{lem:density-reduction}, equation~\eqref{eq:ipm3} is transformed into
    \[
    \partial_t h+U_1\partial_1 h-\cR_1^2h=G(h).
    \]
    {Here $U:=u\circ X$ and $G(h):=U_2-\cR_1^2h$.}
    The crucial point is that the vertical transport has been canceled by this change of coordinates: the only remaining transport term is the horizontal one, $U_1\partial_1 h$.

    \item This reduction makes it possible to exploit the anisotropic structure. In Proposition~\ref{hcom}, we establish a commutator estimate specifically adapted to the horizontal transport. It allows us to close the estimates for \(h\) in
$\cB^{\frac32,\frac12}\cap\cB^{\frac12,\frac12},$
    together with the dissipative control in $L^1_T(\cB^{-\frac12,\frac52}).$

    \item Finally, we return to Eulerian coordinates. Lemma~\ref{dcmp} shows that {the norms of \(b\) and \(h\) are in fact equivalent} in the critical isotropic space $\dot B^1_{2,1}\cap\dot B^2_{2,1}.$
   Such an equivalence is not available in the anisotropic spaces, which is why the anisotropic analysis is performed in isopycnal coordinates while the final bootstrap is closed in the isotropic critical norm. Combining these estimates yields the required {bounds for $u_2$ in $L^1_T(W^{1,\infty})$},
    and therefore closes the global a priori estimates.

\end{itemize}
The decay estimates require additional control near $\xi_1 = 0$,
where the dissipation {degenerates. This is ensured by}
$b_0 \in W^{4,1}$,
or in fact by a slightly weaker assumption.}
A more detailed description of the proof strategy is provided in Sections~\ref{sec:stratGWP} and~\ref{sec:stratDecay}.

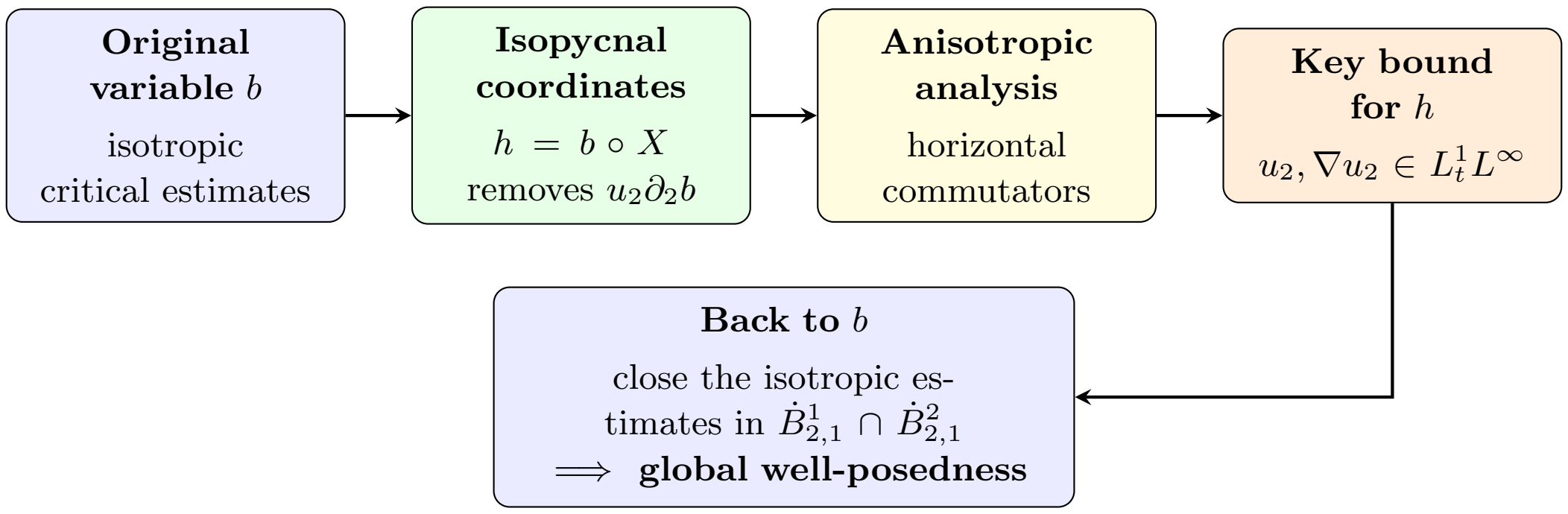
\begin{figure}[htbp]
\centering
\begin{tikzpicture}[
     >=stealth,
     every node/.style={align=center,font=\small},
     box/.style={
        draw,
         rounded corners,
        minimum height=1.35cm,
         text width=2.65cm,
         inner sep=5pt
     }
 ]

 \node[box,fill=blue!8] (b) at (-5.4,0)
 {\textbf{Original variable $b$}\\[1mm]
 isotropic \\critical estimates};

 \node[box,fill=green!10] (h) at (-1.8,0)
 {\textbf{Isopycnal coordinates}\\[1mm]
 $h=b\circ X$\\
 removes $u_2\partial_2 b$};

 \node[box,fill=yellow!15] (aniso) at (1.8,0) {\textbf{Anisotropic analysis}\\[1mm]
 horizontal commutators};

 \node[box,fill=orange!15] (lip) at (5.4,0)
 {\textbf{Key bound for $h$}\\[1mm]
 $u_2,\nabla u_2\in L^1_tL^\infty$};

 \node[
     box,
     fill=blue!8,
     text width=4.8cm
 ] (back) at (0,-2.5)
 {\textbf{Back to $b$}\\[1mm]
 close the isotropic estimates in
 $\dot B^1_{2,1}\cap\dot B^2_{2,1}$\\
 \textbf{$\implies$ global well-posedness}};
 \draw[->,thick] (b) -- (h);
 \draw[->,thick] (h) -- (aniso);
 \draw[->,thick] (aniso) -- (lip);
 \draw[->,thick] (lip.south) |- (back.east);

 \end{tikzpicture}
 \caption{Sketch of the strategy of the proof.}
 \label{fig:overall-strategy}
 \end{figure}

\subsection{Global well-posedness}\label{sec:stratGWP}
\subsubsection{The key quantities to control}
We now explain the mechanism behind the global estimate in Theorem~\ref{them1}. The key point is to derive a priori estimates yielding
\begin{equation}\label{vel1-strategy}
u_2\in L^1_T(L^\infty)
\qquad\text{and}\qquad
\nabla u_2\in L^1_T(L^\infty),
\end{equation}
with bounds that are uniform with respect to $T$.

The second estimate controls the vertical deformation generated by the transport term. 
The first one is also key. Indeed, since no $L^2$ control is imposed at low frequencies (in the most general Theorem \ref{them1}), the $L^\infty$ norm of $u_2$ cannot be recovered from $\nabla u_2$ alone and must therefore be estimated separately. Together, the two bounds in \eqref{vel1-strategy} allow to close the nonlinear estimates and to justify the change to isopycnal coordinates. 

\subsubsection{The role of anisotropic Besov spaces}\label{sec:roleofaniso}

In our critical framework, the aim is therefore to establish the two bounds in \eqref{vel1}. The starting point is the anisotropic dissipation of $-\cR_1^2$. On the Fourier support of $\dot\Delta_j\dot\Delta_k^h b$, one has
\[
|\xi|\sim 2^j,
\qquad
|\xi_1|\sim 2^k.
\]
Hence, by Plancherel's identity and Bernstein's inequalities,
\[
\begin{aligned}
-\int_{\R^2}\cR_1^2\dot\Delta_j\dot\Delta_k^h b\,
\dot\Delta_j\dot\Delta_k^h b\,dx
=
\|\cR_1\dot\Delta_j\dot\Delta_k^h b\|_{L^2}^2\gtrsim
2^{2k-2j}
\|\dot\Delta_j\dot\Delta_k^h b\|_{L^2}^2.
\end{aligned}
\]
Thus, for the linear equation
\[
\partial_t b-\cR_1^2b=0,
\]
the dissipation yields $L^1$-in-time control of each anisotropic block with gain $2^{2k-2j}$. Summing with the appropriate weights gives
\[
\|b\|_{L^1_T(\cB^{-\frac12,\frac52})}
\lesssim
\|b_0\|_{\cB^{\frac32,\frac12}}
\lesssim
\|b_0\|_{\dot B^2_{2,1}},
\]
where $\cB^{s_1,s_2}$ denotes the anisotropic Besov space introduced in Section~\ref{sectionabv}. Since
\[
u_2=\cR_1^2 b,
\]
the anisotropic Bernstein inequalities then imply
\[
\|\nabla u_2\|_{L^1_T(L^\infty)}
\lesssim
\|b\|_{L^1_T(\cB^{-\frac12,\frac52})}.
\]

Similarly
\[
\|u_2\|_{L^1_T(L^\infty)}
\lesssim
\|b\|_{L^1_T(\cB^{-\frac32,\frac52})}
\lesssim
\|b_0\|_{\cB^{\frac12,\frac12}}
\lesssim
\|b_0\|_{\dot B^1_{2,1}}.
\] 

The main difficulty is to propagate these anisotropic estimates in the presence of the transport term
\[
u_2\partial_{x_2}b.
\]
Indeed, frequency localization introduces the commutator
$[u_2,\dot\Delta_j\dot\Delta_k^h]\partial_{x_2}b$.
For example, if $|\xi|\simeq2^j$ and
$|\xi_1|\simeq2^k$, the horizontal commutator estimate
in \eqref{eq:esr11} contains $2^{-k}$, while
the vertical derivative is estimated by $2^j$.
The resulting estimate therefore contains
\[
 2^{-k}2^j=2^{j-k}.
\]
This becomes arbitrarily large when $k\ll j$. In other words, the vertical transport generates an \emph{angular derivative loss} near the degenerate region $\xi_1=0$.
In contrast, the horizontal transport can be handled within the anisotropic energy estimates. This is the reason for treating the two directions differently: we develop new commutator estimates for the horizontal transport and remove the problematic term $u_2\partial_{x_2}b$ through the change to isopycnal coordinates.

\subsubsection{The role of isopycnal coordinates}
To this end, we use the identity
\[
 (\partial_t+u\cdot\nabla_x)(x_2-b)=0
\]
and introduce \emph{isopycnal coordinates}, in the spirit of nonlinear changes of
variables used in the analysis of \emph{eddy diffusivity} models \cite{BD2024} and in the study of viscous fluids in critical spaces (see, for instance,
\cite{CMZ0,RLG}). We set
\[
 Y(t,x){:=}(x_1,x_2-b(t,x)),\qquad X{:=}Y^{-1},\qquad h{:=}b\circ X.
\]
{Since $y_2=x_2-b={1-\rho}$, points with the same
density have the same $y_2$ coordinate. This is the isopycnal
coordinate transformation used in \cite[Section~2.1]{BD2024}.

When $b=0$, the surface of density $\rho={1-y_2}$ is the
horizontal line $x_2=y_2$. In the presence of a perturbation,
the same density surface is located at $x_2=y_2+h(t,y)$.
Thus, $h$ measures how far this surface has moved vertically.
Moreover, $J{:=}\partial_2X_2=1+\partial_2h$ measures the change
in physical height corresponding to a unit change in $y_2$
\cite{GM1990,Griffies2004}.}
The second coordinate is constant along fluid trajectories, so the
transformed equation \eqref{dc6} contains only horizontal transport. Precisely, as in Lemma \ref{lem:density-reduction}, defining  $U{:=}u\circ X$, we can rewrite the IPM system \eqref{eq:ipm3} as
\[
 \partial_t h+U_1\partial_1h-\cR_1^2h=G(h),
\]
with the change-of-coordinates error
\[
 G(h){:=}U_2-\cR_1^2h
 =\bigl[\cR_1^2(h\circ Y)\bigr]\circ X-\cR_1^2h.
\]
A key novelty is the development of new
anisotropic commutator estimates \eqref{hc1}--\eqref{hc2}
in Proposition~\ref{hcom}, which
allow us to establish maximal regularity in anisotropic spaces and to obtain the  $L^1$-in-time Lipschitz bound for the vertical velocity. However, a new difficult term $G(h)$ appears, since the nonlocal Riesz operator 
$\mathcal R_1^2$ does not commute with composition by $X$. In fact, the classical estimates for
operators in isotropic spaces, e.g., \cite{CMZ0} do not estimate $G(h)$ in the \emph{anisotropic} norm. 

To overcome this difficulty, we prove the anisotropic estimate \eqref{dcb} for $G(h)$, which allows us to absorb the coordinate
error into the dissipation
(see Appendix~\ref{app:dc-riesz}). This leads to the
anisotropic estimate \eqref{dc2} for $h$ and to the
uniform $L^1$-in-time Lipschitz estimate for $u_2$. Using these estimates, we derive the isotropic bounds
\eqref{dien}.
Returning to the original variables,
we obtain the desired a priori estimates
\eqref{db2}-\eqref{db1}.

\subsection{Decay estimates}\label{sec:stratDecay}

We develop a new energy method for the IPM equation adapted
to its anisotropic dissipation. For the linear problem, the Fourier multiplier
$e^{-t\frac{\xi_1^2}{|\xi|^2}}$ of $e^{t\mathcal R_1^2}$ shows that
the damping depends on the angular ratio $\frac{|\xi_1|}{|\xi|}$.
The key observation is the bound
\[
 \left(\frac{|\xi_1|}{|\xi|}\right)^\sigma
 e^{-t\frac{\xi_1^2}{|\xi|^2}}
 \lesssim (1+t)^{-\frac{\sigma}{2}},
 \qquad t\geq0,\quad \sigma>0.
\]
An angular weight of order $\sigma$ in the initial data yields
the decay rate $(1+t)^{-\frac{\sigma}{2}}$ (and does not change the regularity requirements).
Using this multiplier and the anisotropic Littlewood-Paley
decomposition, we obtain the optimal decay estimates in
Theorem~\ref{thm:linear-decay}.

For the nonlinear problem, we seek a lower bound for the
dissipation in terms of the $L^2$ norm.
To this end, we prove the following anisotropic interpolation estimate:
\[
 \|h\|_{L^2}
 \lesssim
 \|h\|_{\mathcal X^{\frac12,-\frac12}}^{\frac23}
 \|\cR_1h\|_{L^2}^{\frac{1}{3}}.
\]
Transferring this estimate to $b$ and using the  basic energy identity, we obtain
\[
 \frac{d}{dt}\|b(t)\|_{L^2}^2
 +\frac{c}{\|h(t)\|_{\mathcal X^{\frac12,-\frac12}}^4}
  \|b(t)\|_{L^2}^6\leq0.
\]
We need to propagate the angular bound in
$\mathcal X^{\frac12,-\frac12}$ uniformly in time.
Integrating the resulting differential inequality yields
$\|b(t)\|_{L^2}\lesssim(1+t)^{-\frac14}$.
This argument is inspired by classical decay methods for fluid
equations (see, e.g., \cite{GuoWang2012}).

The difficulty in the higher-order estimates is to control
the transport commutators in the anisotropic norms.
Further localized energy estimates establish the gradient
decay \eqref{velocity-gradient-decay}
and the decay of $(h,U)$ in $\dot B^2_{2,1}$.
The composition estimate \eqref{dcn} transfers these bounds
to $(b,u)$.
Together with the $L^2$ estimates, this proves the desired
$H^2$ decay. See Sections~\ref{sec:linear-decay} and~\ref{sec:decay}
for the detailed proofs.

\section{Anisotropic functional spaces and anisotropic commutator estimates}\label{sectionabv}
We define the frequency spaces adapted to the directional
dissipation and state the product and commutator laws used below.
Their proofs, except the ones that can be found in \cite{BCD}, are collected in Appendices~\ref{sec:appendixBesov}
and~\ref{app:com}.

\subsection{Isotropic and anisotropic homogeneous Besov spaces}
Choose a
smooth, radial, nonincreasing function $\chi$, supported in
$B(0,\frac{4}{3})$ and equal to one on $B(0,\frac{3}{4})$, and set
\[
 \varphi(\xi):=\chi(\frac{\xi}{2})-\chi(\xi).
\]
Then $\varphi$ satisfies
\[
 \sum_{j\in\mathbb Z}\varphi(2^{-j}\xi)=1\quad(\xi\neq0),
 \qquad
 \operatorname{supp}\varphi\subset
 \left\{\xi\in\mathbb R^d:\frac34\leq|\xi|\leq\frac83\right\}.
\]
For $j,k,q\in\mathbb Z$, we use the following notation:
\begin{enumerate}
 \item $\dD_j$ is the Littlewood-Paley block associated with the full
 Fourier variable $\xi$.
 \item $\dot\Delta_k^h$ is the block associated with the horizontal
 Fourier variable $\xi_h=(\xi_1,\ldots,\xi_{d-1})$.
 \item $\dot\Delta_q^v$ is the block associated with the vertical
 Fourier variable $\xi_d$.
\end{enumerate}
For $f\in\mathcal S'(\R^d)$, we set
\begin{equation}\label{eq:def}
    \begin{aligned}
      &\dot\Delta_ jf:=\F^{-1}(\varphi(2^{-j}\abs{\xi})\hat{f}) \andf \dot S_j f:=\F^{-1}(\chi(2^{-j}\abs{\xi})\hat{f}),\\
    &\dot\Delta_k^h f:=\F^{-1}(\varphi(2^{-k}|\xi_h|)\hat{f}) \andf \dot S_k^h f:=\F^{-1}(\chi(2^{-k}|\xi_h|)\hat{f}),\\
      &\dot\Delta_q^v f:=\F^{-1}(\varphi(2^{-q}|\xi_d|)\hat{f}) \andf \dot S_q^v f:=\F^{-1}(\chi(2^{-q}|\xi_d|)\hat{f}).\\
\end{aligned}
\end{equation}
These dyadic operators satisfy the usual almost-orthogonality properties:
there exists $N_0\in\mathbb N$ such that
\begin{equation}\label{eq:orth}
\left\{\begin{aligned}
& \dot \Delta_j \dot\Delta_k f\equiv 0  \ifa \abs{k-j}\geq 2, \\
&\dot\Delta_k(\dot S_{j-1}f \dot\Delta_j f)\equiv 0 \ifa \abs{k-j}\geq 5,\\
& \dot\Delta_k(\dot S_{j+2}f \dot\Delta_j f)\equiv 0 \ifa k> j+N_0,
\end{aligned}\right.
\end{equation}
and similarly for the horizontal and vertical blocks. Their frequency
localization also implies, after increasing $N_0$ if necessary, that
\begin{equation}\label{eq:supp-mixed}
\dot\Delta_j\dot\Delta_k^h f= \dot\Delta_j\dot\Delta_q^v f\equiv0
\qquad\text{for } k>j+N_0 \quad \text{and } q>j+N_0 .
\end{equation}
\begin{remark}
\label{rem:supjm1}
 $\mathcal{F}(\dS_{j-1}u\,\dD_j v)$ is supported in $2^j\tilde{\mathcal{C}}$, where $\tilde{\mathcal{C}}{:=}B(0,\frac{2}{3})+\mathcal{C}$ and $\mathcal{C}{:=}\{\xi\in\mathbb{R}^d \;:\; \frac{3}{4} \leq |\xi| \leq \frac{8}{3}\}$. Note that the annulus $\tilde{\mathcal{C}}$ has an inner radius $\frac{1}{12}$ and an outer radius $\frac{10}{3}$.
   $ \mathcal{F}(\dS_{j+2}u\dD_j v)$ is supported in $2^j\tilde{\mathcal{B}}$, where $\tilde{\mathcal{B}}{:=}B(0, \frac{16}{3})+\mathcal{C}$ and the ball $\tilde{\mathcal{B}}$ has the radius $8$.
\end{remark}
We next recall the classical homogeneous Besov spaces (see \cite{BCD}).
  \begin{definition}[Homogeneous Besov spaces]
Let $(p,r)\in[1,\infty]^2$, $s\in\R$ and
$f\in\mathcal S'(\R^d)$ such that
$\lim_{j\to-\infty}\|\dot S_jf\|_{L^\infty}=0$. We set
\[
 \|f\|_{\dot B^s_{p,r}}
 :=\left\|\left(2^{js}\|\dot\Delta_jf\|_{L^p}\right)_{j\in\Z}
 \right\|_{\ell^r}
\]
and denote by $\dot B^s_{p,r}$ the space associated with this norm.
\end{definition}
We shall repeatedly use the equivalent characterization
\[
 \|\dot\Delta_jf\|_{L^p}\leq Cc_j2^{-js}{\|f\|_{\dot B^s_{p,r}}},
 \qquad j\in\Z,
 \qquad \|(c_j)\|_{\ell^r}=1,
\]
for a suitable nonnegative sequence $(c_j)_{j\in\Z}$.

We now introduce homogeneous anisotropic spaces, following \cite{BCD, CPZ, LZCPAM, MP}.

\begin{definition}[Homogeneous anisotropic Besov spaces]\label{def:anib}
Let $(p,r)\in[1,\infty]^2$, $s_1,s_2\in\R$  and
$f\in\mathcal S'(\R^d)$ such that
$\lim_{j\to-\infty}\|\dot S_jf\|_{L^\infty}=0$. We set
\begin{align*}
&\norm{f}_{\cB^{s_1,s_2}_{p,r}(\R^{d})}{:=}\norm{(2^{
    js_1}2^{ k s_2}\norm{\dot{\Delta}_{j}\dD^h_k  f}_{L^{p}(\R^{d})})_{j,k\in \Z}}_{\mathit{\ell}^{r}(\Z^2)}.
\end{align*}
and denote by $\cB^{s_1,s_2}_{p,r}(\R^d)$ the space associated to this norm. Throughout
the rest of the paper, we use the notation
\[
\cB^{s_1,s_2} {:=} \cB^{s_1,s_2}_{2,1}.
\]
 \end{definition}

We next recall the Bernstein inequalities (see Lemma~6.10 in \cite{BCD}).
\begin{lemma}[Bernstein inequalities]
\label{lem:BIA}
Let $\cB_h$ and $\cB_v$ be balls in $\R_h^{d-1}$ and $\R_v$,
respectively, and let $\cC_h$ and $\cC_v$ be annuli in the corresponding
spaces. Assume that $1\leq p_2\leq p_1\leq\infty$ and
$1\leq q_2\leq q_1\leq\infty$. Then the following estimates hold:
\begin{enumerate}
    \item If $\operatorname{supp}\widehat f\subset (2^k\cB_h)\times\mathbb R_v$, then, for every multi-index $\alpha\in\mathbb N_0^{d-1}$, we have
    \begin{equation*}
        \|\d_{x_h}^{\alpha} f\|_{L^{p_1}_{h}(L^{q_1}_{v})}\lesssim 2^{k(\abs{\alpha}+(d-1)(\frac{1}{p_2}-\frac{1}{p_1}))}\|f\|_{L^{p_2}_{h}(L_{v}^{q_1})}\cdotp
    \end{equation*}
    \item If $\operatorname{supp}\widehat f\subset \mathbb R_h^{d-1}\times(2^k\cB_v)$, then, for every integer $\beta\geq0$, we have
    \begin{equation*}
        \|\d_{x_d}^{\beta} f\|_{L^{p_1}_{h}(L^{q_1}_{v})}\lesssim 2^{k(\abs{\beta}+(\frac{1}{q_2}-\frac{1}{q_1}))}\|f\|_{L^{p_1}_{h}(L_{v}^{q_2})}\cdotp
    \end{equation*}
   \item  If $\operatorname{supp}\widehat f\subset (2^k\cC_h)\times\mathbb R_v$, then, for $N\in \N$, we have
    \begin{equation*}
        \|f\|_{L^{p_1}_{h}(L^{q_1}_{v})}\lesssim 2^{-kN}\sum_{\abs{\alpha}=N}\|\d^{\alpha}_{x_h} f\|_{L^{p_1}_{h}(L_{v}^{q_1})}\lesssim \|f\|_{L^{p_1}_{h}(L^{q_1}_{v})}\cdotp
    \end{equation*}
 \item   If $\operatorname{supp}\widehat f\subset \mathbb R_h^{d-1}\times(2^k\cC_v)$, then, for $N\in \N$, we have
    \begin{equation*}
        \|f\|_{L^{p_1}_{h}(L^{q_1}_{v})}\lesssim 2^{-kN}\|\d^{N}_{x_d} f\|_{L^{p_1}_{h}(L_{v}^{q_1})}\lesssim \|f\|_{L^{p_1}_{h}(L^{q_1}_{v})}\cdotp
    \end{equation*}
\end{enumerate}
\end{lemma}
Additional properties are collected in Appendix~\ref{sec:appendixBesov}.
Here we just recall that, for $d=2$,
it follows from Lemma~\ref{lem:emb} that
    \begin{align*}  \cB^{\frac12,\frac12}_{2,1}\hookrightarrow L^{\infty} \quad \text{and} \quad \cB^{\frac32,\frac12}_{2,1}\hookrightarrow \dot{W}^{1,\infty}.
    \end{align*}
Using Lemma~\ref{prop:esu}, we obtain, for
$f\in\mathcal S'(\R^2)$,
\[
 \|\cR_1f\|_{\cB^{s_1,s_2}_{p,r}}
 \simeq\|f\|_{\cB^{s_1-1,s_2+1}_{p,r}} \quad \text{and} \quad
 \|\cR_2f\|_{\cB^{s_1,s_2}_{p,r}}
 \lesssim\|f\|_{\cB^{s_1,s_2}_{p,r}}.
\]

Next, we recall the Chemin--Lerner norms
(\cite{chemin1}) and introduce their anisotropic version.
\begin{definition}[Chemin--Lerner norms]\label{defntimespace}
Let $T>0$, $\varrho,p,r\in[1,\infty]$ and
$s,s_1,s_2\in\mathbb R$.
We define the isotropic and anisotropic Chemin--Lerner norms by
\[
 \|f\|_{\widetilde L_T^\varrho(\dot B^s_{p,r})}
 :=\left\|
 \left(2^{js}
 \|\dot\Delta_jf\|_{L_T^\varrho(L^p)}
 \right)_{j\in\mathbb Z}
 \right\|_{\ell^r(\mathbb Z)}
\]
and
\[
 \|f\|_{\widetilde L_T^\varrho(\cB^{s_1,s_2}_{p,r})}
 :=\left\|
 \left(2^{js_1+ks_2}
 \|\dot\Delta_j\dot\Delta_k^hf\|_{L_T^\varrho(L^p)}
 \right)_{j,k\in\mathbb Z}
 \right\|_{\ell^r(\mathbb Z^2)}.
\]
The corresponding spaces consist of time-dependent distributions
with finite norms, using the homogeneous representatives of
the respective spatial spaces.
\end{definition}

For $p=2$, $r=1$ and $1\leq\varrho\leq\infty$,
Minkowski's inequality implies
\[
\begin{aligned}
 \|f\|_{L_T^\varrho(\dot B^s_{2,1})}
 &\leq
 \|f\|_{\widetilde L_T^\varrho(\dot B^s_{2,1})},\qquad \|f\|_{L_T^\varrho(\cB^{s_1,s_2})}
 \leq
 \|f\|_{\widetilde L_T^\varrho(\cB^{s_1,s_2})}.
\end{aligned}
\]
Both inequalities are equalities when $\varrho=1$.

{
For the decay estimates, we recall the mixed space $\mathcal X^{s,r}$
introduced in \eqref{def:mixed-X}.
}
For $\mathcal X^{s,r}$, the time norm is taken on each mixed block
before the horizontal supremum. For $1\leq p\leq\infty$, we define
\begin{equation}\label{decay-mixed-time}
 \|f\|_{\widetilde L_T^p(\mathcal X^{s,r})}
 :=\sum_j2^{sj}\sup_k
       \bigl(2^{rk}\|\dot\Delta_j\dot\Delta_k^hf\|_{L_T^p(L^2)}\bigr).
\end{equation}

{{One of the main contributions of this work is the development of
new anisotropic product and commutator estimates, which play a
key role in our analysis and may be of independent interest.}}

\subsection{Anisotropic product estimates}\label{section:com}

{
\begin{proposition}[Anisotropic product estimates]\label{AnisoPL}\label{lem:low-product}
For $p=2$ and $r=1${ in dimension two}, the following estimates hold whenever the
right-hand sides are finite.

\smallskip\noindent\textup{(i)} \emph{(Positive radial index).}
Let $s_1,s_2>0$, and let $f,g\in\mathcal S'_h(\mathbb R^2)$ {have finite norms
on the right}. Then we have
\begin{equation}\label{app-positive-product}
 \begin{aligned}
 \|fg\|_{\cB^{s_1,s_2}}\lesssim{}&
 \|f\|_{\cB^{s_1,s_2}}\|g\|_{L^\infty}
 +\|g\|_{\cB^{s_1,s_2}}\|f\|_{L^\infty}\\
 &\quad+\|f\|_{\cB^{\frac{1}{2},s_2}}\|g\|_{\cB^{s_1,\frac{1}{2}}}
 +\|g\|_{\cB^{\frac{1}{2},s_2}}\|f\|_{\cB^{s_1,\frac{1}{2}}}.
 \end{aligned}
\end{equation}
In particular, the embedding $\cB^{\frac 12,\frac 12}\hookrightarrow L^\infty$
implies
\begin{equation}\label{app-critical-algebra}
 \|fg\|_{\cB^{\frac 12,\frac 12}}
 \lesssim\|f\|_{\cB^{\frac 12,\frac 12}}\|g\|_{\cB^{\frac 12,\frac 12}}.
\end{equation}

\smallskip\noindent\textup{(ii)} \emph{(Negative radial index).}
For $f\in\cB^{-\frac 12,\frac 12}$ and $g\in\dot B^1_{2,1}$, we have
\begin{equation}\label{dmltapp}
 \|fg\|_{\cB^{-\frac 12,\frac 12}}
 \lesssim\|f\|_{\cB^{-\frac 12,\frac 12}}\|g\|_{\dot B^1_{2,1}}.
\end{equation}
\end{proposition}

}

\begin{proposition}[Endpoint mixed product estimate]
\label{lem:endpoint-product}
{Let} {$s>0$ and} $f,g\in\mathcal S'_h(\mathbb R^2)$. We have
\begin{equation}\label{decay-mixed-product}
 \|fg\|_{\mathcal X^{{s},-\frac 12}}
 \lesssim\|f\|_{\cB^{\frac 12,\frac 12}}\|g\|_{\mathcal X^{{s},-\frac 12}}
        +\|g\|_{\cB^{\frac 12,\frac 12}}\|f\|_{\mathcal X^{{s},-\frac 12}}.
\end{equation}
\end{proposition}

We prove Proposition~\ref{AnisoPL} and
Proposition~\ref{lem:endpoint-product} in Appendix~\ref{app:product-proof}.
{The same estimates hold with time integrability exponents satisfying Hölder’s relation.}

\subsection{Commutator estimates}\label{section:commutator}

We first recall the classical commutator estimate in isotropic Besov spaces (see \cite[Chapter~2]{BCD}).

\begin{proposition}
\label{lem:besov-commutator}
Let $0\leq s\leq2$ and $f,g\in\mathcal S'_h(\mathbb R^2)$. We have
\[
 \sum_{j\in\mathbb Z}2^{js}\|[\dot\Delta_j,f]g\|_{L^2}
 \leq C\|\nabla f\|_{\dot B^1_{2,1}}
          \|g\|_{\dot B^{s-1}_{2,1}},
\]
where $C$ depends only on $s$ and the dyadic decomposition.
\end{proposition}

To establish anisotropic dissipation estimates, we need the following horizontal commutator bounds.

\begin{proposition}[Horizontal anisotropic commutator]\label{hcom}
Let \(f,g\in\mathcal S'_h(\mathbb R^2)\). We have
\begin{equation}\label{hc1}
\begin{aligned}
 \sum_{j,k\in\mathbb Z}2^{\frac32j+\frac12k}
 \bigl\|[f,\dot\Delta_j\dot\Delta_k^h]\partial_1g\bigr\|_{L^2}&\lesssim
 \|\partial_2f\|_{\cB^{\frac12,\frac12}}
 \|g\|_{\cB^{\frac12,\frac32}}
 +\|\partial_1f\|_{\cB^{\frac12,\frac12}}
 \|g\|_{\cB^{\frac32,\frac12}}.
\end{aligned}
\end{equation}
Moreover, we have
\begin{equation}\label{hc2}
\begin{aligned}
 \sum_{j,k\in\mathbb Z}2^{\frac12j+\frac12k}
 \bigl\|[f,\dot\Delta_j\dot\Delta_k^h]\partial_1g\bigr\|_{L^2}&\lesssim
 \|\partial_2f\|_{\dot B^1_{2,1}}
 \|\cR_1g\|_{\dot B^1_{2,1}}
 +\|\partial_1f\|_{\cB^{-\frac12,\frac12}}
 \|g\|_{\dot B^2_{2,1}}.
\end{aligned}
\end{equation}
\end{proposition}

\begin{proposition}[Endpoint mixed commutator estimate]
\label{lem:endpoint-commutator}
Let $s>0$ and $f,g\in\mathcal S'_h(\mathbb R^2)$. We have
\begin{equation}\label{decay-mixed-commutator}
 \begin{aligned}
 &\sum_j2^{{js}}\sup_k2^{-\frac{k}{2}}
       \|[f,\dot\Delta_j\dot\Delta_k^h]\partial_1g\|_{L^2}\\
 &\quad\qquad\lesssim\|f\|_{\cB^{\frac 32,\frac 12}}\|g\|_{\mathcal X^{{s-1},\frac 12}}
       +\|\partial_1f\|_{\cB^{\frac 12,\frac 12}}
                                   \|g\|_{\mathcal X^{{s},-\frac 12}}+\|f\|_{\mathcal X^{{s},-\frac 12}}
                                   \|g\|_{\cB^{\frac 12,\frac 32}}.
 \end{aligned}
\end{equation}
\end{proposition}

\begin{corollary}[Time decaying commutator estimate]
\label{lem:decay-damped-commutator}
Fix $c>0$. For $f,g\in\mathcal S'_h(\mathbb R^2)$ and every $\tau\geq0$, we have
\begin{equation}\label{decay-damped-commutator}
\begin{aligned}
 &\sum_{j,k}2^{-\frac{j}{2}+\frac{5k}{2}}e^{-c2^{2(k-j)}\tau}
       \|[f,\dot\Delta_j\dot\Delta_k^h]\partial_1g\|_{L^2}\\
 &\quad\qquad\leq\frac{C}{1+\tau}\Bigl(
       \|\partial_2f\|_{\cB^{\frac 12,\frac 12}}\|g\|_{\cB^{\frac 12,\frac 32}}
       +\|\partial_1f\|_{\cB^{\frac 12,\frac 12}}\|g\|_{\cB^{\frac 32,\frac 12}}
       \Bigr).
\end{aligned}
\end{equation}
The constant $C>0$ is independent of $\tau$.
\end{corollary}

{{The above product and commutator estimates also hold in the
corresponding time-dependent norms, with
$1\leq p,p_1,p_2\leq\infty$ satisfying $\frac{1}{p}=\frac{1}{p_1}+\frac{1}{p_2}$,
by applying H\"older's inequality in time to each dyadic product.}}

We prove Proposition~\ref{hcom}, Proposition~\ref{lem:endpoint-commutator} and Corollary~\ref{lem:decay-damped-commutator} in
Appendix~\ref{app:commutator-proof}.

\section{A priori estimates for the IPM equation}\label{dcsec}
\subsection{A priori estimates}
In this section, we derive the a priori estimates needed for the global stability of solutions to the stratified IPM equation. The dissipation
$-\cR_1^2$ is anisotropic: on a block localized at radial frequency $2^j$
and horizontal frequency $2^k$, it acts with strength $2^{2(k-j)}$ and anisotropic commutator estimates do not resolve this issue
(see Lemma \ref{lem:2.97}). 

Our argument has two parts. We first use isopycnal coordinates to eliminate
the vertical transport and exploit the anisotropic smoothing, then
combine the resulting bounds with isotropic Besov energy
estimates to propagate the norms
$\dot B^1_{2,1}\cap\dot B^2_{2,1}$.

The IPM system reads
\begin{equation}\label{dc1}
 \begin{cases}
 \partial_t b+u\cdot\nabla b-\cR_1^2b=0,\\
 u=(-\cR_2\cR_1b,\cR_1^2b),\\
 b|_{t=0}=b_0.
 \end{cases}
\end{equation}

\begin{proposition}\label{dc0}
Let $b$ be a smooth solution of \eqref{dc1} on $[0,T]$, with
$b_0\in\dot B^1_{2,1}\cap\dot B^2_{2,1}$. There exist constants
$\delta_0>0$ and $C>0$, independent of $T$, such that, if
\begin{equation}\label{dsmal}
 \|b_0\|_{\dot B^2_{2,1}}\leq\delta_0,
\end{equation}
then we have
\begin{equation}\label{db2}
 \|b\|_{L^\infty_T(\dot B^2_{2,1})}
 +\|\cR_1b\|_{L^2_T(\dot B^2_{2,1})}
 +\int_0^T\|\nabla u_2\|_{L^\infty}\,dt
 \leq C\|b_0\|_{\dot B^2_{2,1}},
\end{equation}
and
\begin{equation}\label{db1}
 \|b\|_{L^\infty_T(\dot B^1_{2,1})}
 +\|\cR_1b\|_{L^2_T(\dot B^1_{2,1})}
 +\int_0^T\|u_2\|_{L^\infty}\,dt
 \leq C\|b_0\|_{\dot B^1_{2,1}}.
\end{equation}
Moreover,
$b\in C([0,T];\dot B^1_{2,1}\cap\dot B^2_{2,1})$.
\end{proposition}

The proof is organized in six steps.

\subsubsection{Step 1. Isopycnal coordinates}

The purpose of the change of variables is to use the transported quantity
$x_2-b$ as the vertical coordinate. This is the content of the following lemma. 

\begin{lemma}[Isopycnal coordinates]\label{lem:density-reduction}
Let $b$ be a solution of \eqref{dc1} on $[0,T]$ such that
\[
 b\in C^1([0,T]\times\R^2),
 \qquad
 \partial_1b,\partial_2b\in L^\infty([0,T]\times\R^2),
\]
and assume that there exists $\kappa>0$ such that
\[
 \inf_{(t,x)\in[0,T]\times\R^2}
 \bigl(1-\partial_2b(t,x)\bigr)\geq\kappa.
\]
We define
\[
 Y(t,x){:=}(x_1,x_2-b(t,x)),\qquad X{:=}Y^{-1},
 \qquad h{:=}b\circ X,\qquad U=(U_1,U_2){:=}u\circ X.
\]
Then, for every $t\in[0,T]$, the map $Y(t,\cdot)$ is a global
$C^1$ bi-Lipschitz diffeomorphism of $\R^2$, and its inverse is given by
\[
 X(t,y)=(y_1,y_2+h(t,y)).
\]
Writing $ J{:=}1+\partial_2h$, we have
\begin{equation}\label{dc4}
\begin{aligned}
 DX&=
 \begin{pmatrix}
 1&0\\
 \partial_1h&J
 \end{pmatrix},
 \qquad
 (DY)\circ X=
 \begin{pmatrix}
 1&0\\
 -\frac{\partial_1h}{J}&\frac{1}{J}
 \end{pmatrix},\qquad J =\frac1{1-(\partial_2b)\circ X}.
\end{aligned}
\end{equation}
In particular, it follows that
\begin{equation}
    \label{eq:bdj}
 \frac{1}{1+\|\partial_2b\|_{L^\infty([0,T]\times\R^2)}}
 \leq J(t,y)\leq \frac1{\kappa}.
\end{equation}

Moreover, $h$ satisfies
\begin{equation}\label{dc6}
 \partial_th+U_1\partial_1h-\cR_1^2h=G(h),
\end{equation}
where
\begin{equation}\label{eq:defG}
 G(h){:=}\bigl[\cR_1^2(h\circ Y)\bigr]\circ X-\cR_1^2h
     =U_2-\cR_1^2h.
\end{equation}
Note that the vertical transport term has disappeared: the only
transport term in \eqref{dc6} is $U_1\partial_1h$.
\end{lemma}

\begin{remark}
The assumptions of Lemma~\ref{lem:density-reduction} are the
geometric conditions needed for the change to isopycnal coordinates.
Indeed, for fixed $t$ and $x_1$, the map
\[
 x_2\longmapsto x_2-b(t,x_1,x_2)
\]
is strictly increasing, since $1-\partial_2b\geq\kappa>0$. As the first coordinate is unchanged,
this yields the invertibility of $Y$. Then, the $L^\infty$ bounds on
$\partial_1b$ and $\partial_2b$, together with the lower bound on
$1-\partial_2b$, give the Lipschitz bounds for $Y$ and $X$.

In the applications below, these assumptions follow from the smallness
condition imposed on $b$. In particular, a sufficiently small
$L^\infty(0,T;\cB^{\frac 32,\frac 12})$ norm guarantees the required Lipschitz
control and ensures that $1-\partial_2b$ remains uniformly bounded away
from zero. The stronger isotropic regularity $\dot B^1_{2,1}\cap \dot B^2_{2,1}$ is not needed for the invertibility of the isopycnal change of variables.
It is used in the next lemma to control compositions with $X$ and $Y$, and hence to
transfer the estimates obtained for $(h,U)$ to the original
variables $(b,u)$.
\end{remark}

\begin{proof}[Proof of Lemma~\ref{lem:density-reduction}]
Since $u_2=\cR_1^2b$, equation \eqref{dc1} gives
\[
(\partial_t+u\cdot\nabla)b=u_2,
\qquad
(\partial_t+u\cdot\nabla)(x_2-b)=0.
\]
Thus $x_2-b$ is transported by the flow. Fix $t\in[0,T]$ and $x_1\in\mathbb R$, and set
\[
\psi_{t,x_1}(x_2):=x_2-b(t,x_1,x_2).
\]
For $x_2>z_2$, we have
\[
\psi_{t,x_1}(x_2)-\psi_{t,x_1}(z_2)
=
\int_{z_2}^{x_2}
\bigl(1-\partial_2b(t,x_1,s)\bigr)\,ds
\geq \kappa(x_2-z_2).
\]
Hence $\psi_{t,x_1}$ is strictly increasing and tends to $\pm\infty$
as $x_2\to\pm\infty$. Since the first coordinate is unchanged,
$Y(t,x)=(x_1,x_2-b(t,x))$ is a global bijection of $\mathbb R^2$.

The boundedness of $\nabla b$ implies that $Y(t,\cdot)$ is Lipschitz.
Let $y=Y(t,x)$ and $y'=Y(t,x')$. Since $x_1-x_1'=y_1-y_1'$, we get
\[
\begin{aligned}
|x_2-x_2'|
&\leq
\frac1\kappa
\left|
\psi_{t,x_1}(x_2)-\psi_{t,x_1}(x_2')
\right|\\
&\leq
\frac1\kappa
\left(
|y_2-y_2'|
+
\left|
\psi_{t,x_1'}(x_2')-\psi_{t,x_1}(x_2')
\right|
\right)\\
&\leq
\frac1\kappa
\left(
|y_2-y_2'|
+
\|\partial_1b(t)\|_{L^\infty}|y_1-y_1'|
\right).
\end{aligned}
\]
Therefore $X=Y^{-1}$ is Lipschitz as well. Since
$\det DY=1-\partial_2b\geq\kappa$, the implicit function theorem
also shows that $X$ is $C^1$ in $(t,y)$.

By definition, one has
\[
X(t,y)=(y_1,y_2+h(t,y)).
\]
Writing $J=1+\partial_2h$ and differentiating $Y\circ X=\operatorname{Id}$,
we obtain
\[
DX=
\begin{pmatrix}
1&0\\
\partial_1h&J
\end{pmatrix},
\qquad
(DY)\circ X=(DX)^{-1}
=
\begin{pmatrix}
1&0\\
-\frac{\partial_1h}{J}&\frac{1}{J}
\end{pmatrix}.
\]
In particular, it follows that
\[
J=\frac1{1-(\partial_2b)\circ X},
\qquad
\frac1{1+\|\partial_2b\|_{L^\infty([0,T]\times\mathbb R^2)}}
\leq J\leq\frac1\kappa.
\]

We now derive the equation satisfied by $h$. Differentiating
\[
h(t,y)=b\bigl(t,y_1,y_2+h(t,y)\bigr)
\]
with respect to $t$ leads to
\[
\partial_th
=
(\partial_tb)\circ X
+
\bigl((\partial_2b)\circ X\bigr)\partial_th.
\]
Moreover, one has
\[
(\partial_tb)\circ X=\frac{\partial_th}{J},
\qquad
(\partial_1b)\circ X=\frac{\partial_1h}{J},
\qquad
(\partial_2b)\circ X=\frac{\partial_2h}{J}.
\]
Composing $\partial_tb+u_1\partial_1b+u_2\partial_2b=u_2$ with $X$,
we obtain
\[
\partial_th+U_1\partial_1h+U_2\partial_2h=JU_2.
\]
Since $J=1+\partial_2h$, this reduces to
\[
\partial_th+U_1\partial_1h=U_2.
\]
Finally, $b=h\circ Y$ and $u_2=\cR_1^2b$ imply
\[
U_2=\bigl[\cR_1^2(h\circ Y)\bigr]\circ X.
\]
This proves \eqref{dc6} and \eqref{eq:defG}.
\end{proof}

We now provide composition estimates needed to transfer isotropic norms
between the Eulerian and isopycnal variables. This lemma will be a key ingredient in the following.

\begin{lemma}[From Eulerian to isopycnal coordinates and back]
\label{dcmp}
Assume that $b\in\dot B^1_{2,1}\cap\dot B^2_{2,1}$ and that
$\|b\|_{\dot B^2_{2,1}}$ is sufficiently small. Define $X,Y,h$ as above.
Then, for each $s\in\{0,1,2\}$ and every
$f\in\dot B^s_{2,1}$, it holds that
\begin{equation}\label{dcn}
 \|f\circ X\|_{\dot B^s_{2,1}}
 +\|f\circ Y\|_{\dot B^s_{2,1}}
 \leq C\|f\|_{\dot B^s_{2,1}}.
\end{equation}
In particular, taking $f=b$ and $f=h$, respectively, we have
\[
 \|h\|_{\dot B^s_{2,1}}\simeq\|b\|_{\dot B^s_{2,1}},
 \qquad s=1,2,
\]
and
\begin{equation}\label{dc3}
 \|h(0)\|_{\cB^{\frac32,\frac12}}
 \leq C\|b_0\|_{\dot B^2_{2,1}},\quad
 \|h(0)\|_{\dot B^s_{2,1}}
 \leq C\|b_0\|_{\dot B^s_{2,1}},
 \qquad s=1,2.
\end{equation}
Conversely, if $h\in\dot B^1_{2,1}\cap\dot B^2_{2,1}$ and
$\|h\|_{\dot B^2_{2,1}}$ is sufficiently small, we define
$X(y){:=}(y_1,y_2+h(y))$, $Y{:=}X^{-1}$ and $b{:=}h\circ Y$.
Then \eqref{dcn} remains valid. 
\end{lemma}

 The proof of Lemma~\ref{dcmp} is given in Appendix~\ref{app:dc-composition}.
\begin{remark}
Once $X$ and $Y$ are known to be bi-Lipschitz, one can use the change of variables in one of the two directions: $h=b\circ X$ and $b=h\circ Y$.  The role of the isotropic regularity $\dot B^1_{2,1}\cap\dot B^2_{2,1}$ is to control the above compositions in our functional spaces. More precisely, Lemma~\ref{dcmp} ensures that the maps $a\longmapsto a\circ X$ and $a\longmapsto a\circ Y$ are bounded on the isotropic spaces $\dot B^s_{2,1}$, for
$s=0,1,2$. This allows us to transfer quantitative estimates obtained
for $h$ in the new coordinates to corresponding estimates for $b$ in
the original coordinates.
Note that this point is less immediate in the anisotropic space
$\cB^{\frac 32,\frac 12}$ because the change of variables mixes the horizontal
and vertical directions, e.g.
\[
\partial_1(h\circ Y)
=
(\partial_1h)\circ Y
-
\partial_1b\,(\partial_2h)\circ Y.
\]
\end{remark}

\subsubsection{Step 2. Anisotropic estimates}
\begin{lemma}\label{dctechnical}
Assume that $\|h\|_{\cB^{\frac32,\frac12}}$ is sufficiently small. Then we have
\begin{equation}\label{dcb}
 \|G(h)\|_{\cB^{\frac32,\frac12}}
 \leq C\left(
 \|h\|_{\cB^{\frac12,\frac32}}^2
 +\|h\|_{\cB^{\frac32,\frac12}}
       \|h\|_{\cB^{-\frac12,\frac52}}
 \right).
\end{equation}
Moreover, we have
\begin{align}
 \|\partial_2U_1\|_{\cB^{\frac12,\frac12}}
 &\leq C\|h\|_{\cB^{\frac12,\frac32}},\label{dch}\\
 \|\partial_1U_1\|_{\cB^{\frac12,\frac12}}
 &\leq C\|h\|_{\cB^{-\frac12,\frac52}},\label{dci}\\
 \|\nabla U_2\|_{L^\infty}
 &\leq C\left(
 \|h\|_{\cB^{-\frac12,\frac52}}
 +\|h\|_{\cB^{\frac12,\frac32}}^2
 +\|h\|_{\cB^{\frac32,\frac12}}
       \|h\|_{\cB^{-\frac12,\frac52}}
 \right).\label{dcu2lip}
\end{align}
\end{lemma}

The proof is postponed to Appendix~\ref{app:dc-riesz}.

\begin{lemma}\label{dchm}
Let $h$ solve \eqref{dc6} on $[0,T]$ and assume
\[
 \|h\|_{L^\infty_T(\dot B^2_{2,1})}\leq c_0.
\]
If $c_0$ is sufficiently small, then we have
\begin{equation}\label{dc2}
 \begin{split}
 &\|h\|_{\widetilde L^\infty_T(\cB^{\frac32,\frac12})}
 +\|h\|_{\widetilde L^2_T(\cB^{\frac12,\frac32})}
 +\|h\|_{L^1_T(\cB^{-\frac12,\frac52})}
 +\int_0^T\|\nabla u_2\|_{L^\infty}\,dt
 \leq C\|h(0)\|_{\cB^{\frac32,\frac12}},
 \end{split}
\end{equation}
where $C$ is independent of $T$.
\end{lemma}

\begin{proof}
On the support of $\dot\Delta_j\dot\Delta_k^h h$, we have
\[
 -\langle\cR_1^2\dot\Delta_j\dot\Delta_k^h h,
 \dot\Delta_j\dot\Delta_k^h h\rangle
 \geq c\,2^{2(k-j)}\|\dot\Delta_j\dot\Delta_k^h h\|_{L^2}^2.
\]
We apply $\dot\Delta_j\dot\Delta_k^h$ to \eqref{dc6}, take the
$L^2$ scalar product with $\dot\Delta_j\dot\Delta_k^h h$ to obtain, for almost every $t$,
\begin{equation}\label{dck}
\begin{aligned}
 &\frac{d}{dt}\|\dot\Delta_j\dot\Delta_k^h h\|_{L^2}
 +c2^{2(k-j)}\|\dot\Delta_j\dot\Delta_k^h h\|_{L^2}\\
 &\quad\leq
 \|\dot\Delta_j\dot\Delta_k^h G(h)\|_{L^2}
 +\|[U_1,\dot\Delta_j\dot\Delta_k^h]\partial_1h\|_{L^2}
 +C\|\partial_1U_1\|_{L^\infty}
 \|\dot\Delta_j\dot\Delta_k^h h\|_{L^2}.
\end{aligned}
\end{equation}
We integrate \eqref{dck} over $[0,t]$, take the supremum over
$0\leq t\leq T$ and sum with weight
$2^{\frac32j+\frac12k}$. By \eqref{dcb}--\eqref{dci} and Proposition \ref{hcom}, we obtain
\begin{equation}\label{dcl}
 \begin{split}
 &\|h\|_{\widetilde L^\infty_T(\cB^{\frac32,\frac12})}
 +c\|h\|_{L^1_T(\cB^{-\frac12,\frac52})}\\
 &\quad \leq
 C\|h(0)\|_{\cB^{\frac32,\frac12}}
 +C\int_0^T\left(
 \|h\|_{\cB^{\frac12,\frac32}}^2
 +\|h\|_{\cB^{\frac32,\frac12}}
       \|h\|_{\cB^{-\frac12,\frac52}} \right)\,dt\\
 & \quad \lesssim \|h(0)\|_{\cB^{\frac32,\frac12}}+\|h\|_{L^\infty_T(\cB^{\frac32,\frac12})}
   \|h\|_{L^1_T(\cB^{-\frac12,\frac52})}.
 \end{split}
\end{equation}

\noindent Here we used the interpolation of
Proposition~\ref{prop:inter} with $\theta=\frac{1}{2}$ yielding 
$\|h\|_{\cB^{\frac12,\frac32}}^2
\lesssim\|h\|_{\cB^{\frac32,\frac12}}
\|h\|_{\cB^{-\frac12,\frac52}}$.
By the assumption of the lemma, we have $ \|h\|_{L^\infty_T(\cB^{\frac32,\frac12})}
 \lesssim\|h\|_{L^\infty_T(\dot B^2_{2,1})}\leq c_0$. Choosing $c_0$ sufficiently small, we absorb the last term of
\eqref{dcl} into its left-hand side.
Arguing as in the proof of Proposition~\ref{prop:inter},
with H\"older's inequality in time, we obtain
\[
 \|h\|_{\widetilde L^2_T(\cB^{\frac12,\frac32})}
 \leq
 \|h\|_{\widetilde L^\infty_T(\cB^{\frac32,\frac12})}^{\frac12}
 \|h\|_{L^1_T(\cB^{-\frac12,\frac52})}^{\frac12}
 \lesssim\|h(0)\|_{\cB^{\frac32,\frac12}}.
\]

Finally, it follows from \eqref{dcu2lip} that
\[
 \int_0^T\|\nabla_yU_2\|_{L^\infty}\,dt
 \leq C\|h(0)\|_{\cB^{\frac32,\frac12}}.
\]
Since $u_2=U_2\circ Y$ and
$\nabla_xu_2=(DY)^T(\nabla_yU_2\circ Y)$, the boundedness of $DY$ yields the
last term in \eqref{dc2}.
\end{proof}

\subsubsection{Step 3. Isotropic estimates}

We now propagate the two isotropic norms $\dot{B}^2_{2,1}$ and $\dot{B}^{1}_{2,1}$. The key is to
estimate the transformed transport term using the horizontal
commutator and the anisotropic bounds.
Using $G=\cR_1G_1+G_2$ { from \eqref{dc-G-decomposition}}, where $G_1$ and $G_2$ are defined in
\eqref{def:g1g2} and $\cR_1^*=-\cR_1$, we apply $\dot\Delta_j$ to
\eqref{dc6} and obtain
\begin{equation}\label{disbk}
 \begin{split}
 \frac12\|\dot\Delta_jh(t)\|_{L^2}^2
 +\int_0^t\|\cR_1\dot\Delta_jh\|_{L^2}^2\,d\tau
 ={}&\frac12\|\dot\Delta_jh(0)\|_{L^2}^2
 +\frac12\int_0^t\int\partial_1U_1|\dot\Delta_jh|^2\,dy\,d\tau\\
 &-\int_0^t
 \langle[\dot\Delta_j,U_1]\partial_1h,\dot\Delta_jh\rangle\,d\tau\\
 &-\int_0^t
 \langle\dot\Delta_jG_1,\cR_1\dot\Delta_jh\rangle\,d\tau
 +\int_0^t
 \langle\dot\Delta_jG_2,\dot\Delta_jh\rangle\,d\tau.
 \end{split}
\end{equation}
For $s=1,2$, Proposition~\ref{lem:besov-commutator} gives rise to
\begin{equation*}
 \sum_j2^{sj}
 \|[\dot\Delta_j,f]\partial_1g\|_{L^2}
 \leq C\|\nabla f\|_{\dot B^1_{2,1}}
          \|\cR_1g\|_{\dot B^s_{2,1}}.
\end{equation*}
Then, we take the supremum in time in \eqref{disbk}, 
take square roots and sum with weight $2^{js}$ to obtain
\begin{equation}\label{dien}
 \begin{split}
 \|h\|_{\widetilde L^\infty_T(\dot B^s_{2,1})}^2
 +\|\cR_1h\|_{\widetilde L^2_T(\dot B^s_{2,1})}^2
 &\leq C\|h(0)\|_{\dot B^s_{2,1}}^2
 +C\|\cR_1h\|_{\widetilde L^2_T(\dot B^s_{2,1})}
       \|G_1\|_{\widetilde L^2_T(\dot B^s_{2,1})}\\
 &\quad
 +C\|h\|_{\widetilde L^\infty_T(\dot B^s_{2,1})}
       \|G_2\|_{L^1_T(\dot B^s_{2,1})}\\
 &\quad
 +C\|\partial_1U_1\|_{L^1_T(L^\infty)}
       \|h\|_{\widetilde L^\infty_T(\dot B^s_{2,1})}^2\\
 &\quad
 +C\|h\|_{\widetilde L^\infty_T(\dot B^s_{2,1})}
       \|\nabla U_1\|_{L^2_T(\dot B^1_{2,1})}
       \|\cR_1h\|_{L^2_T(\dot B^s_{2,1})}.
 \end{split}
\end{equation}

For the terms on the right-hand side of \eqref{dien}, we use the
following estimates, which we prove in Appendix~\ref{app:dc-isotropic}.
\begin{lemma}\label{diso}
Assume that $\|h\|_{\cB^{\frac32,\frac12}}$ is sufficiently small. Then, for $s=1,2$, we have
\begin{equation}\label{dgt}
 \begin{aligned}
 \|G_1(h)\|_{\dot B^s_{2,1}}
 &\leq C\|h\|_{\cB^{\frac32,\frac12}}
          \|\cR_1h\|_{\dot B^s_{2,1}}{
       +C\|h\|_{\cB^{\frac32,\frac12}}
          \|h\|_{\dot B^s_{2,1}}
          \|h\|_{\cB^{\frac12,\frac32}}},\\
 \|G_2(h)\|_{\dot B^s_{2,1}}
 &\leq C\|h\|_{\cB^{\frac12,\frac32}}
          \|\cR_1h\|_{\dot B^s_{2,1}}
       +C\|h\|_{\dot B^s_{2,1}}
          \|h\|_{\cB^{-\frac12,\frac52}},\\
 \|U\|_{\dot B^s_{2,1}}
 &\leq C\|\cR_1h\|_{\dot B^s_{2,1}}
       +C\|h\|_{\dot B^s_{2,1}}
          \|h\|_{\cB^{\frac12,\frac32}}.
 \end{aligned}
\end{equation}
In particular, for $s=1,2$, we have
\begin{equation}\label{dm12}
 \|U_1\|_{\dot B^s_{2,1}}
 \leq C\left(
 \|\cR_1h\|_{\dot B^s_{2,1}}
 +\|h\|_{\dot B^s_{2,1}}
  \|h\|_{\cB^{\frac12,\frac32}}
 \right),
\end{equation}
and
\begin{equation}\label{dm13}
 \|G(h)\|_{\dot B^s_{2,1}}
 \leq C\left(
 \|h\|_{\cB^{\frac32,\frac12}}
       \|\cR_1h\|_{\dot B^s_{2,1}}
 +\|h\|_{\dot B^s_{2,1}}
       \|h\|_{\cB^{\frac12,\frac32}}
 \right).
\end{equation}
Moreover, we have
\[
 \|G_1(h)\|_{\widetilde L^2_T(\dot B^s_{2,1})}
 \leq C
 \|h\|_{\widetilde L^\infty_T(\cB^{\frac32,\frac12})}
 \|\cR_1h\|_{\widetilde L^2_T(\dot B^s_{2,1})}
 {+C
 \|h\|_{\widetilde L^\infty_T(\cB^{\frac32,\frac12})}
 \|h\|_{\widetilde L^\infty_T(\dot B^s_{2,1})}
 \|h\|_{\widetilde L^2_T(\cB^{\frac12,\frac32})}},
\]
and
\[
 \begin{split}
 \|G_2(h)\|_{L^1_T(\dot B^s_{2,1})}
 \leq C\Bigl(
 &\|h\|_{\widetilde L^2_T(\cB^{\frac12,\frac32})}
  \|\cR_1h\|_{\widetilde L^2_T(\dot B^s_{2,1})}+\|h\|_{\widetilde L^\infty_T(\dot B^s_{2,1})}
  \|h\|_{L^1_T(\cB^{-\frac12,\frac52})}
 \Bigr).
 \end{split}
\]
\end{lemma}

\subsubsection{Step 4. Closure of the isotropic estimates}\label{diclo}

To close the estimates, we combine the isotropic and anisotropic inequalities and absorb
their nonlinear terms on a small bootstrap interval. We first work on
an interval where the isopycnal coordinates are defined, with
$\sup_{0\leq t\leq T}\|\nabla b(t)\|_{L^\infty}<\frac{1}{2}$, and where
\begin{equation}\label{dbs}
 \|h\|_{\widetilde L^\infty_T(\dot B^2_{2,1})}\leq c_*,
\end{equation}
where $c_*>0$ is fixed and small. These two conditions hold for a
short time by continuity and \eqref{dsmal}. In Step~6, we improve both of
them. Recall the elementary embeddings { in \eqref{eq:embh}}
\begin{equation}\label{dca}
 \|f\|_{\cB^{\frac32,\frac12}}
 \leq C\|f\|_{\dot B^2_{2,1}},
 \qquad
 \|f\|_{\cB^{\frac12,\frac12}}
 \leq C\|f\|_{\dot B^1_{2,1}}.
\end{equation}
The embeddings \eqref{dca}, together with Lemmas~\ref{dcmp} and \ref{dchm}, imply
\begin{equation}\label{dahi}
 \begin{split}
 &\|h\|_{\widetilde L^\infty_T(\cB^{\frac32,\frac12})}
 +\|h\|_{\widetilde L^2_T(\cB^{\frac12,\frac32})}
 +\|h\|_{L^1_T(\cB^{-\frac12,\frac52})}
 +\|\partial_1U_1\|_{L^1_T(L^\infty)}\leq C\|b_0\|_{\dot B^2_{2,1}}.
 \end{split}
\end{equation}

{
For $s=1,2$, we use \eqref{dahi} to bound the contribution of the
additional term in the estimate of $G_1$ to \eqref{dien} by
\[
 \begin{aligned}
 &C\|h\|_{\widetilde L^\infty_T(\cB^{\frac32,\frac12})}
   \|h\|_{\widetilde L^2_T(\cB^{\frac12,\frac32})}
   \|h\|_{\widetilde L^\infty_T(\dot B^s_{2,1})}
   \|\cR_1h\|_{\widetilde L^2_T(\dot B^s_{2,1})}\leq C\delta_0^2\left(
   \|h\|_{\widetilde L^\infty_T(\dot B^s_{2,1})}^2
   +\|\cR_1h\|_{\widetilde L^2_T(\dot B^s_{2,1})}^2
   \right).
 \end{aligned}
\]
For sufficiently small $\delta_0$, we absorb this term into the
left-hand side of \eqref{dien}, uniformly in $T$.
}

We first take $s=2$ in \eqref{dien}. Applying Lemma~\ref{diso},
\eqref{dahi} and Young's inequality, we have
\begin{equation}\label{di2}
 \begin{split}
 &\|h\|_{\widetilde L^\infty_T(\dot B^2_{2,1})}^2
 +\|\cR_1h\|_{\widetilde L^2_T(\dot B^2_{2,1})}^2\leq C\|b_0\|_{\dot B^2_{2,1}}^2
 +C(c_*+\delta_0)\left(
 \|h\|_{\widetilde L^\infty_T(\dot B^2_{2,1})}^2
 +\|\cR_1h\|_{\widetilde L^2_T(\dot B^2_{2,1})}^2
 \right).
 \end{split}
\end{equation}
One may choose $c_*$ and $\delta_0$ so that the last term is
absorbed. Taking $\delta_0$ smaller if necessary, we obtain
\[
 \|h\|_{\widetilde L^\infty_T(\dot B^2_{2,1})}
 \leq C\|b_0\|_{\dot B^2_{2,1}}
 \leq C\delta_0\leq \frac{c_*}{2},
\]
which improves \eqref{dbs} by a bootstrap argument.

Using \eqref{dahi} and
$\|\nabla U_1\|_{L^2_T(\dot B^1_{2,1})}\lesssim \delta_0 \ll 1$ 
in \eqref{dien}, we obtain, for $s=1,2$,
\begin{equation}\label{di12}
 \|h\|_{\widetilde L^\infty_T(\dot B^s_{2,1})}
 +\|\cR_1h\|_{\widetilde L^2_T(\dot B^s_{2,1})}
 \leq C\|b_0\|_{\dot B^s_{2,1}}.
\end{equation}

\subsubsection{Step 5. Time integrability of the vertical velocity}

We recover the low-order time integral of $u_2$ needed in
\eqref{db1}. To this end, we need the following lemma.

{{\begin{lemma}\label{lem:dc-loworder}
If $\|h\|_{\dot B^2_{2,1}}$ is sufficiently small, then we have
\begin{equation}\label{dglo}
 \|G(h)\|_{\cB^{\frac12,\frac12}}
 \leq C\|\cR_1h\|_{\dot B^1_{2,1}}
          \|\cR_1h\|_{\dot B^2_{2,1}}
 +C\|h\|_{\dot B^2_{2,1}}
          \|h\|_{\cB^{-\frac32,\frac52}}.
\end{equation}
\end{lemma}
The proof of Lemma \ref{lem:dc-loworder} is given in Appendix~\ref{app:dc-loworder}. Integrating \eqref{dck} over $[0,T]$ and summing with the lower weight
$2^{\frac12j+\frac12k}$, we estimate the source term $G(h)$ by
\eqref{dglo} and the commutator by \eqref{hc2}, \eqref{dm12} with
$s=2$ and \eqref{dc-low-gradient}. Together with \eqref{dbs} and
the identity
\[
 J\partial_1U_1=\partial_1h\,\partial_2U_1-\partial_2U_2,
\]
this yields
\begin{equation}\label{dllo}
 \begin{split}
 \|h\|_{L^1_T(\cB^{-\frac32,\frac52})}
 \leq{}&C\|h(0)\|_{\dot B^1_{2,1}}
 +C\|h\|_{\widetilde L^\infty_T(\dot B^2_{2,1})}
       \|h\|_{L^1_T(\cB^{-\frac32,\frac52})}\\
 &+C\|\cR_1h\|_{\widetilde L^2_T(\dot B^1_{2,1})}
 \left(
 \|\cR_1h\|_{\widetilde L^2_T(\dot B^2_{2,1})}
 +c_*\|h\|_{\widetilde L^2_T(\cB^{\frac12,\frac32})}
 \right)\\
 &+C\|\partial_1U_1\|_{L^1_T(L^\infty)}
       \|h\|_{\widetilde L^\infty_T(\dot B^1_{2,1})}\\
 \leq{}&C\|b_0\|_{\dot B^1_{2,1}}
 +Cc_*\|h\|_{L^1_T(\cB^{-\frac32,\frac52})}
 +C\|b_0\|_{\dot B^1_{2,1}}\|b_0\|_{\dot B^2_{2,1}}.
 \end{split}
\end{equation}
The $L^2_T$ norms are bounded by \eqref{dahi} and \eqref{di12},
which also control the last term. Together with \eqref{dcn} at $t=0$ and
\eqref{dbs}, this proves the second inequality.
Since $c_*$ in \eqref{dbs} is chosen so that $Cc_*\leq\frac{1}{2}$,
it follows from \eqref{dsmal} that
\[
 \|h\|_{L^1_T(\cB^{-\frac32,\frac52})}
 \leq \frac{C}{1-Cc_*}
 \bigl(1+\delta_0\bigr)
 \|b_0\|_{\dot B^1_{2,1}}
 \leq C\|b_0\|_{\dot B^1_{2,1}}.
\]
Finally, due to $U_2=\cR_1^2h+G(h)$, 
$\cB^{\frac12,\frac12}\hookrightarrow L^\infty$, 
\eqref{dbs}, \eqref{di12} and \eqref{dglo}, one has
\begin{equation}\label{du0}
 \begin{aligned}
 \int_0^T\|U_2\|_{L^\infty}\,dt
 &\leq C\|h\|_{L^1_T(\cB^{-\frac32,\frac52})}
       +C\|G(h)\|_{L^1_T(\cB^{\frac12,\frac12})}\\
 &\leq C(1+c_*)\|h\|_{L^1_T(\cB^{-\frac32,\frac52})}
       +C\|\cR_1h\|_{\widetilde L^2_T(\dot B^1_{2,1})}
          \|\cR_1h\|_{\widetilde L^2_T(\dot B^2_{2,1})}\\
 &\leq C\|b_0\|_{\dot B^1_{2,1}}.
 \end{aligned}
\end{equation}}}

\subsubsection{Step 6. Returning to the original variables and time continuity}
\label{dcret}

We transfer the estimates for $h,U$ back to $b,u$ through the coordinate
change.
Since $b=h\circ Y$, Lemma~\ref{dcmp} and \eqref{di12} lead to the
$L^\infty_T(\dot B^s_{2,1})$ estimates for $b$.  Moreover, \eqref{dc4} leads to
\begin{equation}
    \label{rel:d1bh}
    (\partial_{x_1}b)\circ X=\frac{\partial_1h}{J}.
\end{equation}
Since $ \|J^{-1}-1\|_{\dot B^1_{2,1}}
 \lesssim \|h\|_{\dot B^2_{2,1}}$, we have
\[
\|J^{-1}\partial_1h\|_{\dot B^s_{2,1}}
\lesssim
(1+\|h\|_{\dot B^2_{2,1}})
\|\partial_1h\|_{\dot B^s_{2,1}}
\lesssim
\|\partial_1h\|_{\dot B^s_{2,1}},
\qquad s=0,1.
\]
Using this bound together with \eqref{eq:bdj}{, \eqref{dcn}},
\eqref{rel:d1bh}, Bernstein's inequality and the
product law in $\dot B^1_{2,1}$, one obtains, for $s=1,2$,
\[
 \|\cR_1b\|_{\dot B^s_{2,1}}
 \lesssim
 \|\partial_{x_1}b\|_{\dot B^{s-1}_{2,1}}
 \leq
 C\|\partial_1h\|_{\dot B^{s-1}_{2,1}}
 =
 C\|\Lambda\Lambda^{-1}\partial_1 h\|_{\dot B^{s-1}_{2,1}}
 \lesssim
 \|\cR_1h\|_{\dot B^s_{2,1}}.
\]
{{Combining the preceding estimate with \eqref{di12} proves the
dissipation bounds in \eqref{db2}--\eqref{db1}.
The bounds for $\|\nabla u_2\|_{L^1_T(L^\infty)}$ and
$\|u_2\|_{L^1_T(L^\infty)}$ follow from
Lemma~\ref{dchm} and \eqref{du0}, respectively.}}

For completeness, we briefly justify time continuity. {{Equation
\eqref{dc6}, together with \eqref{dm12}--\eqref{dm13},
\eqref{dahi} and \eqref{di12},}} implies
$\partial_th\in L^1_T(\dot B^1_{2,1})$ on every finite interval.
Hence, we have $ h\in C([0,T];\dot B^1_{2,1})$. For the higher-order continuity, we split the frequencies to obtain,
for every $J_*\geq1$,
\[
 \begin{split}
 \|h(t)-h(t')\|_{\dot B^2_{2,1}}
 &\leq
 2\sum_{|j|\geq J_*+1}2^{2j}
 \|\dot\Delta_jh\|_{L^\infty_T(L^2)}
 +2^{J_*}\|h(t)-h(t')\|_{\dot B^1_{2,1}}.
 \end{split}
\]
The first term is uniformly small for $J_*$ large, and the second tends to
zero as $t\to t'$. Thus, we have $ h\in C([0,T];\dot B^1_{2,1}\cap\dot B^2_{2,1})$. The uniform composition bounds and the continuity of the inverse map
give the same continuity for $b=h\circ Y$. Moreover, it follows from \eqref{di12} that
\[
 \sup_{0\leq t\leq T}\|\nabla b(t)\|_{L^\infty}
 \leq
 \frac{\|\nabla h\|_{L^\infty_T(L^\infty)}}
      {1-\|\partial_2h\|_{L^\infty_T(L^\infty)}}
 \leq C\|h\|_{\widetilde L^\infty_T(\dot B^2_{2,1})}.
\]
After decreasing $\delta_0$ if necessary, the last quantity is at most
$\frac{1}{4}$. Together with the improvement of \eqref{dbs}, this closes both
bootstrap conditions and completes the proof of Proposition~\ref{dc0}. \qed

\section{Well-posedness of stratified IPM}\label{IPM}

In this section, we prove Theorem \ref{them1:local}, Theorem \ref{them1} and Theorem \ref{thm:ipm-higher-sobolev}.

\subsection{Proof of Theorem~\ref{them1:local}} 
 We first prove Theorem~\ref{them1:local} regarding the local well-posedness of \eqref{eq:ipm3}. Local existence follows
from a standard approximation argument in Besov spaces
{{(see \cite[Chapters~3 and~7]{BCD}), once the uniform estimates below
are established}}. Applying $\dot\Delta_j$ to \eqref{eq:ipm3},  taking the
$L^2$ scalar product with $\dot\Delta_jb$ and using $\operatorname{div}u=0$, we obtain
\begin{equation*}
 \begin{aligned}
 \frac12\frac{d}{dt}\|\dot\Delta_jb\|_{L^2}^2
 +\|\mathcal R_1\dot\Delta_jb\|_{L^2}^2
 \leq
 \|[u\cdot\nabla,\dot\Delta_j]b\|_{L^2}
 \|\dot\Delta_jb\|_{L^2}.
 \end{aligned}
\end{equation*}
For $s=1,2$, we apply Proposition~\ref{lem:besov-commutator}
with $f=u_k$ and $g=\partial_kb$, and sum over $k=1,2$.
Using $u=(-\mathcal R_2\mathcal R_1b,\mathcal R_1^2b)$, we obtain
\begin{equation*}
 \begin{aligned}
 \sum_{j\in\mathbb Z}2^{sj}
 \|[u\cdot\nabla,\dot\Delta_j]b\|_{L^2}
 &\leq C\|u\|_{\dot B^2_{2,1}}\|b\|_{\dot B^s_{2,1}}\leq C\|b\|_{\dot B^2_{2,1}}\|b\|_{\dot B^s_{2,1}}.
 \end{aligned}
\end{equation*}
We discard the nonnegative dissipative term in the block inequality
and divide by a regularized $L^2$ norm. Integrating in time, taking
the time supremum and summing with weights $2^j$ and $2^{2j}$, we obtain
\begin{equation*}
 \begin{aligned}
 \|b\|_{\widetilde L^\infty_T
       (\dot B^1_{2,1}\cap\dot B^2_{2,1})}
 &\leq C\|b_0\|_{\dot B^1_{2,1}\cap\dot B^2_{2,1}}+C\int_0^T
 \|b(t)\|_{\dot B^2_{2,1}}
 \|b(t)\|_{\dot B^1_{2,1}\cap\dot B^2_{2,1}}\,dt.
 \end{aligned}
\end{equation*}
By a continuity argument, we bound the approximate solutions
uniformly on a time interval depending only on
$\|b_0\|_{\dot B^1_{2,1}\cap\dot B^2_{2,1}}$. Using the usual compactness
argument, we then obtain a solution $ b\in\widetilde L^\infty_T
       (\dot B^1_{2,1}\cap\dot B^2_{2,1})$. {{Using \eqref{eq:ipm3} and arguing as in Section~\ref{dcret},
we obtain $b\in C([0,T];\dot B^1_{2,1}\cap\dot B^2_{2,1})$.}}

We prove uniqueness by estimating the difference one derivative below
the critical norm, as in the Besov approach to fluid equations described
in \cite{DanchinFourier2018}. Let $(b^1,u^1)$ and $(b^2,u^2)$ be two
solutions of \eqref{eq:ipm3} with the same initial datum. We define
\[
 \delta b:=b^1-b^2,\qquad \delta u:=u^1-u^2.
\]
The difference $(\delta b,\delta u)$ satisfies
\begin{equation}\label{eq:ipm3-error}
 \begin{cases}
 \partial_t\delta b+u^2\cdot\nabla\delta b-\cR_1^2\delta b
       =-\delta u\cdot\nabla b^1,\\
 \delta b|_{t=0}=0,
 \end{cases}
 \qquad
 \delta u=\cR_1\cR^\perp\delta b.
\end{equation}
{
We apply $\dot\Delta_j$ to \eqref{eq:ipm3-error} and take the
$L^2$ scalar product with $\dot\Delta_j\delta b$. Since
$\operatorname{div}u^2=0$, the classical isotropic commutator estimate
of Proposition~\ref{lem:besov-commutator}, at index $s=1$, gives,
after summation with weight $2^j$,
}
\begin{equation}\label{ipm-uniqueness-est}
 \begin{split}
 &\|\delta b\|_{\widetilde L^\infty_T(\dot B^1_{2,1})}
 +\|\cR_1\delta b\|_{\widetilde L^2_T(\dot B^1_{2,1})}\\
 &\quad \leq C\int_0^T\|b^2(t)\|_{\dot B^2_{2,1}}
    \|\delta b(t)\|_{\dot B^1_{2,1}}\,dt
       +C\int_0^T
    \|\delta u\cdot\nabla b^1\|_{\dot B^1_{2,1}}\,dt.
 \end{split}
\end{equation}
Since $\dot B^1_{2,1}$ is an algebra and
$\delta u=\cR_1\cR^\perp\delta b$, the last term satisfies
\[
 \|\delta u\cdot\nabla b^1\|_{\dot B^1_{2,1}}
 \lesssim \|\delta b\|_{\dot B^1_{2,1}}
          \|b^1\|_{\dot B^2_{2,1}}.
\]
For every finite $T>0$, the functions $b^1$ and $b^2$ belong to
$L^1_T(\dot B^2_{2,1})$. Therefore, \eqref{ipm-uniqueness-est} and
Gronwall's inequality give $\delta b=0$ on $[0,T]$. This completes the
proof of Theorem~\ref{them1:local}.

\subsection{Proof of Theorem \ref{them1}}
To prove the global stability, we rely on the uniform estimates of
Proposition~\ref{dc0}, assuming \eqref{a1}.
Let $T^*>0$ be the maximal lifespan of the solution constructed above.
By the uniform estimates \eqref{db2} and \eqref{db1}, we arrive at
\begin{equation}\label{ipm-global-est}
 \sup_{0\leq t<T^*}
 \left(\|b(t)\|_{\dot B^1_{2,1}}+\|b(t)\|_{\dot B^2_{2,1}}\right)
 \leq C\left(\|b_0\|_{\dot B^1_{2,1}}+\|b_0\|_{\dot B^2_{2,1}}\right).
\end{equation}
Assume that $T^*<\infty$. For any $t_0\in[0,T^*)$, we apply
Theorem~\ref{them1:local} with initial datum $b(t_0)$. By
\eqref{ipm-global-est}, the resulting existence time has a positive lower
bound independent of $t_0$. Choosing $t_0$ sufficiently close to $T^*$
contradicts the maximality of $T^*$. Hence $T^*=\infty$. Moreover,
\eqref{db2} and \eqref{db1} prove \eqref{r1} and \eqref{vel1}, whereas the
time continuity follows from Proposition~\ref{dc0}. This proves
Theorem~\ref{them1}.

\subsection{\texorpdfstring{Proof of Theorem \ref{thm:ipm-higher-sobolev}}
{Proof of Theorem \ref{thm:ipm-higher-sobolev}}}
\label{sec:ipm-higher-sobolev}

We first state the Sobolev commutator estimate needed for the energy
argument. We adapt the argument of \cite[Lemma~4.1]{RW} to $\mathbb R^2$,
retaining the critical Besov norms in the lower-order terms.

\begin{lemma}[Sobolev commutator estimate]\label{lem:ipm-sobolev-commutator}
Let $s>2$ and $\Lambda:=(-\Delta)^{\frac12}$. Let $f\in H^s(\mathbb R^2)$ be
real-valued and set $u=\mathcal R^\perp\mathcal R_1f$. Then
\[
 \big|\langle[\Lambda^s,u\cdot\nabla]f,\Lambda^sf\rangle\big|
 \leq C_s\Bigl(
 \|\partial_2u_2\|_{L^\infty}\|f\|_{H^s}^2
 +\|f\|_{\dot B^2_{2,1}}\|\mathcal R_1f\|_{H^s}^2
 +\|\mathcal R_1f\|_{\dot B^2_{2,1}}\|f\|_{H^s}\|\mathcal R_1f\|_{H^s}\Bigr),
\]
where $C_s$ depends only on $s$.
\end{lemma}
The proof of Lemma \ref{lem:ipm-sobolev-commutator} is given at the end of this subsection. We now prove Theorem~\ref{thm:ipm-higher-sobolev} regarding the uniform propagation of higher-order Sobolev regularity. Fix $s>2$ and let
$b_0\in H^s(\mathbb R^2)$ satisfy
$\|b_0\|_{\dot B^2_{2,1}}\leq\delta_s \ll 1$. Recalling that
$u=\mathcal R^\perp\mathcal R_1b$, we shall use
\[
 \operatorname{div}u=0,\qquad
 u_2=\mathcal R_1^2b,\qquad
 \partial_1b=\Lambda\mathcal R_1b,\qquad
 \langle u_2,b\rangle_{L^2}=-\|\mathcal R_1b\|_{L^2}^2.
\]

Since the linear term $\mathcal R_1^2b$ is bounded on $H^s$ and
dissipative, the standard local existence argument for IPM applies,
see \cite{CGO}. Let $\widetilde b\in C([0,T_s^*);H^s)$ be the local
Sobolev solution with initial datum $b_0$, where $[0,T_s^*)$ is its
maximal interval of existence. Since
\[
 H^s(\mathbb R^2)\hookrightarrow
 \dot B^1_{2,1}(\mathbb R^2)\cap\dot B^2_{2,1}(\mathbb R^2),
 \qquad s>2,
\]
uniqueness in Theorem~\ref{them1:local} gives $\widetilde b(t)=b(t)$ for
$t\in[0,T_s^*)$, where $b$ is the global solution of
Theorem~\ref{them1}. From now on we write $b$ for both.

The energy functional
\[
 E_s(f):=\|f\|_{L^2}^2+\|\Lambda^s f\|_{L^2}^2\simeq_s\|f\|_{H^s}^2
\]
satisfies
\begin{equation}\label{ipm-energy-identity}
 \frac12\frac{d}{dt}E_s(b)+E_s(\mathcal R_1b)
 =-\langle[\Lambda^s,u\cdot\nabla]b,\Lambda^sb\rangle.
\end{equation}
By Lemma~\ref{lem:ipm-sobolev-commutator}, we have
\begin{equation}\label{ipm-tame-sobolev}
\begin{aligned}
 &\frac12\frac{d}{dt}E_s(b)+E_s(\mathcal R_1b)
 \\
 &\quad \leq C_s\Bigl(\|b\|_{\dot B^2_{2,1}}E_s(\mathcal R_1b)
 +\|\mathcal R_1b\|_{\dot B^2_{2,1}}E_s(b)^{\frac12}E_s(\mathcal R_1b)^{\frac12}
 +\|\partial_2u_2\|_{L^\infty}E_s(b)\Bigr).
 \end{aligned}
\end{equation}

Set $a_0:=\|b_0\|_{\dot B^2_{2,1}}$. The estimates \eqref{r1} and
\eqref{vel1} imply
\begin{equation}\label{ipm-sobolev-integrable-coefficients}
 \sup_{t\geq0}\|b(t)\|_{\dot B^2_{2,1}}\leq Ca_0,\qquad
 \int_0^\infty
 \Bigl(\|\mathcal R_1b(t)\|_{\dot B^2_{2,1}}^2
       +\|\partial_2u_2(t)\|_{L^\infty}\Bigr)\,dt
 \leq C(a_0+a_0^2).
\end{equation}
Choose $\delta_s\leq\min\{\delta_0,1\}$, where $\delta_0$ is the threshold
of Theorem~\ref{them1}, so small that $C_sC\delta_s\leq\frac14$, and assume
$a_0\leq\delta_s$. Then the first term on the right-hand side of
\eqref{ipm-tame-sobolev} is bounded by $\frac14E_s(\mathcal R_1b)$ and is
absorbed by the dissipation. Moreover, by Young's inequality, one has
\[
 C_s\|\mathcal R_1b\|_{\dot B^2_{2,1}}E_s(b)^{\frac12}E_s(\mathcal R_1b)^{\frac12}
 \leq\frac14E_s(\mathcal R_1b)+C_s\|\mathcal R_1b\|_{\dot B^2_{2,1}}^2E_s(b).
\]
Consequently, we obtain
\begin{equation}\label{ipm-sobolev-gronwall}
 \frac{d}{dt}E_s(b)+E_s(\mathcal R_1b)
 \leq C_s\Bigl(\|\mathcal R_1b\|_{\dot B^2_{2,1}}^2
               +\|\partial_2u_2\|_{L^\infty}\Bigr)E_s(b).
\end{equation}
For every $T<T_s^*$, it follows from Gronwall's lemma and
\eqref{ipm-sobolev-integrable-coefficients} that
\[
 \sup_{0\leq t\leq T}E_s(b(t))
 +\int_0^TE_s(\mathcal R_1b(t))\,dt
 \leq E_s(b_0)\exp\!\bigl(C_s(a_0+a_0^2)\bigr).
\]
This bound is independent of $T<T_s^*$. Hence the local existence theorem
allows the solution to be restarted sufficiently close to $T_s^*$ and
continued beyond it. By maximality, $T_s^*=\infty$.

Using $E_s(f)\simeq_s\|f\|_{H^s}^2$, $a_0\leq1$ and letting $T\to\infty$,
we conclude that $b\in C([0,\infty);H^s)$ and
\[
 \sup_{t\geq0}\|b(t)\|_{H^s}^2
 +\int_0^\infty\|\mathcal R_1b(t)\|_{H^s}^2\,dt
 \leq C_s\|b_0\|_{H^s}^2.
\]
This proves \eqref{ipm-uniform-sobolev} and
Theorem~\ref{thm:ipm-higher-sobolev}. The smallness threshold $\delta_s$
may depend on $s$. No smallness is required of the initial $H^s$ norm.

\begin{proof}[Proof of Lemma~\ref{lem:ipm-sobolev-commutator}]
By density, we may assume that
$\widehat f\in C_c^\infty(\mathbb R^2\setminus\{0\})$.
We decompose
\[
 [\Lambda^s,u\cdot\nabla]f
 =[\Lambda^s,u_1]\partial_1f
  +[\Lambda^s,u_2]\partial_2f.
\]

\noindent
\emph{Horizontal transport.}
By the Kato--Ponce commutator estimate~\cite{Li2019},
\[
\begin{aligned}
\|[\Lambda^s,u_1]\partial_1f\|_{L^2}
\leq C_s\bigl(
&\|\nabla u_1\|_{L^\infty}
  \|\partial_1f\|_{\dot H^{s-1}}+\|u_1\|_{\dot H^s}
  \|\partial_1f\|_{L^\infty}\bigr).
\end{aligned}
\]
Since $u=\mathcal R^\perp\mathcal R_1f$, we have
\[
\begin{aligned}
\|u\|_{H^s}
 +\|\partial_1f\|_{\dot H^{s-1}}
 \leq C\|\mathcal R_1f\|_{H^s},\\
&\|\nabla u\|_{L^\infty}
 +\|\partial_1f\|_{L^\infty}
 \leq C\|\mathcal R_1f\|_{\dot B^2_{2,1}}.
\end{aligned}
\]
Consequently,
\begin{equation}\label{ipm-sobolev-horizontal-transport}
\big|\langle[\Lambda^s,u_1]\partial_1f,
             \Lambda^sf\rangle\big|
\leq C_s\|\mathcal R_1f\|_{\dot B^2_{2,1}}
 \|f\|_{H^s}\|\mathcal R_1f\|_{H^s}.
\end{equation}

\noindent
\emph{Vertical transport: commutator remainder.}
For smooth functions $f,g$, define
\[
\mathcal C_s(f,g)
:=[\Lambda^s,f]\partial_2g
+s\sum_{k=1}^2(\partial_kf)
  \Lambda^{s-2}\partial_k\partial_2g.
\]
We use the following second-order fractional commutator
estimate: for $s>2$ and $r\in\{0,1\}$,
\begin{equation}\label{ipm-sobolev-remainder-bounds}
\begin{aligned}
\|\Lambda^r\mathcal C_s(f,g)\|_{L^2}
\leq C_s\Bigl(
&\|\nabla f\|_{\dot H^s}
  \|g\|_{\dot B^{r+1}_{2,1}}+\|\nabla f\|_{\dot B^2_{2,1}}
  \|g\|_{\dot H^{s+r-1}}
\Bigr).
\end{aligned}
\end{equation}
This is a standard consequence of the higher-order
fractional Leibniz estimates~\cite{Li2019}.
\noindent
\emph{Vertical transport: exploiting the IPM structure.}
Let $\psi:=\Lambda^{-1}\mathcal R_1f$, so that
$u_2=\partial_1\psi$. We have
\[
\|\nabla\psi\|_{\dot H^s}
 \leq C\|\mathcal R_1f\|_{H^s},
\qquad
\|\nabla\psi\|_{\dot B^2_{2,1}}
 \leq C\|\mathcal R_1f\|_{\dot B^2_{2,1}}.
\]
Since spatial derivatives commute with $\Lambda^s$,
the product rule gives
\[
\mathcal C_s(u_2,f)
=\partial_1\mathcal C_s(\psi,f)
 -\mathcal C_s(\psi,\partial_1f).
\]
Using $\partial_1f=\Lambda\mathcal R_1f$
and integrating by parts, we obtain
\[
\begin{aligned}
\langle\mathcal C_s(u_2,f),\Lambda^sf\rangle
={}&-\langle\Lambda\mathcal C_s(\psi,f),
             \Lambda^s\mathcal R_1f\rangle-\langle\mathcal C_s(\psi,\Lambda\mathcal R_1f),
             \Lambda^sf\rangle.
\end{aligned}
\]
By \eqref{ipm-sobolev-remainder-bounds} with $r=1$,
\[
\begin{aligned}
\|\Lambda\mathcal C_s(\psi,f)\|_{L^2}
\leq C_s\Bigl(
&\|\mathcal R_1f\|_{H^s}\|f\|_{\dot B^2_{2,1}}+\|\mathcal R_1f\|_{\dot B^2_{2,1}}\|f\|_{H^s}
\Bigr).
\end{aligned}
\]
Similarly, using \eqref{ipm-sobolev-remainder-bounds}
with $r=0$ and
\[
\|\Lambda\mathcal R_1f\|_{\dot B^1_{2,1}}
 =\|\mathcal R_1f\|_{\dot B^2_{2,1}},
\qquad
\|\Lambda\mathcal R_1f\|_{\dot H^{s-1}}
 \leq\|\mathcal R_1f\|_{H^s},
\]
we get
\[
\|\mathcal C_s(\psi,\Lambda\mathcal R_1f)\|_{L^2}
\leq C_s\|\mathcal R_1f\|_{\dot B^2_{2,1}}
          \|\mathcal R_1f\|_{H^s}.
\]
Combining these estimates yields
\[
\begin{aligned}
|\langle\mathcal C_s(u_2,f),\Lambda^sf\rangle|
\leq C_s\Bigl(
&\|f\|_{\dot B^2_{2,1}}\|\mathcal R_1f\|_{H^s}^2+\|\mathcal R_1f\|_{\dot B^2_{2,1}}
  \|f\|_{H^s}\|\mathcal R_1f\|_{H^s}
\Bigr).
\end{aligned}
\]

\noindent
\emph{Vertical transport: conclusion.}
By definition of the remainder,
\[
[\Lambda^s,u_2]\partial_2f
=-s\sum_{k=1}^2(\partial_ku_2)
  \Lambda^{s-2}\partial_k\partial_2f
 +\mathcal C_s(u_2,f).
\]
Moreover,
\[
\begin{aligned}
\|\Lambda^{s-2}\partial_1\partial_2f\|_{L^2}
 \leq\|\mathcal R_1f\|_{H^s}, \quad \|\Lambda^{s-2}\partial_2^2f\|_{L^2}
 \leq\|f\|_{H^s}, \quad \|\partial_1u_2\|_{L^\infty}
 \leq C\|\mathcal R_1f\|_{\dot B^2_{2,1}}.
\end{aligned}
\]
Hence, by H\"older's inequality and the preceding
remainder estimate,
\begin{equation*}
\big|\langle[\Lambda^s,u_2]\partial_2f,
                 \Lambda^sf\rangle\big|\leq C_s\Bigl(
 \|\partial_2u_2\|_{L^\infty}\|f\|_{H^s}^2
 +\|f\|_{\dot B^2_{2,1}}\|\mathcal R_1f\|_{H^s}^2
 +\|\mathcal R_1f\|_{\dot B^2_{2,1}}
  \|f\|_{H^s}\|\mathcal R_1f\|_{H^s}
 \Bigr).
\end{equation*}
Adding \eqref{ipm-sobolev-horizontal-transport}
and this last inequality completes the proof.
\end{proof}

\section{Linear decay}\label{sec:linear-decay}

In this section, we prove optimal decay estimates for the linearized IPM system.

\begin{theorem}[Linear decay estimates]\label{thm:linear-decay}
Let ${s\geq0}$, $\sigma>0$ and
$b_0\in{\dot B^s_{2,1}\cap\mathcal X^{s+\sigma,-\sigma}}$.
Let $(b,u)$ satisfy
\begin{equation}\label{linearIPM}
\left\{
\begin{aligned}
 &\partial_t b-\cR_1^2b=0,\\ 
 & u=\cR^\perp\cR_1b,\\
 &b|_{t=0}=b_0.
 \end{aligned}
\right.
\end{equation}
For every $t\geq0$, we have the following bounds with constants independent of $t$:
\begin{equation}\label{linear-angular-enhanced}
\left\{
\begin{alignedat}{2}
{\|\Lambda^s b(t)\|_{L^2}}
 &\leq C(1+t)^{-\frac{\sigma}{2}},
&\qquad {\|\nabla\Lambda^{s-2}b(t)\|_{L^\infty}}
 &\leq C(1+t)^{-\frac{1}{4}-\frac{\sigma}{2}},\\[2pt]
{\|\Lambda^s u_1(t)\|_{L^2}}
 &\leq C(1+t)^{-\frac{1}{2}-\frac{\sigma}{2}},
&\qquad {\|\nabla\Lambda^{s-2}u_1(t)\|_{L^\infty}}
 &\leq C(1+t)^{-\frac{3}{4}-\frac{\sigma}{2}},\\[2pt]
{\|\Lambda^s u_2(t)\|_{L^2}}
 &\leq C(1+t)^{-1-\frac{\sigma}{2}},
&\qquad {\|\nabla\Lambda^{s-2}u_2(t)\|_{L^\infty}}
 &\leq C(1+t)^{-\frac{5}{4}-\frac{\sigma}{2}}.
\end{alignedat}
\right.
\end{equation}
The six exponents are sharp within this data class. More precisely,
there is one nonzero datum $b_0$ in the same space such that, writing
$f_0=b$, $f_1=u_1$ and $f_2=u_2$, one has, for $m=0,1,2$ and $t\geq0$,
\begin{equation}\label{linear-fixed-L2}
 {\|\Lambda^s f_m(t)\|_{L^2}}
 \geq {C}(1+t)^{-\frac{m+\sigma}{2}},
 \qquad
 {\|\nabla\Lambda^{s-2}f_m(t)\|_{L^\infty}}
 \geq {C}(1+t)^{-\frac{m+\sigma}{2}-\frac14}.
\end{equation}
\end{theorem}

\subsection{The spectrally localized semigroup estimate}
The semigroup $S(t):=e^{t\cR_1^2}$ has the Fourier symbol $ e^{-t\frac{\xi_1^2}{|\xi|^2}}$. For $t\geq1$, the relevant angular scale is determined by
\[
t\frac{\xi_1^2}{|\xi|^2}\simeq1,
\qquad
\frac{|\xi_1|}{|\xi|}\simeq t^{-\frac12}.
\]
The dissipation is weak when
$\frac{|\xi_1|}{|\xi|}\ll t^{-\frac12}$, and decay in this region
requires additional angular control of the initial data.
For larger angular frequencies, the semigroup provides
stronger decay. More precisely, we have
\[
\left(\frac{|\xi_1|}{|\xi|}\right)^\sigma
e^{-t\frac{\xi_1^2}{|\xi|^2}}
\lesssim (1+t)^{-\frac{\sigma}{2}},
\qquad \sigma>0.
\]

On the Fourier support of $\dot\Delta_j\dot\Delta_k^h f$, 
$|\xi|\simeq 2^j$ and $|\xi_1|\simeq 2^k$ follow. Moreover, $\dot\Delta_j\dot\Delta_k^h$
vanishes unless $k\leq j+N_0$. Using Plancherel's theorem, we obtain,
uniformly with respect to $j$ and $k$,
\begin{equation}\label{linear-block-decay}
 \|\dot\Delta_j\dot\Delta_k^h\cR_1^m S(t)f\|_{L^2}
 \leq C2^{m(k-j)}e^{-ct2^{2(k-j)}}
       \|\dot\Delta_j\dot\Delta_k^h f\|_{L^2},
 \qquad m=0,1,2.
\end{equation}

We shall use the following elementary summation estimate:
\begin{equation}\label{linear-elementary-decay}
 \sum_{\ell\leq N_0}2^{\beta\ell}e^{-ct2^{2\ell}}
 \leq C_\beta(1+t)^{-\frac{\beta}{2}},\qquad \beta>0.
\end{equation}
Indeed, the estimate is immediate when $0\leq t\leq1$. Suppose now
that $t\geq1$ and choose $\ell_t\in\mathbb Z$ such that $ 2^{\ell_t}\leq t^{-\frac12}<2^{\ell_t+1}$. We split the sum at $\ell_t$. For $\ell\leq\ell_t$, we obtain
\[
 \sum_{\ell\leq\ell_t}
 2^{\beta\ell}e^{-ct2^{2\ell}}
 \leq C2^{\beta\ell_t}
 \leq Ct^{-\frac{\beta}{2}}.
\]
Since $t2^{2\ell_t}\geq\frac{1}{4}$, one has
\[
\begin{aligned}
\sum_{\ell_t<\ell\leq N_0}
2^{\beta\ell}e^{-ct2^{2\ell}}
&\leq
2^{\beta\ell_t}
\sum_{\ell>\ell_t}
2^{\beta(\ell-\ell_t)}
e^{-c4^{\ell-\ell_t-1}}\leq C_\beta 2^{\beta\ell_t}
\leq C_\beta t^{-\frac{\beta}{2}}.
\end{aligned}
\]
Consequently, \eqref{linear-elementary-decay} is proved.

\subsection{Upper bounds}
\begin{proof}[Proof of the upper bounds]
By the definition of $\mathcal X^{s+\sigma,-\sigma}$,
for $k=j+\ell$, we have
\[
\begin{aligned}
2^{sj}\|\dot\Delta_j\dot\Delta_{j+\ell}^h b_0\|_{L^2}
&\leq 2^{\sigma\ell}2^{(s+\sigma)j}
\sup_k\left(
2^{-\sigma k}
\|\dot\Delta_j\dot\Delta_k^h b_0\|_{L^2}
\right).
\end{aligned}
\]
Using \eqref{linear-block-decay} and summing over the
horizontal blocks, we obtain
\[
\begin{aligned}
\|\cR_1^mS(t)b_0\|_{\dot B^s_{2,1}}
&\lesssim \sum_j\sum_{\ell\leq N_0}
2^{sj+m\ell}e^{-ct2^{2\ell}}
\|\dot\Delta_j\dot\Delta_{j+\ell}^h b_0\|_{L^2}\lesssim \|b_0\|_{\mathcal X^{s+\sigma,-\sigma}} 
\sum_{\ell\leq N_0}2^{(m+\sigma)\ell}
e^{-ct2^{2\ell}}.
\end{aligned}
\]
Similarly, since
$2^{(s-\frac{1}{2})j+\frac{1}{2}k}=2^{sj+\frac{1}{2}\ell}$,
we have
\[
\begin{aligned}
\|\cR_1^mS(t)b_0\|_{\cB^{s-\frac 12,\frac 12}}
&\lesssim \sum_j\sum_{\ell\leq N_0}
2^{sj+(m+\frac12)\ell}e^{-ct2^{2\ell}}
\|\dot\Delta_j\dot\Delta_{j+\ell}^h b_0\|_{L^2}\\
&\lesssim \|b_0\|_{\mathcal X^{s+\sigma,-\sigma}} 
\sum_{\ell\leq N_0}
2^{(m+\sigma+\frac12)\ell}e^{-ct2^{2\ell}}.
\end{aligned}
\]
Applying \eqref{linear-elementary-decay} with
$\beta=m+\sigma$ and $\beta=m+\sigma+\frac{1}{2}$,
respectively, we obtain
\begin{equation}\label{linear-angular-norms}
\begin{aligned}
&\|\cR_1^mS(t)b_0\|_{\dot B^s_{2,1}}
\lesssim (1+t)^{-\frac{m+\sigma}{2}}
\|b_0\|_{\mathcal X^{s+\sigma,-\sigma}},\\
&\|\cR_1^mS(t)b_0\|_{\cB^{s-\frac 12,\frac 12}}
\lesssim (1+t)^{-\frac{m}{2}-\frac14-\frac{\sigma}{2}}
\|b_0\|_{\mathcal X^{s+\sigma,-\sigma}}.
\end{aligned}
\end{equation}
Since $u_1=-\cR_2\cR_1b$, $u_2=\cR_1^2b$ and $\cR_2$ is
bounded on these block norms, we take $m=0,1,2$ for $b,u_1,u_2$,
respectively. Using $\dot B^s_{2,1}\hookrightarrow\dot H^s$
and Bernstein's inequality in the form
$\|\nabla\Lambda^{s-2}f\|_{L^\infty}
\lesssim\|f\|_{\cB^{s-\frac 12,\frac 12}}$,
one obtains \eqref{linear-angular-enhanced}
from \eqref{linear-angular-norms}.
\end{proof}

\subsection{Sharpness for a fixed datum}

To prove the optimality of the decay exponents in \eqref{linear-angular-enhanced}, we prove the lower bounds for a class of data.

\begin{proof}
{{Fix $\sigma>0$ as in Theorem~\ref{thm:linear-decay}.
Choose even, nonnegative cutoff functions
$\chi_h\in C_c^\infty((-\frac{1}{2},\frac{1}{2}))$ and
$\chi_v\in C_c^\infty((-2,-1)\cup(1,2))$,
acting on the horizontal and vertical Fourier variables
$\xi_1$ and $\xi_2$, respectively.
Assume that $\chi_h=1$ on $[-\frac{1}{4},\frac{1}{4}]$ and
$\chi_v$ is strictly positive on some interval contained in $(1,2)$. Define
\begin{align}\label{linear-sharp-datum-physical}
b_0(x_1,x_2)
:=\mathcal F^{-1}\bigl(
|\xi_1|^{\sigma-\frac{1}{2}}\chi_h(\xi_1)\chi_v(\xi_2)
\bigr)(x_1,x_2).
\end{align}
By the support properties of $\chi_h$ and $\chi_v$,
the Fourier support of $b_0$ is contained in
$\{|\xi_1|<\frac{1}{2},\ 1<|\xi_2|<2\}$, which consists of
two rectangles lying in the fixed annulus
$\{1\leq|\xi|\leq3\}$. Indeed, every $\xi$ in these rectangles satisfies
$1<|\xi_2|\leq|\xi|<\frac{\sqrt{17}}{2}<3$. Since $\sigma>0$, the function
$|\xi_1|^{\sigma-\frac{1}{2}}$ is locally square-integrable
near $\xi_1=0$. Hence $\widehat b_0\in L^2$.
Moreover, the Fourier support of $b_0$ intersects only
finitely many radial dyadic annuli, which implies
$b_0\in\dot B^s_{2,1}$. Denote $J=\{j\in\Z: \log_2 (\frac 38) \leq j\leq 2\}.$ Based on the Fourier support of $b_0$ and the support of the Littlewood--Paley block, one has\begin{align*}
    \|b_0\|_{\dot B^s_{2,1}}=&\sum_{j\in \Z}2^{js}\|\dD_j b_0\|_{L^2}=\sum_{j\in J}2^{js}\|\dD_j b_0\|_{L^2}\lesssim \|\hat{b_0}\|_{L^2}.
\end{align*}

To verify the additional angular condition, Plancherel's
theorem and the horizontal frequency localization give,
for every $j,k\in\mathbb Z$,
\[
\begin{aligned}
\|\dot\Delta_j\dot\Delta_k^h b_0\|_{L^2}^2
&\lesssim
\int_{|\xi_1|\simeq2^k}
|\xi_1|^{2\sigma-1}\,d\xi_1
\lesssim 2^{2\sigma k}.
\end{aligned}
\]
Consequently, it follows from the definition of $\mathcal X^{s+\sigma,-\sigma}$ that
\[
\begin{aligned}
\|b_0\|_{\mathcal X^{s+\sigma,-\sigma}}
&=\sum_j2^{(s+\sigma)j}\sup_k
2^{-\sigma k}
\|\dot\Delta_j\dot\Delta_k^h b_0\|_{L^2}\lesssim\sum_{j\in J}2^{(s+\sigma)j}<\infty,
\end{aligned}
\]
Thus $b_0\in\dot B^s_{2,1}
\cap\mathcal X^{s+\sigma,-\sigma}$.}}

We first prove the $L^2$ lower bounds. By Plancherel's theorem,
for $m=0,1,2$, we have
\[
 {\|\Lambda^s\cR_1^mS(t)b_0\|_{L^2}^2}
 \simeq
 \int_{\mathbb R^2}
 {|\xi|^{2s-2m}}|\xi_1|^{2m}
 e^{-2t\frac{\xi_1^2}{|\xi|^2}}
 |\widehat b_0(\xi)|^2\,d\xi.
\]
We restrict the integral to $0<\xi_1<\frac{1}{4}$ and $1<\xi_2<2$.
On this region, we have $\chi_h(\xi_1)=1$ and $|\xi|\geq1$.
Substituting the definition of $b_0$ and using {$|\xi|^{2s-2m}\geq c_{s,m}>0$ on this fixed annulus},
we obtain
\[
\begin{aligned}
 {\|\Lambda^s\cR_1^mS(t)b_0\|_{L^2}^2}
 &\geq c
 \left(\int_1^2\chi_v(\xi_2)^2\,d\xi_2\right)
 \int_0^{\frac{1}{4}}r^{2(m+\sigma)-1}e^{-2tr^2}\,dr\geq c\int_0^{\frac{1}{4}}r^{2(m+\sigma)-1}e^{-2tr^2}\,dr,
\end{aligned}
\]
where we used $\int_1^2\chi_v(\xi_2)^2\,d\xi_2>0$.
For $t\geq1$, we further restrict the integral to
$0<r<(4\sqrt t)^{-1}$. On this interval, we have
$e^{-2tr^2}\geq e^{-\frac{1}{8}}$, and hence
\[
\begin{aligned}
 {\|\Lambda^s\cR_1^mS(t)b_0\|_{L^2}^2}
 &\geq c\int_0^{\frac{1}{4\sqrt t}}r^{2(m+\sigma)-1}\,dr=\frac{c}{2(m+\sigma)}(4\sqrt t)^{-2(m+\sigma)}
 \gtrsim t^{-(m+\sigma)}.
\end{aligned}
\]

For $u_1=-\cR_2\cR_1b$, the additional operator $\cR_2$ changes
the estimate only by a fixed constant. Therefore, writing
$f_0:=b$, $f_1:=u_1$ and $f_2:=u_2$, we have
\[
 {\|\Lambda^sf_m(t)\|_{L^2}}
 \geq c(1+t)^{-\frac{m+\sigma}{2}},\qquad m=0,1,2.
\]
For $0\leq t\leq1$, the semigroup is uniformly bounded
below on the support of $\widehat b_0$, so the same conclusion follows
after decreasing $c$.
{This proves the first lower bound in \eqref{linear-fixed-L2}.}

We next prove the $L^\infty$ lower bounds. Choose nonnegative,
nonzero cutoffs $\eta_h,\eta_v\in C_c^\infty(\mathbb R)$ such that
$\operatorname{supp}\eta_h\subset(0,\frac{1}{2})$, $\chi_h=1$ on
$\operatorname{supp}\eta_h$ and
$\operatorname{supp}\eta_v\subset\{\xi_2>0:\chi_v(\xi_2)>0\}$.
We set $r_t:=c_0(1+t)^{-\frac12}$, where $c_0>0$ is sufficiently small,
and define $P_t$ by
$$
\widehat{P_tf}(\xi):=\eta_h(\frac{\xi_1}{r_t})\eta_v(\xi_2)\widehat f(\xi).
$$
Its kernel is, up to a fixed normalization constant,
$K_t(x_1,x_2)=r_t\check\eta_h(r_tx_1)\check\eta_v(x_2)$.
We have $\|K_t\|_{L^1}\leq C$ uniformly in $t$, and hence
$\|P_tf\|_{L^\infty}\leq C\|f\|_{L^\infty}$.

On the support of $P_t$, we have
$|\xi_1|\simeq r_t$, $|\xi_2|\simeq1$ and
$e^{-t\frac{\xi_1^2}{|\xi|^2}}\geq c$. After multiplication by fixed complex
numbers of modulus one, the Fourier integrands corresponding to
${\partial_2\Lambda^{s-2}b}$, ${\partial_2\Lambda^{s-2}u_1}$ and ${\partial_2\Lambda^{s-2}u_2}$ are nonnegative
on this support. Hence no cancellation occurs in the Fourier inversion
formula at the origin, and we obtain
\[
 {|P_t\partial_2\Lambda^{s-2}f_m(t)(0)|}
 \geq cr_t^{m+\sigma+\frac{1}{2}},\qquad m=0,1,2.
\]
Using the uniform $L^\infty$ boundedness of $P_t$,
we conclude that
\[
 {\|\nabla\Lambda^{s-2}f_m(t)\|_{L^\infty}}
 \geq C^{-1}{|P_t\partial_2\Lambda^{s-2}f_m(t)(0)|}
 \geq c(1+t)^{-\frac{m+\sigma}{2}-\frac14}.
\]
Together with the preceding $L^2$ estimates, this proves
\eqref{linear-fixed-L2} for one fixed nonzero initial datum.
\end{proof}

\section{Nonlinear decay}\label{sec:decay}

In this section, we prove decay estimates at lower order and at critical
regularity. Throughout, $h=b\circ X$ and $U=u\circ X$ are the isopycnal
coordinate variables of Lemma~\ref{lem:density-reduction}. From \eqref{dahi} and
$\|b_0\|_{\dot B^2_{2,1}}\leq\delta_0$, we have
\begin{equation}\label{decay-critical-coefficients}
 \|h\|_{\widetilde L_T^\infty(\cB^{\frac 32,\frac 12})}
 +\|h\|_{\widetilde L_T^2(\cB^{\frac 12,\frac 32})}
 +\|h\|_{L_T^1(\cB^{-\frac 12,\frac 52})}\leq C\delta_0\ll 1.
\end{equation}
Using \eqref{dahi}, we have
$\|\partial_1U_1\|_{L_T^1(L^\infty)}\leq C\delta_0$.

\begingroup
\subsection{$L^2$ decay of the density}

\subsubsection{Energy and interpolation estimates}

We prove the $L^2$ decay estimate for $b$ under the additional assumption
$h_0\in\mathcal X^{\frac 12,-\frac 12}$. 
Since $\operatorname{div}u=0$, the transport term vanishes in
the $L^2$ energy identity:
\begin{align}\label{ddtL2}
 \frac12\frac d{dt}\|b\|_{L^2}^2+\|\cR_1b\|_{L^2}^2=0.
\end{align}

To quantify how the dissipation controls the energy, we establish
the following interpolation estimate:
\begin{align}\label{interL22}
\|h\|_{L^2}
\lesssim
\|h\|_{\mathcal X^{\frac 12,-\frac 12}}^{\frac23}
\|\cR_1h\|_{L^2}^{\frac{1}{3}}.
\end{align}
{
To prove \eqref{interL22}, we use the anisotropic
Littlewood--Paley decomposition and the frequency localization
$\frac{|\xi_1|}{|\xi|}\simeq2^{k-j}$.
Using the support restriction $k\leq j+N_0$,  and the equivalence $\cB^{0,0}_{2,2}\sim L^2$, we obtain that, for every
$\eta>0$,
\[
\begin{aligned}
 \|h\|_{L^2}^2
 &\lesssim
 \sum_{j\in\mathbb Z}
 \sum_{\substack{k\leq j+N_0\\2^{k-j}\leq\eta}}
 \|\dot\Delta_j\dot\Delta_k^h h\|_{L^2}^2
 +\sum_{j\in\mathbb Z}
 \sum_{\substack{k\leq j+N_0\\2^{k-j}>\eta}}
 \|\dot\Delta_j\dot\Delta_k^h h\|_{L^2}^2\\
 &\lesssim
 \sum_{j\in\mathbb Z}
 \left(2^{\frac{j}{2}}\sup_{k\in\mathbb Z}2^{-\frac{k}{2}}
       \|\dot\Delta_j\dot\Delta_k^h h\|_{L^2}\right)^2
 \sum_{\substack{k\leq j+N_0\\2^{k-j}\leq\eta}}2^{k-j}
 +\eta^{-2}\|\cR_1h\|_{L^2}^2\\
 &\lesssim
 \eta\|h\|_{\mathcal X^{\frac12,-\frac12}}^2
 +\eta^{-2}\|\cR_1h\|_{L^2}^2.
\end{aligned}
\]
For the first sum, we use $ \sum_{\substack{k\leq j+N_0,\, 2^{k-j}\leq\eta}}2^{k-j}
 \leq2\eta$ for $j\in\mathbb Z$ and $\ell^1\hookrightarrow\ell^2$.
For the second sum, frequency localization and Plancherel's theorem imply
\[
\begin{aligned}
\sum_{j\in\mathbb Z}
 \sum_{\substack{k\leq j+N_0\\2^{k-j}>\eta}}
 \|\dot\Delta_j\dot\Delta_k^h h\|_{L^2}^2\lesssim  \sum_{j\in\mathbb Z}
 \sum_{\substack{k\leq j+N_0\\2^{k-j}>\eta}} \eta^{-2}\|\dD_j\dD_k^h\cR_1h\|_{L^2}^2\lesssim \eta^{-2} \|\cR_1 h\|^2_{\cB^{0,0}_{2,2}}
 \lesssim
 \eta^{-2}\|\mathcal R_1h\|_{L^2}^2.
\end{aligned}
\]

The case $h=0$ is clear. 
 For $h\neq0$, choose
\[
 \eta=
 \left(
 \frac{\|\cR_1h\|_{L^2}}
      {\|h\|_{\mathcal X^{\frac12,-\frac12}}}
 \right)^{\frac23},
\]
which yields \eqref{interL22}. 
}

The identities $\partial_1h=J(\partial_1b)\circ X$ and
$dy=J^{-1}dx$ imply, for every test function $\varphi$,
\[
 \begin{aligned}
 |\langle\partial_1h,\varphi\rangle|
 &=|\langle\partial_1b,\varphi\circ Y\rangle|\leq\|\cR_1b\|_{L^2}\|\nabla(\varphi\circ Y)\|_{L^2}
 \lesssim\|\cR_1b\|_{L^2}\|\nabla\varphi\|_{L^2},
 \end{aligned}
\]
where the last inequality follows from the uniform bounds on $DY$
and the Jacobian. Thus, by duality,
\[
 \|\cR_1h\|_{L^2}
 =\|\partial_1h\|_{\dot H^{-1}}
 \lesssim\|\cR_1b\|_{L^2}.
\]
Since $\|h\|_{L^2}\simeq\|b\|_{L^2}$, we conclude that
\[
 \|b\|_{L^2}
 \lesssim
 \|h(t)\|_{\mathcal X^{\frac 12,-\frac 12}}^{\frac23}
 \|\cR_1b\|_{L^2}^{\frac{1}{3}},
 \qquad
 \|\cR_1b\|_{L^2}^2
 \gtrsim
 \frac{\|b\|_{L^2}^6}{\|h(t)\|_{\mathcal X^{\frac 12,-\frac 12}}^4}.
\]
Combining this bound with the energy identity, we obtain
\begin{equation}\label{int:bh}
 \frac d{dt}\|b\|_{L^2}^2
 +\frac{c}{\|h(t)\|_{\mathcal X^{\frac 12,-\frac 12}}^4}\|b\|_{L^2}^6
 \leq0.
\end{equation}
We next control the angular norm appearing in this inequality.

\subsubsection{Propagation of the angular norm}

To propagate regularity in $\mathcal X^{\frac12,-\frac12}$, we estimate
$h$ in isopycnal coordinates. For this purpose, we need the following
bounds on $G(h)$ in \eqref{dc6}.

\begin{lemma}\label{lem:mixed-source-estimates}
Let $G_1,G_2$ be defined by \eqref{def:g1g2}.
Then, under the assumption \eqref{decay-critical-coefficients} with $\delta_0\ll 1$, we have
\begin{equation}\label{zero-mixed-source}
 \begin{aligned}
 &\|G_1\|_{\widetilde L_T^2(\mathcal X^{\frac12,-\frac12})}
 +\|G_2\|_{\widetilde L_T^1(\mathcal X^{\frac12,-\frac12})}\lesssim\delta_0\left(
 \|h\|_{\widetilde L_T^\infty(\mathcal X^{\frac12,-\frac12})}
 +\|h\|_{\widetilde L_T^2(\mathcal X^{-\frac12,\frac12})}
 \right)
 \end{aligned}
\end{equation}
and
\begin{equation}\label{decay-mixed-source}
 \begin{aligned}
 \|G(h)\|_{\widetilde L_T^1(\mathcal X^{\frac 52,-\frac 12})}
 \lesssim\delta_0\bigl(
 \|h\|_{\widetilde L_T^\infty(\mathcal X^{\frac 52,-\frac 12})}+\|h\|_{\widetilde L_T^2(\mathcal X^{\frac32,\frac12})}
 +\|h\|_{\widetilde L_T^1(\mathcal X^{\frac12,\frac32})}\bigr).
 \end{aligned}
\end{equation}
\end{lemma}
The proof of Lemma \ref{lem:mixed-source-estimates} is given in Appendix~\ref{app:mixed-source}.

\begin{lemma}\label{lem:zero-angular-propagation}
Under the assumptions of Theorem~\ref{thm:decay}, uniformly for $T>0$, we have
\begin{equation}\label{zero-angular-propagation}
 \|h\|_{\widetilde L_T^\infty(\mathcal X^{\frac 12,-\frac 12})}
 +\|h\|_{\widetilde L_T^2(\mathcal X^{-\frac12,\frac12})}
 \leq C\|h_0\|_{\mathcal X^{\frac 12,-\frac 12}}.
\end{equation}
\end{lemma}

{
\begin{proof}
We estimate each mixed block as in \eqref{disbk}, using
\eqref{dc-G-decomposition}. Taking the supremum in time, then
square roots and combining the estimates in the
$\mathcal X^{\frac12,-\frac12}$ norm, we obtain
\begin{equation}\label{zero-angular-energy}
 \begin{aligned}
 &\|h\|_{\widetilde L_T^\infty(\mathcal X^{\frac12,-\frac12})}
 +\|h\|_{\widetilde L_T^2(\mathcal X^{-\frac12,\frac12})}\\
 &\quad\leq C\|h_0\|_{\mathcal X^{\frac12,-\frac12}}
 +C\|G_1\|_{\widetilde L_T^2(\mathcal X^{\frac12,-\frac12})}^{\frac12}
    \|h\|_{\widetilde L_T^2(\mathcal X^{-\frac12,\frac12})}^{\frac12}\\
 &\qquad
 +C\|G_2\|_{\widetilde L_T^1(\mathcal X^{\frac12,-\frac12})}^{\frac12}
    \|h\|_{\widetilde L_T^\infty(\mathcal X^{\frac12,-\frac12})}^{\frac12}\\
 &\qquad
 +C\left(
    \sum_j2^{\frac{j}{2}}\sup_k2^{-\frac{k}{2}}
    \|[U_1,\dot\Delta_j\dot\Delta_k^h]\partial_1h
      \|_{L_T^1(L^2)}
    \right)^{\frac12}
    \|h\|_{\widetilde L_T^\infty(\mathcal X^{\frac12,-\frac12})}^{\frac12}\\
 &\qquad
 +C\|\partial_1U_1\|_{L_T^1(L^\infty)}^{\frac12}
    \|h\|_{\widetilde L_T^\infty(\mathcal X^{\frac12,-\frac12})}.
 \end{aligned}
\end{equation}
We first control the last term in \eqref{zero-angular-energy}.
Using the fourth term in \eqref{dahi}, we have
$\|\partial_1U_1\|_{L_T^1(L^\infty)}\leq C\delta_0$.
For sufficiently small $\delta_0$, we absorb the last term of
\eqref{zero-angular-energy} into its left-hand side.

Note that $G_1$ and $G_2$ have been addressed in \eqref{zero-angular-energy} by \eqref{zero-mixed-source}. To handle the commutator term, we apply \eqref{decay-mixed-commutator}
with $s=\frac{1}{2}$ and time exponents $(2,2)$, $(1,\infty)$
and $(2,2)$ to its three terms, respectively. We obtain
\begin{equation}\label{zero-mixed-commutator-time}
 \begin{aligned}
 &\sum_j2^{\frac{j}{2}}\sup_k2^{-\frac{k}{2}}
 \|[U_1,\dot\Delta_j\dot\Delta_k^h]\partial_1h
   \|_{L_T^1(L^2)}\\
 &\quad\lesssim
 \|U_1\|_{\widetilde L_T^2(\cB^{\frac32,\frac12})}
 \|h\|_{\widetilde L_T^2(\mathcal X^{-\frac12,\frac12})}\\
 &\qquad+
 \|\partial_1U_1\|_{L_T^1(\cB^{\frac12,\frac12})}
 \|h\|_{\widetilde L_T^\infty(\mathcal X^{\frac12,-\frac12})}+
 \|U_1\|_{\widetilde L_T^2(\mathcal X^{\frac12,-\frac12})}
 \|h\|_{\widetilde L_T^2(\cB^{\frac12,\frac32})}\\
 &\quad\lesssim\delta_0\left(
 \|h\|_{\widetilde L_T^\infty(\mathcal X^{\frac12,-\frac12})}
 +\|h\|_{\widetilde L_T^2(\mathcal X^{-\frac12,\frac12})}
 \right),
 \end{aligned}
\end{equation}
where the last inequality follows from
\eqref{decay-critical-coefficients}, \eqref{zero-mixed-gradient-time} and \eqref{zero-mixed-velocity-time}.

Substituting \eqref{zero-mixed-source} and \eqref{zero-mixed-commutator-time} into 
\eqref{zero-angular-energy}, squaring and using Young's
inequality, we obtain
\begin{equation}\label{zero-angular-closure}
 \begin{aligned}
 &\|h\|_{\widetilde L_T^\infty(\mathcal X^{\frac12,-\frac12})}^2
 +\|h\|_{\widetilde L_T^2(\mathcal X^{-\frac12,\frac12})}^2\\
 &\quad\leq C\|h_0\|_{\mathcal X^{\frac12,-\frac12}}^2+C\delta_0\left(
 \|h\|_{\widetilde L_T^\infty(\mathcal X^{\frac12,-\frac12})}^2
 +\|h\|_{\widetilde L_T^2(\mathcal X^{-\frac12,\frac12})}^2
 \right).
 \end{aligned}
\end{equation}
For sufficiently small $\delta_0$, we absorb the last term of
\eqref{zero-angular-closure} and obtain
\eqref{zero-angular-propagation}.
\end{proof}
}

\subsubsection{Decay of the density}

The case $b_0=0$ is immediate, so we assume that $b_0\neq0$. We prove that, for every $t\geq0$,
\begin{align}\label{L214}
 \|b(t)\|_{L^2}
 \leq
 \left(
 \|b_0\|_{L^2}^{-4}
 +\frac{ct}{\|h_0\|_{\mathcal X^{\frac 12,-\frac 12}}^4}
 \right)^{-\frac14}
 \lesssim(1+t)^{-\frac14}.
\end{align}

\begin{proof}[Proof of \eqref{L214}] Combining the above estimate \eqref{int:bh} with \eqref{zero-angular-propagation}, we obtain
\[
 \|b\|_{L^2}
 \lesssim
 \|h_0\|_{\mathcal X^{\frac 12,-\frac 12}}^{\frac23}
 \|\cR_1b\|_{L^2}^{\frac{1}{3}},
 \qquad
 \|\cR_1b\|_{L^2}^2
 \gtrsim
 \frac{\|b\|_{L^2}^6}{\|h_0\|_{\mathcal X^{\frac 12,-\frac 12}}^4}.
\]
Combining this bound with \eqref{ddtL2}, we obtain
\[
 \frac d{dt}\|b\|_{L^2}^2
 +\frac{c}{\|h_0\|_{\mathcal X^{\frac 12,-\frac 12}}^4}\|b\|_{L^2}^6
 \leq0.
\]
Solving this differential inequality proves \eqref{L214}.
\end{proof}

\endgroup

\subsection{Anisotropic decay and gradient estimates}

The interpolation and dissipation mechanism is the same as for the
$L^2$ estimate. The additional issue is the transport commutator in
the localized energy estimates: the transport term vanishes in the
unlocalized $L^2$ energy identity, whereas its localized contribution
must be estimated using the commutator estimate \eqref{hc1}.

We sum \eqref{dck} with weight
$2^{\frac{3j}{2}+\frac{k}{2}}$. Using \eqref{hc1}, \eqref{dcb}--\eqref{dci} and
{Proposition~\ref{prop:inter}}, we bound the nonlinear terms by
$C\|h\|_{\cB^{\frac 32,\frac 12}}\|h\|_{\cB^{-\frac 12,\frac 52}}$.
Since \eqref{dahi} implies
$\sup_{t\geq0}\|h(t)\|_{\cB^{\frac 32,\frac 12}}\leq C\delta_0\ll 1$, we arrive at
\begin{align}\label{ddthigh}
 \frac{d}{dt}\|h(t)\|_{\cB^{\frac 32,\frac 12}}
 +c\|h(t)\|_{\cB^{-\frac 12,\frac 52}}\leq0.
\end{align}
Applying Proposition~\ref{prop:real-interpolation}, we obtain
\[
 \|h(t)\|_{\cB^{-\frac 12,\frac 52}}
 \gtrsim
 \frac{\|h(t)\|_{\cB^{\frac 32,\frac 12}}^3}
      {\|h(t)\|_{\mathcal X^{\frac 52,-\frac 12}}^2}.
\]
Substituting this bound into \eqref{ddthigh}, we obtain
\[
 \frac{d}{dt}\|h(t)\|_{\cB^{\frac 32,\frac 12}}
 +\frac{c\|h(t)\|_{\cB^{\frac 32,\frac 12}}^3}
       {\|h(t)\|_{\mathcal X^{\frac 52,-\frac 12}}^2}
 \leq0.
\]
This motivates us to propagate the angular norm in the denominator uniformly in time.

\subsubsection{Propagation of the angular norm}

We prove a uniform bound for $h$ in
$\mathcal X^{\frac 52,-\frac 12}$.

\begin{lemma}\label{lem:decay-mixed-propagation}
Under the assumptions of Theorem~\ref{thm:decay}, the density satisfies,
uniformly for $T>0$,
\begin{equation}\label{decay-mixed-propagation}
 \|h\|_{\widetilde L_T^\infty(\mathcal X^{\frac 52,-\frac 12})}
 +\|h\|_{\widetilde L_T^2(\mathcal X^{\frac32,\frac12})}
 +\|h\|_{\widetilde L_T^1(\mathcal X^{\frac12,\frac32})}
 \leq C\|h_0\|_{\mathcal X^{\frac 52,-\frac 12}}.
\end{equation}
\end{lemma}

\begin{proof}
Integrating \eqref{dck} over $[0,t]$ and taking the supremum
over $t\in[0,T]$, we obtain
\[
 \begin{aligned}
 &\|\dot\Delta_j\dot\Delta_k^h h\|_{L_T^\infty(L^2)}
 +c2^{2(k-j)}\|\dot\Delta_j\dot\Delta_k^h h\|_{L_T^1(L^2)}\\
 &\quad\leq C\|\dot\Delta_j\dot\Delta_k^h h_0\|_{L^2}
 +C\|\dot\Delta_j\dot\Delta_k^hG(h)\|_{L_T^1(L^2)}\\
 &\qquad+C\|[U_1,\dot\Delta_j\dot\Delta_k^h]\partial_1h
                                      \|_{L_T^1(L^2)}+C\|\partial_1U_1\|_{L_T^1(L^\infty)}
              \|\dot\Delta_j\dot\Delta_k^h h\|_{L_T^\infty(L^2)}.
 \end{aligned}
\]
Then, one multiplies by $2^{\frac{5}{2}j-\frac{1}{2}k}$, takes the horizontal supremum and sums
over $j$. The dissipation acquires the weight $2^{\frac{1}{2}j+\frac{3}{2}k}$, so we obtain
\begin{equation}\label{decay-angular-energy}
 \begin{aligned}
 &\|h\|_{\widetilde L_T^\infty(\mathcal X^{\frac 52,-\frac 12})}
 +c\|h\|_{\widetilde L_T^1(\mathcal X^{\frac12,\frac32})}\\
 &\quad\leq C\|h_0\|_{\mathcal X^{\frac 52,-\frac 12}}
 +C\|G(h)\|_{\widetilde L_T^1(\mathcal X^{\frac 52,-\frac 12})}\\
 &\qquad+C\sum_j2^{\frac{5}{2}j}\sup_k2^{-\frac{k}{2}}
    \|[U_1,\dot\Delta_j\dot\Delta_k^h]\partial_1h\|_{L_T^1(L^2)}+C\|\partial_1U_1\|_{L_T^1(L^\infty)}
                 \|h\|_{\widetilde L_T^\infty(\mathcal X^{\frac 52,-\frac 12})}.
 \end{aligned}
\end{equation}

{
By \eqref{a1} and \eqref{dahi}, choosing $\delta_0$ sufficiently
small ensures $C\|\partial_1U_1\|_{L_T^1(L^\infty)}\leq\frac{1}{2}$.
We thus absorb the last term on the right-hand side of
\eqref{decay-angular-energy} into the first term on its left-hand side.
}

We estimate $G(h)$ in \eqref{decay-angular-energy}
by \eqref{decay-mixed-source}. 
Regarding the commutator term, as in \eqref{zero-mixed-commutator-time}, one applies
\eqref{decay-mixed-commutator} to obtain
\[
 \begin{aligned}
 \sum_j2^{\frac{5}{2}j}\sup_k2^{-\frac{k}{2}}
 \|[U_1,\dot\Delta_j\dot\Delta_k^h]\partial_1h\|_{L_T^1(L^2)}
 &\lesssim
 \|U_1\|_{\widetilde L_T^2(\cB^{\frac 32,\frac 12})}
 \|h\|_{\widetilde L_T^2(\mathcal X^{\frac32,\frac12})}\\
 &+
 \|\partial_1U_1\|_{L_T^1(\cB^{\frac 12,\frac 12})}
 \|h\|_{\widetilde L_T^\infty(\mathcal X^{\frac 52,-\frac 12})}\\
 &+
 \|U_1\|_{\widetilde L_T^2(\mathcal X^{\frac 52,-\frac 12})}
 \|h\|_{\widetilde L_T^2(\cB^{\frac 12,\frac 32})}.
 \end{aligned}
\]
{
Using \eqref{decay-critical-coefficients}, \eqref{zero-mixed-gradient-time} and \eqref{decay-mixed-velocity-time}, we obtain
}
\[
 \begin{aligned}
 \sum_j2^{\frac{5}{2}j}\sup_k2^{-\frac{k}{2}}
 \|[U_1,\dot\Delta_j\dot\Delta_k^h]\partial_1h\|_{L_T^1(L^2)}\lesssim\delta_0\bigl(
 \|h\|_{\widetilde L_T^\infty(\mathcal X^{\frac 52,-\frac 12})}
 +\|h\|_{\widetilde L_T^2(\mathcal X^{\frac32,\frac12})}\bigr).
 \end{aligned}
\]
Substituting the bounds for $G(h)$ and the commutator into the
summed block estimate, we obtain
\begin{equation}\label{decay-mixed-closure}
 \begin{aligned}
 &\|h\|_{\widetilde L_T^\infty(\mathcal X^{\frac 52,-\frac 12})}
 +c\|h\|_{\widetilde L_T^1(\mathcal X^{\frac12,\frac32})}\\
 &\quad\leq C\|h_0\|_{\mathcal X^{\frac 52,-\frac 12}}
 +C\delta_0\bigl(
 \|h\|_{\widetilde L_T^\infty(\mathcal X^{\frac 52,-\frac 12})}
 +\|h\|_{\widetilde L_T^2(\mathcal X^{\frac32,\frac12})}
 +\|h\|_{\widetilde L_T^1(\mathcal X^{\frac12,\frac32})}\bigr).
 \end{aligned}
\end{equation}
By Proposition~\ref{prop:real-interpolation}, we have
\[
 \|h\|_{\widetilde L_T^2(\mathcal X^{\frac32,\frac12})}
 \leq\|h\|_{\widetilde L_T^\infty(\mathcal X^{\frac 52,-\frac 12})}^{\frac12}
      \|h\|_{\widetilde L_T^1(\mathcal X^{\frac12,\frac32})}^{\frac12}.
\]
Using the smallness of $\delta_0$ and Young's inequality in \eqref{decay-mixed-closure}, one concludes
\eqref{decay-mixed-propagation}.
\end{proof}

\subsubsection{Decay in the critical anisotropic norms}

Under the assumptions of Theorem~\ref{thm:decay}, the following estimate holds for every $t\geq0$:
\begin{equation}\label{decay-mixed-h}
 \|h(t)\|_{\cB^{\frac 32,\frac 12}}
 \leq C\|h_0\|_{\mathcal X^{\frac 52,-\frac 12}}(1+t)^{-\frac12}.
\end{equation}

\begin{proof}[Proof of \eqref{decay-mixed-h}]
If $\|h_0\|_{\mathcal X^{\frac 52,-\frac 12}}=0$,
uniqueness implies $h=0$.
Otherwise, the differential inequality \eqref{ddthigh} and
\eqref{decay-mixed-propagation} imply
\[
 \frac{d}{dt}\|h(t)\|_{\cB^{\frac 32,\frac 12}}
 +\frac{c\|h(t)\|_{\cB^{\frac 32,\frac 12}}^3}
       {\|h_0\|_{\mathcal X^{\frac 52,-\frac 12}}^2}
 \leq0.
\]
Solving this differential inequality and using 
$\|h(0)\|_{\cB^{\frac 32,\frac 12}}
 \leq C\|h_0\|_{\mathcal X^{\frac 52,-\frac 12}}$,
we obtain
\[
 \|h(t)\|_{\cB^{\frac 32,\frac 12}}
 \leq
 \left(
 \|h_0\|_{\cB^{\frac 32,\frac 12}}^{-2}
 +\frac{2ct}{\|h_0\|_{\mathcal X^{\frac 52,-\frac 12}}^2}
 \right)^{-\frac12}
 \lesssim
 \|h_0\|_{\mathcal X^{\frac 52,-\frac 12}}(1+t)^{-\frac12}.
\]
This proves \eqref{decay-mixed-h}.
\end{proof}

\subsubsection{Dissipation and gradient decay}

Under the assumptions of Theorem~\ref{thm:decay}, the following estimates hold for every $t\geq0$:
\begin{equation}\label{velocity-high-window-decay}
 \|h(t)\|_{\cB^{-\frac 12,\frac 52}}
 \leq C_{h_0}(1+t)^{-\frac32},\qquad
 \|h(t)\|_{\cB^{\frac 12,\frac 32}}
 \leq C_{h_0}(1+t)^{-1}.
\end{equation}
Consequently, we have
\begin{equation}\label{velocity-gradient-decay}
 \begin{aligned}
 \|\nabla b(t)\|_{L^\infty}
   &\leq C\|h_0\|_{\mathcal X^{\frac 52,-\frac 12}}(1+t)^{-\frac12},\\
 \|\nabla u_1(t)\|_{L^\infty}&\leq C_{h_0}(1+t)^{-1},\\
 \|\nabla u_2(t)\|_{L^\infty}&\leq C_{h_0}(1+t)^{-\frac32}.
 \end{aligned}
\end{equation}
The constant $C_{h_0}$ is independent of $t$.

\begin{proof}[Proof of \eqref{velocity-high-window-decay} and \eqref{velocity-gradient-decay}]
Multiplying \eqref{dck} by $e^{c2^{2(k-j)}t}$ and integrating
from $a$ to $t$, we obtain
\begin{equation}\label{velocity-block-exponential}
 \begin{aligned}
 \|\dot\Delta_j\dot\Delta_k^h h(t)\|_{L^2}
 &\leq e^{-c2^{2(k-j)}(t-a)}
       \|\dot\Delta_j\dot\Delta_k^h h(a)\|_{L^2}\\
 &\quad+\int_a^t e^{-c2^{2(k-j)}(t-s)}
  \Bigl(\|\dot\Delta_j\dot\Delta_k^hG(h)(s)\|_{L^2}\Bigr.+\|[U_1(s),\dot\Delta_j\dot\Delta_k^h]
                         \partial_1h(s)\|_{L^2}\\
 &\qquad\qquad\Bigl.
 +C\|\partial_1U_1(s)\|_{L^\infty}
 \|\dot\Delta_j\dot\Delta_k^h h(s)\|_{L^2}\Bigr)\,ds.
 \end{aligned}
\end{equation}
For $G(h)$ and the transport contribution, we use \eqref{dcb},
\eqref{dci}, the embedding $\cB^{\frac 12,\frac 12}\hookrightarrow L^\infty$
and {Proposition~\ref{prop:inter}} to obtain
\begin{equation}\label{velocity-high-source-time}
 \begin{aligned}
 &\sum_{j,k}2^{\frac{3j}{2}+\frac{k}{2}}
       \|\dot\Delta_j\dot\Delta_k^hG(h)(s)\|_{L^2}\\
 &\quad
 +\|\partial_1U_1(s)\|_{L^\infty}\|h(s)\|_{\cB^{\frac 32,\frac 12}}\lesssim
 \|h(s)\|_{\cB^{\frac 32,\frac 12}}\|h(s)\|_{\cB^{-\frac 12,\frac 52}}.
 \end{aligned}
\end{equation}

We take $a=\frac{t}{2}$ in \eqref{velocity-block-exponential}, multiply
by $2^{-\frac{1}{2}j+\frac{5}{2}k}$ and sum over $j,k$.
For the initial term, we use \eqref{decay-damped-weight} with
$\tau=\frac{t}{2}$ to gain the rate $(1+t)^{-1}$.
We estimate the commutator integral by applying
Corollary~\ref{lem:decay-damped-commutator} with $f=U_1(s)$,
$g=h(s)$ and $\tau=t-s$, using \eqref{dch}--\eqref{dci} and
interpolation. For the other two terms, we use
\eqref{velocity-high-source-time} and \eqref{decay-damped-weight}
with $\tau=t-s$. Combining the three bounds, we obtain
\begin{equation}\label{decay-D-volterra}
 \begin{aligned}
 \|h(t)\|_{\cB^{-\frac 12,\frac 52}}
 &\lesssim \frac{\|h(\frac{t}{2})\|_{\cB^{\frac 32,\frac 12}}}{1+t}+\int_{\frac{t}{2}}^t
 \frac{\|h(s)\|_{\cB^{\frac 32,\frac 12}}
       \|h(s)\|_{\cB^{-\frac 12,\frac 52}}}{1+t-s}\,ds,
 \qquad t\geq2.
 \end{aligned}
\end{equation}
We now use a time-weighted estimate.
On nonzero blocks, $k\leq j+N_0$, so the ratio
$2^{2(k-j)}$ of the dissipative and density weights is bounded.
By \eqref{dahi}, we have
\[
 \|h(t)\|_{\cB^{-\frac 12,\frac 52}}
 \lesssim\|h(t)\|_{\cB^{\frac 32,\frac 12}}
 \leq C\delta_0.
\]
By continuity, we extend this bound to every $t\geq0$.
In particular, using \eqref{decay-mixed-h} on $[\frac{t}{2},t]$ and
$\int_{\frac{t}{2}}^t(1+t-s)^{-1}\,ds\leq\log(1+t)$ in
\eqref{decay-D-volterra}, we obtain, for $2\leq t\leq T$,
\[
 \begin{aligned}
 (1+t)^{\frac32}\|h(t)\|_{\cB^{-\frac 12,\frac 52}}\lesssim\|h_0\|_{\mathcal X^{\frac 52,-\frac 12}}
  +\|h_0\|_{\mathcal X^{\frac 52,-\frac 12}}
    \frac{\log(1+t)}{\sqrt{1+t}}
    \sup_{1\leq s\leq T}
    (1+s)^{\frac32}\|h(s)\|_{\cB^{-\frac 12,\frac 52}}.
 \end{aligned}
\]
Since the coefficient of the supremum tends to zero, we choose
$t_*\geq2$, depending on $h_0$ but not on $T$, so that it is at most
$\frac12$ for $t\geq t_*$. The weighted norm is bounded on $[0,t_*]$.
Taking the supremum over $t_*\leq t\leq T$ and absorbing the last term,
we obtain a bound independent of $T$. Together with the bound on
$[0,t_*]$, this proves the first estimate in
\eqref{velocity-high-window-decay}. 
For the second estimate, we interpolate with \eqref{decay-mixed-h} { using
Proposition~\ref{prop:inter}}:
\[
 \|h(t)\|_{\cB^{\frac 12,\frac 32}}^2
 \leq\|h(t)\|_{\cB^{\frac 32,\frac 12}}
       \|h(t)\|_{\cB^{-\frac 12,\frac 52}}
 \leq C_{h_0}(1+t)^{-2}.
\]

Using the chain rule and \eqref{dc4}, we transfer
\eqref{dc-full-gradient}--\eqref{dc-good-gradient} to the original variables:
\[
\begin{aligned}
\|\nabla b(t)\|_{L^\infty}\lesssim\|h(t)\|_{\cB^{\frac 32,\frac 12}},\qquad \|\nabla u_1(t)\|_{L^\infty}\lesssim\|h(t)\|_{\cB^{\frac 12,\frac 32}},\qquad \|\nabla u_2(t)\|_{L^\infty}\lesssim\|h(t)\|_{\cB^{-\frac 12,\frac 52}}.
\end{aligned}
\]
This proves \eqref{velocity-gradient-decay}.
\end{proof}

\subsection{$L^2$ decay of the velocity}

Using \eqref{L214} and \eqref{velocity-gradient-decay}, we prove
\begin{align}\label{uL2decay11}
 \|u(t)\|_{L^2}\leq C(1+t)^{-\frac{3}{4}},
 \qquad
 \|u_2(t)\|_{L^2}\leq C(1+t)^{-\frac{5}{4}}
\end{align}
for every $t\geq0$, where $C$ depends on the initial data
but is independent of $t$.

\begin{proof}[Proof of \eqref{uL2decay11}]
Applying $\partial_t+u\cdot\nabla$ to Darcy's law
$u=-\nabla\pi-be_2$, where $e_2=(0,1)^T$ and
$\pi$ is the modified pressure, and using
$(\partial_t+u\cdot\nabla)b=u_2$ from \eqref{eq:ipm3},
we obtain
\[
\begin{aligned}
(\partial_t+u\cdot\nabla)u
&=-\nabla\bigl((\partial_t+u\cdot\nabla)\pi\bigr)
+(\nabla u)^T\nabla\pi-u_2e_2\\
&=-\nabla\left((\partial_t+u\cdot\nabla)\pi
+\frac12|u|^2\right)
-b\nabla u_2-u_2e_2.
\end{aligned}
\]
Applying the Leray projector
\[
\mathbb P=
\begin{pmatrix}
-\cR_2^2&\cR_1\cR_2\\
\cR_1\cR_2&-\cR_1^2
\end{pmatrix}
\]
and using $u=-\mathbb P(be_2)$, which implies
$-\mathbb P(u_2e_2)=\cR_1^2u$, we obtain
\begin{equation}\label{zero-momentum}
\partial_tu-\cR_1^2u
=-\mathbb P(u\cdot\nabla u+b\nabla u_2).
\end{equation}

Taking the $L^2$ scalar product with $u$, using
$\operatorname{div}u=0$ and $\|\cR_1u\|_{L^2}=\|u_2\|_{L^2}$, one obtains
\[
 \begin{aligned}
 \frac12\frac d{dt}\|u\|_{L^2}^2+\|u_2\|_{L^2}^2
 &=-\int_{\mathbb R^2}b\,\nabla u_2\cdot u\,dx\\
 &\leq\|b\|_{L^2}\|\nabla u_2\|_{L^\infty}\|u\|_{L^2}
 \lesssim(1+t)^{-\frac{7}{4}}\|u\|_{L^2},
 \end{aligned}
\]
where we used the known decay estimates \eqref{L214} and \eqref{velocity-gradient-decay}.
Since
\[
 \|u\|_{L^2}^4=\|\cR_1b\|_{L^2}^4
 \leq\|b\|_{L^2}^2\|u_2\|_{L^2}^2,
\]
it follows that
\[
 \frac d{dt}\|u\|_{L^2}^2
 +c(1+t)^{\frac12}\|u\|_{L^2}^4
 \leq C(1+t)^{-\frac{7}{4}}\|u\|_{L^2}.
\]
Set $Y(t):=M(1+t)^{-\frac32}$. For $M$ sufficiently large, we have
$Y(0)\geq\|u(0)\|_{L^2}^2$ and
\[
 \begin{aligned}
 Y'(t)+c(1+t)^{\frac12}Y(t)^2
 &=\left(cM^2-\frac32M\right)(1+t)^{-\frac{5}{2}}\geq C(1+t)^{-\frac{7}{4}}Y(t)^{\frac12}.
 \end{aligned}
\]
The comparison principle therefore yields
\[
 \|u(t)\|_{L^2}^2\leq Y(t)=M(1+t)^{-\frac32}.
\]

We now establish the improved decay of $u_2$. Taking the $L^2$ inner product of the second component of
\eqref{zero-momentum} with $u_2$, we obtain
\[
 \frac12\frac{d}{dt}\|u_2\|_{L^2}^2
 +\|\mathcal R_1u_2\|_{L^2}^2
 \lesssim
 \|u\cdot\nabla u+b\nabla u_2\|_{L^2}
 \|\mathcal R_1u_2\|_{L^2}.
\]
Here we used the expression for $\mathbb P$ and the $L^2$
boundedness of the Riesz transforms. Young's inequality then gives
\[
 \frac{d}{dt}\|u_2\|_{L^2}^2
 +\|\mathcal R_1u_2\|_{L^2}^2
 \lesssim
 \|u\cdot\nabla u+b\nabla u_2\|_{L^2}^2.
\]
The source term satisfies
\[
 \begin{aligned}
 \|u\cdot\nabla u+b\nabla u_2\|_{L^2}
 &\leq\|u\|_{L^2}\|\nabla u\|_{L^\infty}
       +\|b\|_{L^2}\|\nabla u_2\|_{L^\infty}\lesssim(1+t)^{-\frac{7}{4}}.
 \end{aligned}
\]
Using
$\|u_2\|_{L^2}^4\leq\|u\|_{L^2}^2\|\cR_1u_2\|_{L^2}^2$,
it follows that
\[
 \frac d{dt}\|u_2\|_{L^2}^2
 +c(1+t)^{\frac{3}{2}}\|u_2\|_{L^2}^4
 \leq C(1+t)^{-\frac{7}{2}}.
\]
Set $Y_2(t):=M_2(1+t)^{-\frac{5}{2}}$. For $M_2$ sufficiently large, we have
$Y_2(0)\geq\|u_2(0)\|_{L^2}^2$ and
\[
 Y_2'(t)+c(1+t)^{\frac{3}{2}}Y_2(t)^2
 =\left(cM_2^2-\frac52M_2\right)(1+t)^{-\frac{7}{2}}
 \geq C(1+t)^{-\frac{7}{2}}.
\]
The comparison principle therefore yields
\[
 \|u_2(t)\|_{L^2}^2\leq Y_2(t)=M_2(1+t)^{-\frac{5}{2}}.
\]
Eventually, we end up with \eqref{uL2decay11}.
\end{proof}

\subsection{$H^2$ decay of the density and velocity}

We prove that, for every $t\geq0$,
\begin{equation}\label{velocity-component-besov-decay}
 \|b(t)\|_{\dot B^2_{2,1}}
 +(1+t)^{\frac12}\|u_1(t)\|_{\dot B^2_{2,1}}
 +(1+t)\|u_2(t)\|_{\dot B^2_{2,1}}
 \leq C_{h_0}(1+t)^{-\frac14}.
\end{equation}

To estimate the second derivatives, we use the following interpolation
inequality.

\begin{lemma}[Angular interpolation]\label{lem:decay-angular-interpolation}
For $m\in\{0,1,2\}$ and every $f$ with finite norms on the right,
the following interpolation estimate holds:
\begin{equation}\label{velocity-angular-interpolation}
 \|\cR_1^mf\|_{\dot B^2_{2,1}}
 \lesssim\|f\|_{\mathcal X^{\frac 52,-\frac 12}}^{\frac{1}{2m+2}}
 \|f\|_{\cB^{\frac{3}{2}-m,\frac{1}{2}+m}}^{\frac{2m+1}{2m+2}}.
\end{equation}
\end{lemma}

\begin{proof}
Fix $\eta>0$ and split the horizontal blocks according to
$2^{k-j}\leq\eta$ and $2^{k-j}>\eta$. For $2^{k-j}\leq\eta$, note that, due to the horizontal almost orthogonality,
\[
 \begin{aligned}
 &\sum_j2^{2j}
 \left(\sum_{k:\,2^{k-j}\leq\eta}
       2^{2m(k-j)}
       \|\dot\Delta_j\dot\Delta_k^h f\|_{L^2}^2\right)^{\frac12}\lesssim
 \left(\sum_{2^\ell\leq\eta}2^{(2m+1)\ell}\right)^{\frac12}
 \|f\|_{\mathcal X^{\frac52,-\frac12}}
 \lesssim
 \eta^{m+\frac{1}{2}}\|f\|_{\mathcal X^{\frac52,-\frac12}}.
 \end{aligned}
\]
For $2^{k-j}>\eta$, we sum the blocks in $\ell^1$ and use $ 2^{2j}2^{m(k-j)}
 =2^{-\frac{k-j}{2}}
  2^{(\frac32-m)j+(\frac12+m)k}$. This contribution is bounded by 
  $$
  \eta^{-\frac12}\|f\|_{\cB^{\frac{3}{2}-m,\,\frac{1}{2}+m}}.
  $$
  Choosing
\[
 \eta=
 \left(
 \frac{\|f\|_{\cB^{\frac{3}{2}-m,\,\frac{1}{2}+m}}}
      {\|f\|_{\mathcal X^{\frac52,-\frac12}}}
 \right)^{\frac{1}{m+1}}
\]
balances the two bounds and proves
\eqref{velocity-angular-interpolation}.
\end{proof}

\begin{proof}[Proof of \eqref{velocity-component-besov-decay}]
Applying Lemma~\ref{lem:decay-angular-interpolation} to $h(t)$, we
combine the uniform bound \eqref{decay-mixed-propagation} with
the decay estimates \eqref{decay-mixed-h} and
\eqref{velocity-high-window-decay} to obtain
\begin{equation}\label{velocity-riesz-besov-decay}
 \|\cR_1^mh(t)\|_{\dot B^2_{2,1}}
 \leq C_{h_0}(1+t)^{-\frac{2m+1}{4}},\qquad m=0,1,2.
\end{equation}

Using the div--curl system {\eqref{dc-errors}--\eqref{dc-divcurl}} and the
algebra property of $\dot B^1_{2,1}$, we obtain
\[
 \|\nabla U\|_{\dot B^1_{2,1}}
 \lesssim\|\cR_1h\|_{\dot B^2_{2,1}}
       +\|h\|_{\dot B^2_{2,1}}\|\nabla U\|_{\dot B^1_{2,1}}.
\]
Since $\|h\|_{\dot B^2_{2,1}}\leq C\delta_0$, we absorb the last term.
For $\partial_1U_1$ and $\nabla U_2$, we use \eqref{dc-errors}
and \eqref{dc-hodge} to obtain
\[
 \begin{aligned}
 &\|\partial_1U_1\|_{\dot B^1_{2,1}}
  +\|\nabla U_2\|_{\dot B^1_{2,1}}\\
 &\quad\lesssim\|\cR_1^2h\|_{\dot B^2_{2,1}}
  +\|h\|_{\dot B^2_{2,1}}
      \bigl(\|\partial_1U_1\|_{\dot B^1_{2,1}}
            +\|\nabla U_2\|_{\dot B^1_{2,1}}\bigr)
  +\|\cR_1h\|_{\dot B^2_{2,1}}\|\nabla U\|_{\dot B^1_{2,1}}.
 \end{aligned}
\]
We estimate $\partial_1h\,\partial_2U_1$ by the last term on the
right-hand side. Every other product in $F,Q$ contains $\partial_1U_1$
or a component of $\nabla U_2$.
Absorbing the term multiplied by $\|h\|_{\dot B^2_{2,1}}$ and
using the preceding bound for $\nabla U$, we obtain
\begin{equation}\label{velocity-isotropic-reconstruction}
 \|U_1\|_{\dot B^2_{2,1}}
 \lesssim\|\cR_1h\|_{\dot B^2_{2,1}},\qquad
 \|U_2\|_{\dot B^2_{2,1}}
 \lesssim\|\cR_1^2h\|_{\dot B^2_{2,1}}
          +\|\cR_1h\|_{\dot B^2_{2,1}}^2.
\end{equation}
By \eqref{velocity-riesz-besov-decay}, the quadratic term is bounded
by $C_{h_0}(1+t)^{-\frac{3}{2}}$. Applying the composition estimate
\eqref{dcn} to $b=h\circ Y$ and $u=U\circ Y$, we obtain
\eqref{velocity-component-besov-decay}.
\end{proof}

\begin{proof}[Proof of Theorem~\ref{thm:decay}]
Combining \eqref{L214}, \eqref{uL2decay11} and
\eqref{velocity-component-besov-decay} with
$\dot B^2_{2,1}\hookrightarrow\dot H^2$, we obtain
the $H^2$ decay estimates.
The gradient decay estimates are \eqref{velocity-gradient-decay}.
\end{proof}

{
\begin{proof}[Proof of Theorem~\ref{thm:ipm-higher-sobolev-intro}]
For $s>2$, the embedding
$H^s\hookrightarrow\dot B^1_{2,1}\cap\dot B^2_{2,1}$ and
Theorem~\ref{thm:ipm-higher-sobolev} give global well-posedness
and the uniform Sobolev bound \eqref{ipm-uniform-sobolev-intro}.

If, in addition, $b_0\in W^{4,1}$, set
$h_0=b_0\circ X(0,\cdot)$ as in Theorem~\ref{thm:decay}.
Since $\|\nabla b_0\|_{L^\infty}<\frac{1}{2}$ and
$W^{4,1}(\mathbb R^2)\hookrightarrow W^{2,\infty}(\mathbb R^2)$,
differentiating $h_0(y)=b_0(y_1,y_2+h_0(y))$ up to order four
and using the boundedness of the Jacobians gives $h_0\in W^{4,1}$.
By \eqref{decay-localized-data} in Lemma~\ref{lem:mixed-embeddings},
$h_0\in\mathcal X^{\frac12,-\frac12}\cap\mathcal X^{\frac52,-\frac12}$.
The decay estimates follow from Theorem~\ref{thm:decay}.
\end{proof}
}

\appendix

\section{Properties of anisotropic Besov spaces}\label{sec:appendixBesov}
We justify the frequency embeddings used in the main text, with
special attention to the radial and horizontal index shift.

\subsection{Bernstein-type inequalities and embeddings}

We prove the Bernstein estimates and record the embeddings that
follow from the support restriction $k\leq j+N_0$.
\begin{lemma}\label{lem:emb}
  The following statements hold:
\begin{enumerate}
    \item \textbf{Anisotropic Bernstein embedding:} Let $1\leq r_1\leq r_2\leq \infty,$ $1\leq p_1\leq p_2\leq \infty.$ Then, for any $s_1,s_2\in\R$, we have
    \begin{align*}
        \cB^{s_1,s_2}_{p_1,r_1}\hookrightarrow \cB^{s_1-(\frac{1}{p_1}-\frac{1}{p_2}),s_2-(d-1)(\frac{1}{p_1}-\frac{1}{p_2})}_{p_2,r_2}.
    \end{align*}
    \item \textbf{Anisotropic index shift}: For $s_1,s_2\in \R$,
    $\alpha \geq 0$ and $1\leq p,r\leq\infty$, we have
    \begin{equation}
    \label{eq:emlh}
    \cB^{s_1+\alpha,s_2-\alpha}_{p,r}\hookrightarrow\cB^{s_1,s_2}_{p,r}.
    \end{equation}

    \item \textbf{Embedding in Lebesgue space:} For any $(p_1,p_2)\in[1,\infty]^2$ with $p_1\leq p_2$, we have 
    \begin{align*}
    \cB^{\frac{1}{p_1}-\frac{1}{p_2},(d-1)(\frac{1}{p_1}-\frac{1}{p_2})}_{p_1,1}\hookrightarrow L^{p_2}.
    \end{align*}
\end{enumerate}
\end{lemma}
We next relate the anisotropic Besov spaces to isotropic Besov spaces.
\begin{lemma}[{\cite[Lemma~2.4]{LZCPAM}}]
Let $s_1,s_2,t_1,t_2\in\R$ satisfy
$t_2>0$ and $s_2>t_1+t_2>s_1$. If
$f\in\dot H^{s_1}(\R^d)\cap\dot H^{s_2}(\R^d)$, then
$f\in\cB^{t_1,t_2}(\R^d)$, and we have
    \begin{equation}
    \label{eq:embh}
    \|f\|_{\cB^{t_1,t_2}(\R^d)}\lesssim \|f\|_{\dot B^{t_1+t_2}_{2,1}(\R^d)}\lesssim \|f\|_{\dot H^{s_1}(\R^d)\cap \dot H^{s_2}(\R^d) }.
    \end{equation}
\end{lemma}

{
\begin{lemma}[Embeddings for the mixed spaces]
\label{lem:mixed-embeddings}
Let $s\in\mathbb R$ and $\sigma>0$.
The following continuous embeddings hold:
\begin{equation}\label{mixed-besov-embedding}
 \mathcal X^{s+\sigma,-\sigma}
 \hookrightarrow\dot B^s_{2,1}
 \hookrightarrow\cB^{s-\sigma,\sigma}
 \hookrightarrow\mathcal X^{s-\sigma,\sigma}.
\end{equation}
In particular, we have 
$\mathcal X^{s+\sigma,-\sigma}\hookrightarrow\dot H^s$. Moreover, the space $W^{4,1}(\mathbb R^2)$ is continuously embedded in
$\mathcal X^{\frac12,-\frac12}\cap\mathcal X^{\frac52,-\frac12}$.
\end{lemma}

\begin{proof}
Since the mixed blocks vanish for $k>j+N_0$, we have
\[
\begin{aligned}
 \|f\|_{\dot B^s_{2,1}}
 &\lesssim
 \sum_j2^{(s+\sigma)j}
 \sup_k\bigl(2^{-\sigma k}
       \|\dot\Delta_j\dot\Delta_k^h f\|_{L^2}\bigr)
 \sum_{k\leq j+N_0}2^{\sigma(k-j)}\lesssim \|f\|_{\mathcal X^{s+\sigma,-\sigma}},
\end{aligned}
\]
because $\sigma>0$. This proves the first embedding in
\eqref{mixed-besov-embedding}. The second follows from
\eqref{eq:embh}. The third follows by bounding the horizontal
supremum by the corresponding sum.
The embedding into $\dot H^s$ follows from
$\dot B^s_{2,1}\hookrightarrow\dot H^s$.

For $f\in W^{4,1}$, Fourier localization implies
\[
 \|\dot\Delta_j\dot\Delta_k^h f\|_{L^2}\lesssim2^{\frac{1}{2}(j+k)}
 \min\Bigl\{\|f\|_{L^1},
 2^{-4j}\sum_{|\alpha|=4}\|\partial^\alpha f\|_{L^1}\Bigr\}.
\]
Multiplying by $2^{sj-\frac{1}{2}k}$ and taking the supremum in $k$, we obtain
\[
\begin{aligned}
 \|f\|_{\mathcal X^{s,-\frac{1}{2}}}
 &\lesssim
 \sum_{j\leq0}2^{(s+\frac{1}{2})j}\|f\|_{L^1}
 +\sum_{j>0}2^{(s-\frac{7}{2})j}
       \sum_{|\alpha|=4}\|\partial^\alpha f\|_{L^1}.
\end{aligned}
\]
The sums converge for $s=\frac{1}{2}$ and $s=\frac{5}{2}$, so
\begin{equation}\label{decay-localized-data}
 \|f\|_{\mathcal X^{\frac12,-\frac12}}
 +\|f\|_{\mathcal X^{\frac52,-\frac12}}
 \lesssim
 \|f\|_{L^1}+\sum_{|\alpha|=4}\|\partial^\alpha f\|_{L^1}.
\end{equation}
This proves the embedding of $W^{4,1}(\mathbb R^2)$.
\end{proof}
}

We next record the embeddings involving the Riesz transforms that are
adapted to the equations considered here.
\begin{lemma} \label{lemma:Embeds}
Let $1\leq p\leq\infty$ and
$f\in\cB^{\frac{1}{p}+1,\frac{d-1}{p}}_{p,1}(\R^d)$.
Then $\nabla f$, $\nabla\cR_1^2f$, $\nabla\cR_1f$, $\d_1f$ and
$\d_2f$ belong to $L^\infty(\R^d)$. More precisely, we have
    \begin{equation}
          \label{eq:lip}
          \begin{aligned}
 \|\nabla f\|_{L^\infty(\R^d)}+\|\d_2f\|_{L^\infty(\R^d)}
 &\lesssim\|f\|_{\cB^{\frac1p+1,\frac{d-1}{p}}_{p,1}(\R^d)},\\
 \|\d_1f\|_{L^\infty(\R^d)}
 &\lesssim\|f\|_{\cB^{\frac1p,\frac{d-1}{p}+1}_{p,1}(\R^d)},\\
 \|\nabla\cR_1^2f\|_{L^\infty(\R^d)}
 &\lesssim\|f\|_{\cB^{\frac1p-1,\frac{d-1}{p}+2}_{p,1}(\R^d)},\\
 \|\nabla\cR_1f\|_{L^\infty(\R^d)}
 &\lesssim\|f\|_{\cB^{\frac1p,\frac{d-1}{p}+1}_{p,1}(\R^d)}.
          \end{aligned}
    \end{equation}

\end{lemma}
\begin{proof}
By the anisotropic Bernstein inequality in Lemma~\ref{lem:BIA}, we have
   \begin{align*}
       \|\nabla f\|_{L^\infty}&\lesssim \sum_{j,k\in \Z} 2^j\|\dD_{j}\dD^h_{k}f\|_{L^\infty}\lesssim \sum_{j,k\in\Z}
          2^{k\frac{d-1}{p}}2^{j(1+\frac1p)}
          \|\dD_j\dD_k^h f\|_{L^p}=\|f\|_{\cB^{\frac1p+1,\frac{d-1}{p}}_{p,1}}.
   \end{align*}
The same argument gives the estimates for $\nabla\cR_1^2f$,
$\nabla\cR_1f$ and $\d_1f$.
\end{proof}

\begin{lemma}\label{prop:esu}
Let $d=2$, $s_1,s_2\in\R$ and $1\leq p,r\leq\infty$.
For $f\in\cB^{s_1,s_2}_{p,r}$, let
$g=(-\cR_2f,\cR_1f)$. Then the following estimates hold:
\begin{equation}\label{eq:esul}
 \begin{aligned}
 \|\cR_1f\|_{\cB^{s_1,s_2}_{p,r}}
 &\sim\|f\|_{\cB^{s_1-1,s_2+1}_{p,r}}
 \lesssim\|f\|_{\cB^{s_1,s_2}_{p,r}},\\
 \|\cR_2f\|_{\cB^{s_1,s_2}_{p,r}}
 &\lesssim\|f\|_{\cB^{s_1,s_2}_{p,r}}.
 \end{aligned}
\end{equation}
In particular, we have
\begin{equation}
 \|g\|_{\cB^{s_1,s_2}_{p,r}}\lesssim
 \|f\|_{\cB^{s_1,s_2}_{p,r}}.
\end{equation}
\end{lemma}
\begin{proof}
We recall that $\cR_i=(-\Delta)^{-\frac12}\partial_{x_i}$. The vertical
Littlewood--Paley decomposition and the radial support of $\dD_j$ give,
for a sufficiently large $N_0$,
\begin{align*}
 2^{js_1}2^{ks_2}\|\dD_j\dD_k^h(\cR_2f)\|_{L^p}
 &\lesssim\sum_{q\in\Z}2^{js_1}2^{ks_2}
      \|\dD_q^v\dD_j\dD_k^h\cR_2f\|_{L^p}\\
 &\lesssim\sum_{\substack{q\in\Z\\q\leq j+N_0}}
      2^{q-j}2^{js_1}2^{ks_2}
      \|\dD_q^v\dD_j\dD_k^h f\|_{L^p}\\
 &\lesssim2^{js_1}2^{ks_2}\|\dD_j\dD_k^h f\|_{L^p}.
\end{align*}
For $\cR_1f$, we use the Bernstein inequalities in
Lemma~\ref{lem:BIA} and \eqref{eq:supp-mixed} to obtain
\[
 2^{js_1}2^{ks_2}\|\dD_j\dD_k^h(\cR_1f)\|_{L^p}
 \sim2^{js_1}2^{k(s_2+1)}
      \|\dD_j\dD_k^h(\Lambda^{-1}f)\|_{L^p}
 \sim2^{j(s_1-1)}2^{k(s_2+1)}
      \|\dD_j\dD_k^h f\|_{L^p}.
\]
We sum over $(j,k)\in\Z^2$ and use \eqref{eq:emlh} to prove
\eqref{eq:esul}.
\end{proof}

We now state a complex interpolation inequality in anisotropic spaces.
\begin{proposition}[{Complex interpolation}]\label{prop:inter}
Let $t_1,\sigma_1,t_2,\sigma_2\in\R$, $0<\theta<1$ and
$1\le p,r\le\infty$. Set $s_1=\theta t_1+(1-\theta)\sigma_1$ and 
$s_2=\theta t_2+(1-\theta)\sigma_2$. If $f\in\cB^{t_1,t_2}_{p,r}\cap\cB^{\sigma_1,\sigma_2}_{p,r}$,
then $f\in\cB^{s_1,s_2}_{p,r}$, and we have
\[
\|f\|_{\cB^{s_1,s_2}_{p,r}}
\lesssim
\|f\|_{\cB^{t_1,t_2}_{p,r}}^{\theta}
\|f\|_{\cB^{\sigma_1,\sigma_2}_{p,r}}^{1-\theta}.
\]
\end{proposition}
\begin{proof}
{The estimate follows from H\"older's inequality for the
weighted dyadic sequences when $r<\infty$, and from taking the supremum
over $(j,k)$ when $r=\infty$.}
\end{proof}

\begin{proposition}[Real interpolation]
\label{prop:real-interpolation}
Let $t_1,\sigma_1,t_2,\sigma_2\in\mathbb R$ and $0<\theta<1$,
with $t_2\neq\sigma_2$. Set
\[
 s_1=\theta t_1+(1-\theta)\sigma_1,
 \qquad
 s_2=\theta t_2+(1-\theta)\sigma_2.
\]
If $f\in\mathcal X^{t_1,t_2}\cap\mathcal X^{\sigma_1,\sigma_2}$,
then $f\in\cB^{s_1,s_2}$ and
\begin{equation}\label{eq:real-interpolation}
 \|f\|_{\cB^{s_1,s_2}}
 \lesssim
 \|f\|_{\mathcal X^{t_1,t_2}}^\theta
 \|f\|_{\mathcal X^{\sigma_1,\sigma_2}}^{1-\theta},
\end{equation}
where the mixed norm is defined in \eqref{def:mixed-X}.
The implicit constant depends only on $\theta$ and
$|t_2-\sigma_2|$. The estimate also holds with
$\|f\|_{\cB^{\sigma_1,\sigma_2}}$ in place of $\|f\|_{\mathcal X^{\sigma_1,\sigma_2}}$ in \eqref{eq:real-interpolation}.

For $T>0$, the mixed norms defined in \eqref{decay-mixed-time}
satisfy
\begin{align}\label{timedepinter}
 \|f\|_{\widetilde L_T^2(\mathcal X^{\frac32,\frac12})}
 \leq
 \|f\|_{\widetilde L_T^\infty(\mathcal X^{\frac52,-\frac12})}^{\frac12}
 \|f\|_{\widetilde L_T^1(\mathcal X^{\frac{1}{2},\frac{3}{2}})}^{\frac12}.
\end{align}
\end{proposition}
\begin{proof}
We prove \eqref{eq:real-interpolation} by splitting the horizontal
frequencies and estimating the two resulting sums separately.

Fix $\eta>0$. For each $j\in\mathbb Z$, we split the horizontal
indices according to
\[
 2^{(\sigma_1-t_1)j+(\sigma_2-t_2)k}\leq\eta
 \qquad\text{or}\qquad
 2^{(\sigma_1-t_1)j+(\sigma_2-t_2)k}>\eta.
\]

On the first region, the definitions of $s_1$ and $s_2$ imply
\[
 2^{js_1+ks_2}
 =
 2^{jt_1+kt_2}
 2^{(1-\theta)[(\sigma_1-t_1)j+(\sigma_2-t_2)k]}.
\]
We therefore have
\[
 \begin{aligned}
 &\sum_{\substack{k\in\mathbb Z\\
 2^{(\sigma_1-t_1)j+(\sigma_2-t_2)k}\leq\eta}}
 2^{js_1+ks_2}\|\dot\Delta_j\dot\Delta_k^h f\|_2\\
 &\quad\leq
 2^{jt_1}
 \sup_{k\in\mathbb Z}
 \bigl(2^{kt_2}\|\dot\Delta_j\dot\Delta_k^h f\|_2\bigr)
 \sum_{\substack{k\in\mathbb Z\\
 2^{(\sigma_1-t_1)j+(\sigma_2-t_2)k}\leq\eta}}
 2^{(1-\theta)[(\sigma_1-t_1)j+(\sigma_2-t_2)k]}.
 \end{aligned}
\]
For fixed $j$, we sum over integers $k$ bounded above if
$\sigma_2>t_2$ and bounded below if $\sigma_2<t_2$.
In all cases, the largest term is at most $\eta^{1-\theta}$,
and the remaining terms decrease with ratio
$2^{-(1-\theta)|\sigma_2-t_2|}<1$. Hence
\[
 \sum_{\substack{k\in\mathbb Z\\
 2^{(\sigma_1-t_1)j+(\sigma_2-t_2)k}\leq\eta}}
 2^{(1-\theta)[(\sigma_1-t_1)j+(\sigma_2-t_2)k]}
 \leq
 \frac{\eta^{1-\theta}}
      {1-2^{-(1-\theta)|\sigma_2-t_2|}}.
\]
This bound is independent of $j$. Summing the preceding estimate
over $j$ and using \eqref{def:mixed-X}, we obtain
\[
 \sum_{j\in\mathbb Z}
 \sum_{\substack{k\in\mathbb Z\\
 2^{(\sigma_1-t_1)j+(\sigma_2-t_2)k}\leq\eta}}
 2^{js_1+ks_2}\|\dot\Delta_j\dot\Delta_k^h f\|_2
 \lesssim
 \eta^{1-\theta}\|f\|_{\mathcal X^{t_1,t_2}}.
\]

On the second region, we instead use the identity
\[
 2^{js_1+ks_2}
 =
 2^{j\sigma_1+k\sigma_2}
 2^{-\theta[(\sigma_1-t_1)j+(\sigma_2-t_2)k]}.
\]
It follows that
\[
 \begin{aligned}
 &\sum_{\substack{k\in\mathbb Z\\
 2^{(\sigma_1-t_1)j+(\sigma_2-t_2)k}>\eta}}
 2^{js_1+ks_2}\|\dot\Delta_j\dot\Delta_k^h f\|_2\\
 &\quad\leq
 2^{j\sigma_1}
 \sup_{k\in\mathbb Z}
 \bigl(2^{k\sigma_2}\|\dot\Delta_j\dot\Delta_k^h f\|_2\bigr)
 \sum_{\substack{k\in\mathbb Z\\
 2^{(\sigma_1-t_1)j+(\sigma_2-t_2)k}>\eta}}
 2^{-\theta[(\sigma_1-t_1)j+(\sigma_2-t_2)k]}.
 \end{aligned}
\]
Here the largest term is at most $\eta^{-\theta}$, and the
geometric ratio is $2^{-\theta|\sigma_2-t_2|}<1$. Thus
\[
 \sum_{\substack{k\in\mathbb Z\\
 2^{(\sigma_1-t_1)j+(\sigma_2-t_2)k}>\eta}}
 2^{-\theta[(\sigma_1-t_1)j+(\sigma_2-t_2)k]}
 \leq
 \frac{\eta^{-\theta}}
      {1-2^{-\theta|\sigma_2-t_2|}}.
\]
Summing over $j$ now yields
\[
 \sum_{j\in\mathbb Z}
 \sum_{\substack{k\in\mathbb Z\\
 2^{(\sigma_1-t_1)j+(\sigma_2-t_2)k}>\eta}}
 2^{js_1+ks_2}\|\dot\Delta_j\dot\Delta_k^h f\|_2
 \lesssim
 \eta^{-\theta}\|f\|_{\mathcal X^{\sigma_1,\sigma_2}}.
\]

The two regions cover all horizontal indices. Adding their
contributions, we conclude that, for every $\eta>0$,
\[
 \|f\|_{\cB^{s_1,s_2}}
 \lesssim
 \eta^{1-\theta}\|f\|_{\mathcal X^{t_1,t_2}}
 +
 \eta^{-\theta}\|f\|_{\mathcal X^{\sigma_1,\sigma_2}}.
\]
If both mixed norms are positive, we choose
\[
 \eta=
 \frac{\|f\|_{\mathcal X^{\sigma_1,\sigma_2}}}
      {\|f\|_{\mathcal X^{t_1,t_2}}}.
\]
With this choice, we prove \eqref{eq:real-interpolation}.
If either mixed norm vanishes, every mixed block of $f$ vanishes,
so the estimate is immediate.

Finally, if
$\|f\|_{\cB^{\sigma_1,\sigma_2}}$ replaces the second mixed norm, then it follows from
\[
 \begin{aligned}
 \|f\|_{\mathcal X^{\sigma_1,\sigma_2}}
 &=
 \sum_j2^{j\sigma_1}
 \sup_k
 \bigl(2^{k\sigma_2}\|\dot\Delta_j\dot\Delta_k^h f\|_2\bigr)\leq
 \sum_{j,k}2^{j\sigma_1+k\sigma_2}
 \|\dot\Delta_j\dot\Delta_k^h f\|_2
 =
 \|f\|_{\cB^{\sigma_1,\sigma_2}}.
 \end{aligned}
\]
One can get the time-dependent estimate \eqref{timedepinter} by using H\"older's inequality on each
mixed block $(j,k)$. This completes the proof.
\end{proof}

\section{Proofs of the anisotropic product and commutator estimates}\label{app:com}
We give the frequency arguments for Propositions~\ref{AnisoPL}--\ref{lem:endpoint-product}
and Proposition~\ref{hcom}, together with the endpoint mixed
{commutator estimate} used in the decay proof.

\subsection{Proof of the product estimates}\label{app:product-proof}
We use the anisotropic version of Bony's
paraproduct decomposition from \cite{CPZ,LZCPAM}:
\[
 fg=T^{h}(f,g)+R^h(f,g)
 \andf
 fg=T(f,g)+R(f,g),
\]
where the horizontal paraproduct and remainder are defined by
\[
 T^{h}(f,g){:=}\sum_{k}\dS^h_{k-1}f\dD^h_k g,
 \qquad
 R^{h}(f,g){:=}\sum_{k}\dD^h_k f\dS^h_{k+2} g,
\]
and the radial paraproduct and remainder are defined by
\[
 T(f,g){:=}\sum_{j}\dS_{j-1}f\dD_j g,
 \qquad
 R(f,g){:=}\sum_{j}\dD_j f\dS_{j+2} g.
\]
Thus $R(f,g)$ and $R^h(f,g)$ include both the high--low and
comparable-frequency interactions. We use this notation
throughout the appendices.
Applying first the horizontal decomposition and then the radial one,
we obtain
 \[
     fg=TT^h(f,g)+TR^h(f,g)+RT^h(f,g)+RR^h(f,g),
 \]
where the four terms are defined by
 \begin{align*}
 TT^h(f,g)&{:=}\sum_{j,k\in\Z}
 \dS_{j-1}(\dS^h_{k-1}f)\dD_j(\dD^h_kg),
 &RT^h(f,g)&{:=}\sum_{j,k\in\Z}
 \dD_j(\dS^h_{k-1}f)\dS_{j+2}(\dD^h_kg),\\
 TR^h(f,g)&{:=}\sum_{j,k\in\Z}
 \dS_{j-1}(\dD^h_kf)\dD_j(\dS^h_{k+2}g),
 &RR^h(f,g)&{:=}\sum_{j,k\in\Z}
 \dD_j(\dD^h_kf)\dS_{j+2}(\dS^h_{k+2}g).
 \end{align*}

\begin{proof}[Proof of Proposition~\ref{AnisoPL}]
\noindent\textup{(i)} \emph{(Positive radial index).}
We apply the four-term decomposition above to $f$ and $g$.
We estimate the four terms separately. For the first one, we obtain
\begin{align*}
   \|\dD_j\dD_k^h TT^h(f,g)\|_{L^2}&\lesssim \sum_{|j'-j|\leq 4, |k'-k|\leq 4} \|\dS_{j'-1}\dS_{k'-1}^hf\dD_{j'}\dD^h_{k'}g\|_{L^2}
   \\&\lesssim  \sum_{|j'-j|\leq 4, |k'-k|\leq 4} \|\dS_{j'-1}\dS_{k'-1}^hf\|_{L^\infty}\|\dD_{j'}\dD^h_{k'}g\|_{L^2}
   \\& \lesssim c_{j,k}2^{-js_1}2^{-ks_2}\|f\|_{L^\infty}\|g\|_{\cB^{s_1,s_2}}.
\end{align*}
Here $(c_{j,k})_{(j,k)\in\Z^2}$ denotes a sequence with $\|\{c_{j,k}\}_{j,k\in\Z}\|_{\ell^1}\le 1$.

For the fourth term, \eqref{eq:orth} and $s_1,s_2>0$ give
\begin{equation}\label{app-positive-RR}
\begin{aligned}
  &\quad \|\dD_j\dD_k^h RR^h(f,g)\|_{L^2}\\
  &\lesssim \sum_{j'>j-N_0,k'>k-N_0 } \|\dD_{j'}\dD^h_{k'}f\dS_{j'+2}\dS^h_{k'+2}g\|_{L^2}
   \\&\lesssim 2^{-js_1}2^{-ks_2}\sum_{j'>j-N_0,k'>k-N_0} 2^{(j-j')s_1}2^{j's_1}2^{(k-k')s_2}2^{k's_2}\|\dD_{j'}\dD^h_{k'}f\|_{L^2}\|\dS_{j'+2}\dS^h_{k'+2}g\|_{L^\infty}
   \\&\lesssim  c_{j,k}2^{-js_1}2^{-ks_2}\|f\|_{\cB^{s_1,s_2}}\|g\|_{L^\infty}.
\end{aligned}
\end{equation}
For the third term, the anisotropic Bernstein inequality gives
\begin{equation}\label{app-positive-RT}
\begin{aligned}
&\qquad \|\dD_j\dD_k^h RT^h(f,g)\|_{L^2}\\
  &\lesssim \sum_{j'>j-N_0, |k'-k|\leq 4} \|\dD_{j'}\dS^h_{k'-1}f\dS_{j'+2}\dD^h_{k'}g\|_{L^2}
    \\&\lesssim 2^{-js_1}2^{-ks_2}\sum_{j'>j-N_0, |k'-k|\leq 4} 2^{(j-j')s_1}2^{j's_1} \|\dD_{j'}\dS^h_{k'-1}f\|_{L^2_vL^\infty_h}\\
    &\qquad\qquad\qquad\times 2^{(k-k')s_2}2^{k's_2}\|\dS_{j'+2}\dD^h_{k'}g\|_{L^\infty_vL^2_h}
    \\&\lesssim  c_{j,k}2^{-js_1}2^{-ks_2}\|f\|_{\cB^{s_1,\frac12}}\|g\|_{\cB^{\frac12,s_2}}.
\end{aligned}
\end{equation}
By interchanging the radial and horizontal frequency orderings, we obtain
\begin{equation*}
   \|\dD_j\dD_k^h TR^h(f,g)\|_{L^2}
 \lesssim  c_{j,k}2^{-js_1}2^{-ks_2}\|f\|_{\cB^{\frac12,s_2}}\|g\|_{\cB^{s_1,\frac12}}.
\end{equation*}
We sum the four estimates to prove the result.

\medskip\noindent\textup{(ii)} \emph{(Negative radial index).}
We prove \eqref{dmltapp} by decomposing the inputs into radial blocks and
separating the horizontal configurations only when necessary. It is enough to show that, for every $(j,k,\ell)$,
\begin{equation}\label{app-negative-pair}
\|(\dD_j\dD_k^h f)(\dD_\ell g)\|_{\cB^{-\frac12,\frac12}}
\lesssim
2^{-\frac{1}{2}j+\frac{1}{2}k}\|\dD_j\dD_k^h f\|_{L^2}
\,2^\ell\|\dD_\ell g\|_{L^2}.
\end{equation}
Indeed, summing this estimate over $(j,k,\ell)$ immediately yields
\eqref{dmltapp}.

For a sufficiently large fixed $C_0$, we distinguish three radial
interactions.

\smallskip\noindent\emph{Low frequencies: $\ell\leq j-C_0$.}
The product has radial frequency $2^j$. If moreover
$\ell\leq k-C_0$, then its horizontal frequency is also localized at
$2^k$. Hence, it follows by Bernstein's inequality that
\[
\|\dD_\ell g\|_{L^\infty}
\lesssim
2^\ell\|\dD_\ell g\|_{L^2},
\]
and therefore one gets
\[
\|(\dD_j\dD_k^h f)(\dD_\ell g)\|_{\cB^{-\frac12,\frac12}}
\lesssim
2^{-\frac{1}{2}j+\frac{1}{2}k}\|\dD_j\dD_k^h f\|_{L^2}
\,2^\ell\|\dD_\ell g\|_{L^2}.
\]

Assume now that $\ell>k-C_0$. Horizontal cancellation may occur, but
the product is supported in $|\xi_1|\lesssim2^\ell$. Mixed Bernstein
and H\"older inequalities lead to
\[
\|(\dD_j\dD_k^h f)(\dD_\ell g)\|_{L^2}
\lesssim
\|\dD_j\dD_k^h f\|_{L_h^\infty L_v^2}
\|\dD_\ell g\|_{L_h^2L_v^\infty}
\lesssim
2^{\frac{k+\ell}{2}}
\|\dD_j\dD_k^h f\|_{L^2}
\|\dD_\ell g\|_{L^2}.
\]
Since the product has radial frequency comparable to $2^j$ and its horizontal
blocks satisfy $q\leq\ell+C$, summing the horizontal weights yields $
\sum_{q\leq\ell+C}2^{\frac{q}{2}}\lesssim2^{\frac{1}{2}\ell}$.
Thus, we have 
\[
\|(\dD_j\dD_k^h f)(\dD_\ell g)\|_{\cB^{-\frac12,\frac12}}
\lesssim
2^{-\frac{1}{2}j}2^{\frac{1}{2}\ell}2^{\frac{k+\ell}{2}}
\|\dD_j\dD_k^h f\|_{L^2}
\|\dD_\ell g\|_{L^2},
\]
which is exactly \eqref{app-negative-pair}.

\smallskip\noindent\emph{High frequencies: $j\leq\ell-C_0$.}
The output has radial frequency $2^\ell$ and horizontal frequency at
most $C2^\ell$. Since $\sum_{q\leq\ell+C}2^{\frac{q-\ell}{2}}\lesssim1$, the $\cB^{-\frac12,\frac12}$ norm is controlled by the $L^2$ norm.
Using Bernstein, we have
\[
\begin{aligned}
\|(\dD_j\dD_k^h f)(\dD_\ell g)\|_{\cB^{-\frac12,\frac12}}
&\lesssim
2^{\frac{1}{2}(j+k)}
\|\dD_j\dD_k^h f\|_{L^2}
\|\dD_\ell g\|_{L^2}\lesssim
2^{-\frac{1}{2}j+\frac{1}{2}k}\|\dD_j\dD_k^h f\|_{L^2}
\,2^\ell\|\dD_\ell g\|_{L^2},
\end{aligned}
\]
where the last inequality follows from $j\leq\ell+C$.

\smallskip\noindent\emph{Intermediate frequencies: $|j-\ell|\leq C_0$.}
By Young's inequality and Plancherel,
\[
\|\widehat{(\dD_j\dD_k^h f)(\dD_\ell g)}
\|_{L^2_{\xi_1}L^\infty_{\xi_2}}
\lesssim
\|\widehat{\dD_j\dD_k^h f}\|_{L^1_{\xi_1}L^2_{\xi_2}}
\|\widehat{\dD_\ell g}\|_{L^2}
\lesssim
2^{\frac{k}{2}}
\|\dD_j\dD_k^h f\|_{L^2}
\|\dD_\ell g\|_{L^2}.
\]
After applying the radial output projector $\dD_m$, the Fourier support
is contained in $|\xi_2|\lesssim2^m$. Hence, for every horizontal
output block $q$,
\[
\|\dD_m\dD_q^h
((\dD_j\dD_k^h f)(\dD_\ell g))\|_{L^2}
\lesssim
2^{\frac{m+k}{2}}
\|\dD_j\dD_k^h f\|_{L^2}
\|\dD_\ell g\|_{L^2}.
\]
Since only $m\leq j+C$ and $q\leq m+C$ occur,
\[
\sum_{m\leq j+C}\sum_{q\leq m+C}
2^{-\frac{1}{2}m+\frac{1}{2}q}2^{\frac{m+k}{2}}
\lesssim
2^{\frac{1}{2}(j+k)}.
\]
Using $|j-\ell|\leq C_0$, we have $2^{\frac{1}{2}(j+k)}
\lesssim
2^{-\frac{1}{2}j+\frac{1}{2}k}2^\ell$, and \eqref{app-negative-pair} follows. Summing \eqref{app-negative-pair} over $(j,k,\ell)$ completes the proof
of \eqref{dmltapp}.
\end{proof}

\subsubsection{Endpoint mixed product estimate}
\label{app:endpoint-product}
\begin{proof}[Proof of Proposition~\ref{lem:endpoint-product}]
We use the same anisotropic Bony decomposition as above:
\[
fg=TT^h(f,g)+TR^h(f,g)+RT^h(f,g)+RR^h(f,g).
\]
Throughout the proof, $(c_j)_{j\in\Z}$ denotes a nonnegative
sequence, independent of $k$, satisfying $\sum_j c_j\leq1$.
This sequence may change from line to line.

We shall use the following consequence of horizontal Bernstein,
H\"older's inequality and vertical Bernstein:
\begin{equation}\label{decay-horizontal-low-output}
\begin{aligned}
\left\|
\dD_j\dD_k^h
\bigl(\dD_{j'}\dD_{k'}^h f\,
      \dD_{j''}\dD_{k''}^h g\bigr)
\right\|_{L^2}\lesssim
2^{\frac{k}{2}+\frac{\min(j',j'')}{2}}
\|\dD_{j'}\dD_{k'}^h f\|_{L^2}
\|\dD_{j''}\dD_{k''}^h g\|_{L^2}.
\end{aligned}
\end{equation}
Indeed, the horizontal output projector maps $L_h^1(L_v^2)$
to $L^2$ with norm bounded by $C2^{\frac{k}{2}}$, and vertical Bernstein
is applied to the input with the lower radial frequency.

For the first term, the embedding
$\cB^{\frac 12,\frac 12}\hookrightarrow L^\infty$ gives
\[
\begin{aligned}
\|\dD_j\dD_k^h TT^h(f,g)\|_{L^2}
&\lesssim
\sum_{\substack{|j'-j|\leq4\\|k'-k|\leq4}}
\|\dS_{j'-1}\dS_{k'-1}^h f\|_{L^\infty}
\|\dD_{j'}\dD_{k'}^h g\|_{L^2}\\
&\lesssim
c_j2^{{-sj}+\frac{1}{2}k}
\|f\|_{\cB^{\frac 12,\frac 12}}
\|g\|_{\mathcal X^{{s},-\frac 12}}.
\end{aligned}
\]

{{For $RT^h(f,g)$, by using $\dS_{k'-1}^h=\sum_{k''\leq k'-2}\dD_{k''}^h$
and applying the horizontal Bernstein inequality to $f$ and the vertical Bernstein inequality 
to $g$ (cf. Lemma~\ref{lem:BIA}), we have
\begin{equation}\label{app-mixed-product-RT-main}
\begin{aligned}
\|\dD_j\dD_k^h RT^h(f,g)\|_2
&\lesssim
\sum_{\substack{j'\geq j-N_0\\|k'-k|\leq4}}
\|\dD_{j'}\dS_{k'-1}^h f\|_{L_h^\infty(L_v^2)}
\|\dS_{j'+2}\dD_{k'}^h g\|_{L_h^2(L_v^\infty)}\\
&\lesssim
2^{\frac{k}{2}}d_{k}\|g\|_{\cB^{\frac12,\frac12}}
\sum_{j'\geq j-N_0}
\sup_{k''}2^{-\frac{1}{2}k''}
\|\dD_{j'}\dD_{k''}^h f\|_2\\
&\lesssim
c_j2^{-sj+\frac{1}{2}k}
\|f\|_{\mathcal X^{s,-\frac{1}{2}}}
\|g\|_{\cB^{\frac12,\frac12}},
\end{aligned}
\end{equation}
The above estimates hold based on the following facts.}}
{
Using $\dS_{k'-1}^h=\sum_{k''\leq k'-2}\dD_{k''}^h$,
horizontal Bernstein's inequality in Lemma~\ref{lem:BIA} implies
\[
\begin{aligned}
 \|\dD_{j'}\dS_{k'-1}^h f\|_{L_h^\infty(L_v^2)}
 &\lesssim
 \sum_{k''\leq k'-2}2^{\frac{1}{2}k''}
 \|\dD_{j'}\dD_{k''}^h f\|_{L^2}\lesssim
 \sum_{k''\leq k'-2}\sum_{\abs{\ell-k''}\leq 1}2^{\frac{1}{2}k''}
 \|\dD_{j'}\dD_{k''}^h \dD_{\ell}^h f\|_{L^2}\\
 &\lesssim
 2^{k'}\left(\sum_{k''\leq k'-2}2^{k''-k'}\right)
 \sup_{\ell\in\mathbb Z}2^{-\frac{1}{2}\ell}
 \|\dD_{j'}\dD_\ell^h f\|_{L^2}\\
 &\lesssim
 2^{k'}\sup_{\ell\in\mathbb Z}2^{-\frac{1}{2}\ell}
 \|\dD_{j'}\dD_\ell^h f\|_{L^2}.
\end{aligned}
\]
For the norm containing $g$, vertical Bernstein's inequality in Lemma~\ref{lem:BIA} implies
\[
\begin{aligned}
 \|\dS_{j'+2}\dD_{k'}^h g\|_{L_h^2(L_v^\infty)}
 &\lesssim
 \sum_{j''\leq j'+1}2^{\frac{1}{2}j''}
 \|\dD_{j''}\dD_{k'}^h g\|_{L^2}\\
 &=2^{-\frac{1}{2}k'}
 \sum_{j''\leq j'+1}2^{\frac{1}{2}j''+\frac{1}{2}k'}
 \|\dD_{j''}\dD_{k'}^h g\|_{L^2}\\
 &\leq2^{-\frac{1}{2}k'}d_{k'}\|g\|_{\cB^{\frac12,\frac12}}.
\end{aligned}
\]
Here $d_{k'}\geq0$ and $\sum_{k'\in\Z}d_{k'}\leq1$.
Substituting these bounds into the first line of
\eqref{app-mixed-product-RT-main} and using the fact $|k'-k|\leq4$,
we obtain the bound with $2^{\frac{k}{2}}$ and $\sum_{k\in\Z}d_k\lesssim1$.
}
Here and below, we use that
$\sum_{j\leq j'+N_0}2^{{s(j-j')}}\lesssim1$
for $s>0$.

For $TR^h(f,g)$, the frequency restrictions are
$|j'-j|\leq4$ and $k'\geq k-N_0$. By Lemma~\ref{lem:BIA}, one has
\[
\begin{aligned}
\|\dD_j\dD_k^h TR^h(f,g)\|_{L^2}
&\lesssim
2^{\frac{k}{2}}
\sum_{\substack{|j'-j|\leq4\\k'\geq k-N_0}}
\|\dS_{j'-1}\dD_{k'}^h f\|_{L_h^2(L_v^\infty)}
\|\dD_{j'}\dS_{k'+2}^h g\|_{L^2}\\
&\lesssim
2^{\frac{k}{2}}\|f\|_{\cB^{\frac 12,\frac 12}}
\sum_{|j'-j|\leq4}
\sup_{k''}
2^{-\frac{1}{2}k''}\|\dD_{j'}\dD_{k''}^h g\|_{L^2}\\
&\lesssim
c_j2^{{-sj}+\frac{1}{2}k}
\|f\|_{\cB^{\frac 12,\frac 12}}
\|g\|_{\mathcal X^{{s},-\frac 12}}.
\end{aligned}
\]
Here we used the vertical Bernstein inequality in Lemma~\ref{lem:BIA} and
\[
\|\dD_{j'}\dS_{k'+2}^h g\|_{L^2}
\lesssim
2^{\frac{1}{2}k'}\sup_{k''}
2^{-\frac{1}{2}k''}\|\dD_{j'}\dD_{k''}^h g\|_{L^2}.
\]
The resulting sum of
$2^{\frac{1}{2}j''+\frac{1}{2}k'}
\|\dD_{j''}\dD_{k'}^h f\|_{L^2}$
is bounded by $\|f\|_{\cB^{\frac 12,\frac 12}}$.

Finally, for $RR^h(f,g)$, we split
\[
\dS_{k'+2}^h
=
\dS_{k'-1}^h+Q_{k'},
\qquad
Q_{k'}
:=
\sum_{|k''-k'|\leq1}\dD_{k''}^h.
\]
Thus, it follows that
\[
\begin{aligned}
RR^h(f,g)
&=
\sum_{j',k'}
\dD_{j'}\dD_{k'}^h f\,
\dS_{j'+2}\dS_{k'-1}^h g+
\sum_{j',k'}
\dD_{j'}\dD_{k'}^h f\,
\dS_{j'+2}Q_{k'} g.
\end{aligned}
\]
For the first sum, the output satisfies
$j'\geq j-N_0$ and $|k'-k|\leq4$.
Estimating the low-frequency term in $L^\infty$ and using the same
radial summation as above, we obtain
\begin{equation}\label{app-mixed-product-RR-low}
\left\|
\dD_j\dD_k^h
\sum_{j',k'}
\dD_{j'}\dD_{k'}^h f\,
\dS_{j'+2}\dS_{k'-1}^h g
\right\|_{L^2}
\lesssim
c_j2^{{-sj}+\frac{1}{2}k}
\|f\|_{\mathcal X^{{s},-\frac 12}}
\|g\|_{\cB^{\frac 12,\frac 12}}.
\end{equation}
In the second sum, the product may contain arbitrarily low
horizontal frequencies. Using \eqref{decay-horizontal-low-output}, we get
\begin{equation}\label{app-mixed-product-RR-high}
\begin{aligned}
&\left\|
\dD_j\dD_k^h
\sum_{j',k'}
\dD_{j'}\dD_{k'}^h f\,
\dS_{j'+2}Q_{k'} g
\right\|_{L^2}\\
&\quad\lesssim
2^{\frac{k}{2}}
\sum_{\substack{j'\geq j-N_0\\k'\geq k-N_0}}
\|\dD_{j'}\dD_{k'}^h f\|_{L^2}
\sum_{\substack{j''\leq j'+1\\|k''-k'|\leq1}}
2^{\frac{1}{2}j''}\|\dD_{j''}\dD_{k''}^h g\|_{L^2}\\
&\quad\lesssim
2^{\frac{k}{2}}\|g\|_{\cB^{\frac 12,\frac 12}}
\sum_{j'\geq j-N_0}
\sup_{k'}2^{-\frac{1}{2}k'}
\|\dD_{j'}\dD_{k'}^h f\|_{L^2}\\
&\quad\lesssim
c_j2^{{-sj}+\frac{1}{2}k}
\|f\|_{\mathcal X^{{s},-\frac 12}}
\|g\|_{\cB^{\frac 12,\frac 12}}.
\end{aligned}
\end{equation}
Here $|k''-k'|\leq1$ allows us to use the horizontal
$\ell^1$ summation in the norm of $g$.

We now multiply the four
estimates by $2^{{sj}-\frac{1}{2}k}$, so that
$2^{-\frac{k}{2}}$ cancels $2^{\frac{k}{2}}$.
Taking the supremum over $k$ and summing over $j$, we obtain
\eqref{decay-mixed-product}.
\end{proof}

\subsection{Proof of the commutator estimates}\label{app:commutator-proof}
\begin{lemma}\label{lem:2.97}
Let {$d=2$ and let} $p_i,q_i,r_i,m_i,r\in[1,\infty]$ for $i=1,2$.
Let $f$ be a Lipschitz function such that
    $\nabla f\in L^{p_1}_hL^{q_1}_v$ and
    $\partial_1f\in L^{r_1}_hL^{m_1}_v$. For any
    $g\in L^{p_2}_hL^{q_2}_v\cap L^{r_2}_hL^{m_2}_v$, we have
\begin{equation}
      \label{eq:ce2}
        \| [f, \dD_j \dD^h_k]g\|_{L^r}\lesssim 2^{-j}\|\nabla f\|_{L^{p_1}_h L^{q_1}_v}\|\dD^h_kg\|_{L^{p_2}_h L^{q_2}_v}+2^{-k}\|\partial_1 f\|_{L^{r_1}_h L^{m_1}_v}\|g\|_{L^{r_2}_h L^{m_2}_v}   \end{equation}
        with $\dfrac{1}{p_1}+\dfrac{1}{p_2}= \dfrac{1}{q_1}+\dfrac{1}{q_2}=\dfrac{1}{r_1}+\dfrac{1}{r_2}=\dfrac{1}{m_1}+\dfrac{1}{m_2}=\dfrac{1}{r}.$
\end{lemma}

\begin{proof}
We use
\[
 [f,\dD_j\dD_k^h]g
 =\dD_j[f,\dD_k^h]g+[f,\dD_j]\dD_k^h g.
\]
Writing the horizontal commutator with its convolution kernel and
applying the first-order Taylor formula gives
\[
 \begin{aligned}
 |[f,\dD_k^h]g(x)|
 &\leq\int_0^1\int_{\R}
 2^k|\check\varphi^h(2^kz)|\,|z|\,
 |\partial_1f(x_1-\tau z,x_2)|\,|g(x_1-z,x_2)|\,dz\,d\tau.
 \end{aligned}
\]
Then, H\"older and Minkowski inequalities yield
\[
 \|[f,\dD_k^h]g\|_{L^r}
 \lesssim 2^{-k}
 \|\partial_1f\|_{L_h^{r_1}(L_v^{m_1})}
 \|g\|_{L_h^{r_2}(L_v^{m_2})}.
\]
Similarly, one can deal with the kernel estimate for $\dD_j$ as follows:
\[
 \|[f,\dD_j]\dD_k^h g\|_{L^r}
 \lesssim 2^{-j}
 \|\nabla f\|_{L_h^{p_1}(L_v^{q_1})}
 \|\dD_k^h g\|_{L_h^{p_2}(L_v^{q_2})}.
\]
Since $\dD_j$ is bounded on $L^r$, combining the two estimates
proves \eqref{eq:ce2}.
\end{proof}
\begin{remark}
Taking $p_1=q_1$, $p_2=q_2$, $r_1=m_1$ and $r_2=m_2$ in
\eqref{eq:ce2}, we obtain
\begin{equation}\label{eq:esr11}
 \| [f,\dD_j \dD^h_k ]g\|_{L^r}
 \lesssim 2^{-k}\|\partial_1 f\|_{L^{r_1}}\|g\|_{L^{r_2}}
 +2^{-j}\|\nabla f\|_{L^{p_1}}\|\dD^h_kg\|_{L^{p_2}},
\end{equation}
with $\frac{1}{r}=\frac{1}{r_1}+\frac{1}{r_2}
=\frac{1}{p_1}+\frac{1}{p_2}.$
\end{remark}

\begin{proof}[Proof of Proposition~\ref{hcom}]
With the notation of Section~\ref{app:product-proof}, write
$[f,\dD_j\dD_k^h]\partial_1g=\sum_{r=1}^4R^r_{j,k}$, where
\[
\begin{aligned}
 R^1_{j,k}&{:=}TT^h(f,\dD_j\dD_k^h\partial_1g)
              -\dD_j\dD_k^hTT^h(f,\partial_1g),\\
 R^2_{j,k}&{:=}RT^h(f,\dD_j\dD_k^h\partial_1g)
              -\dD_j\dD_k^hRT^h(f,\partial_1g),\\
 R^3_{j,k}&{:=}TR^h(f,\dD_j\dD_k^h\partial_1g)
              -\dD_j\dD_k^hTR^h(f,\partial_1g),\\
 R^4_{j,k}&{:=}RR^h(f,\dD_j\dD_k^h\partial_1g)
              -\dD_j\dD_k^hRR^h(f,\partial_1g).
\end{aligned}
\]
We will rely on the mixed Bernstein bound
\begin{equation}\label{eq:fui}
 \|f\|_{L_v^\infty(\dot B^s_{2,1})_h}
 \leq\sum_{j,k}2^{ks}\|\dD_j\dD_k^hf\|_{L_v^\infty L_h^2}
 \lesssim\|f\|_{\cB^{\frac{1}{2},s}},
\end{equation}
and the embeddings in Appendix~\ref{sec:appendixBesov}.

\smallskip\noindent\emph{The low--low term $R^1_{j,k}$.}
Almost orthogonality gives
\[
 R^1_{j,k}
 =\sum_{\substack{|j'-j|\leq N_0\\|k'-k|\leq N_0}}
 [\dS_{j'-1}\dS_{k'-1}^hf,\dD_j\dD_k^h]
 \dD_{j'}\dD_{k'}^h\partial_1g,
\]
where, throughout the proof, $N_0$ is a
fixed integer depending only on the dyadic cutoffs.
By \eqref{eq:ce2},
\[
\begin{aligned}
 \|R^1_{j,k}\|_{L^2}
 \lesssim\sum_{\substack{|j'-j|\leq N_0\\|k'-k|\leq N_0}}
 \Bigl(2^{k'-j}\|\nabla\dS_{j'-1}\dS_{k'-1}^hf\|_{L^\infty}
       +\|\partial_1\dS_{j'-1}\dS_{k'-1}^hf\|_{L^\infty}\Bigr)\times\|\dD_{j'}\dD_{k'}^hg\|_{L^2}.
\end{aligned}
\]
For the higher-order weight, use
$\cB^{\frac 12,\frac 12}\hookrightarrow L^\infty$ and
$\cB^{\frac 32,\frac 12}\hookrightarrow\cB^{\frac 12,\frac 32}$ to obtain
\[
\begin{aligned}
 \sum_{j,k}2^{\frac{3j}{2}+\frac{k}{2}}\|R^1_{j,k}\|_{L^2}
 \lesssim\|\partial_2f\|_{\cB^{\frac 12,\frac 12}}
             \|g\|_{\cB^{\frac 12,\frac 32}}+\|\partial_1f\|_{\cB^{\frac 12,\frac 12}}
             \|g\|_{\cB^{\frac 32,\frac 12}}.
\end{aligned}
\]
For the lower-order weight, instead use
\[
\begin{aligned}
 \|\partial_2\dS_{j'-1}\dS_{k'-1}^hf\|_{L^\infty}
 &\lesssim\|\partial_2f\|_{\dot B^1_{2,1}},\qquad \|\partial_1\dS_{j'-1}\dS_{k'-1}^hf\|_{L^\infty}
 \lesssim2^{j'}\|\partial_1f\|_{\cB^{-\frac 12,\frac 12}}.
\end{aligned}
\]

This yields
\[
\begin{aligned}
 \sum_{j,k}2^{\frac j2+\frac k2}\|R^1_{j,k}\|_{L^2}
 \lesssim{}&\|\partial_2f\|_{\dot B^1_{2,1}}
             \|g\|_{\cB^{-\frac 12,\frac 32}}+\|\partial_1f\|_{\cB^{-\frac 12,\frac 12}}
             \|g\|_{\cB^{\frac 32,\frac 12}}.
\end{aligned}
\]
By Lemma~\ref{prop:esu} and \eqref{eq:embh}, the last two norms of
$g$ are bounded by $\|\cR_1g\|_{\dot B^1_{2,1}}$ and
$\|g\|_{\dot B^2_{2,1}}$, respectively.

\smallskip\noindent\emph{The mixed term $R^3_{j,k}$.}
Here
\[
 R^3_{j,k}
 =\sum_{\substack{|j'-j|\leq N_0\\k'>k-N_0}}
 [\dS_{j'-1}\dD_{k'}^hf,\dD_j\dD_k^h]
 \dD_{j'}\dS_{k'+2}^h\partial_1g.
\]
We estimate the two products separately. Horizontal localization
allows us to replace $\dD_{k'}^hf$ by
$2^{-k'}\dD_{k'}^h\partial_1f$ in the norm, while the derivative
on $g$ costs $C2^{k'}$. Bernstein inequality \eqref{eq:fui} gives,
for $s=0,1$,
\begin{equation}\label{app-mixed-commutator-sum}
 \sum_{j,k}2^{(s+\frac{1}{2})j+\frac{1}{2}k}\|R^3_{j,k}\|_{L^2}
 \lesssim\|\partial_1f\|_{\cB^{s-\frac 12,\frac 12}}
          \|g\|_{\cB^{\frac 32,\frac 12}}.
\end{equation}
In the case $s=0$, summing the radial low frequencies of $f$ first
produces
\[
 \sum_{j''\leq j'-2}2^{j''-j'}
 \Bigl(2^{-\frac{1}{2}j''+\frac{1}{2}k'}
       \|\dD_{j''}\dD_{k'}^h\partial_1f\|_{L^2}\Bigr).
\]
The summability of $2^{j''-j'}$ allows us to place the extra radial
derivative on $g$. For $s=1$, the same sum is controlled directly by
$\|\partial_1f\|_{\cB^{\frac 12,\frac 12}}$.

\smallskip\noindent\emph{The radial remainders $R^2_{j,k}$ and $R^4_{j,k}$.}
For $r=2,4$, let $H^h=T^h,R^h$, respectively. Since the projectors
commute, both terms have the form
\begin{equation}\label{app-radial-remainders}
\begin{aligned}
 R^r_{j,k}
 =\sum_{j'>j-N_0}\Bigl(&H^h(\dD_{j'}f,
                    \dS_{j'+2}\dD_j\dD_k^h\partial_1g)-\dD_j\dD_k^h H^h(\dD_{j'}f,
                    \dS_{j'+2}\partial_1g)\Bigr).
\end{aligned}
\end{equation}

For the higher-order estimate, the bounds for $RT^h$ and $RR^h$
in Section~\ref{app:product-proof} { (see \eqref{app-positive-RR}
and \eqref{app-positive-RT})} give
\begin{equation}\label{app-unprojected-remainders}
\begin{aligned}
 &\|RT^h(f,\dD_j\dD_k^h\partial_1g)\|_{L^2}
  +\|RR^h(f,\dD_j\dD_k^h\partial_1g)\|_{L^2}\\
 &\qquad\lesssim
 2^{\frac{1}{2}j+k}\|\dD_j\dD_k^hg\|_{L^2}
 \sum_{j'>j-N_0}\sum_{k'}
 2^{\frac{1}{2}k'}\|\dD_{j'}\dD_{k'}^hf\|_{L^2}.
\end{aligned}
\end{equation}
Since $j'>j-N_0$, we have $ 2^{\frac{3j}{2}+\frac{k}{2}}2^{\frac{1}{2}j+k}2^{\frac{1}{2}k'}
 \lesssim
 2^{\frac{1}{2}j+\frac{3}{2}k}2^{\frac{3}{2}j'+\frac{1}{2}k'}$. Multiplying \eqref{app-unprojected-remainders} by
$2^{\frac{3j}{2}+\frac{k}{2}}$ and summing over $j,k$, we therefore obtain
\[
 \sum_{j,k}2^{\frac{3j}{2}+\frac{k}{2}}
 (\|R^2_{j,k}\|_{L^2}+\|R^4_{j,k}\|_{L^2})
 \lesssim\|f\|_{\cB^{\frac 32,\frac 12}}\|g\|_{\cB^{\frac 12,\frac 32}}.
\]
Using
$\|f\|_{\cB^{\frac 32,\frac 12}}
\lesssim\|\partial_1f\|_{\cB^{\frac 12,\frac 12}}
       +\|\partial_2f\|_{\cB^{\frac 12,\frac 12}}$
and \eqref{eq:emlh}, we obtain \eqref{hc1}.

For the lower-order estimate, we first prove
\begin{equation}\label{app-low-remainders}
 \sum_{j,k}2^{\frac j2+\frac k2}
 (\|R^2_{j,k}\|_{L^2}+\|R^4_{j,k}\|_{L^2})
 \lesssim\|f\|_{\dot B^2_{2,1}}\|g\|_{\cB^{-\frac 12,\frac 32}}.
\end{equation}
For either choice of $H^h$, Bernstein gives
\[
 \|H^h(\dD_{j'}f,G)\|_{L^2}
 \lesssim\sum_{k'\leq j'+N_0}
 2^{\frac{j'+k'}{2}}\|\dD_{j'}\dD_{k'}^hf\|_{L^2}\|G\|_{L^2}
 \lesssim2^{j'}\|\dD_{j'}f\|_{L^2}\|G\|_{L^2}.
\]
Thus the first product in \eqref{app-radial-remainders}, after
multiplication by $2^{\frac j2+\frac k2}$, is bounded by
\[
 2^{j-j'}\Bigl(2^{2j'}\|\dD_{j'}f\|_{L^2}\Bigr)
 \Bigl(2^{-\frac{1}{2}j+\frac{3}{2}k}\|\dD_j\dD_k^hg\|_{L^2}\Bigr),
\]
which is summable over $j'>j-N_0$ and $(j,k)$.

For the second product, use that
$\dS_{j'+2}=\sum_{j''\leq j'+1}\dD_{j''}$. The horizontal product
bounds and H\"older's inequality in the vertical variable give rise to
\[
\begin{aligned}
 &\sum_k2^{\frac{k}{2}}
 \|\dD_k^hH^h(\dD_{j'}f,\dD_{j''}\partial_1g)\|_{L_v^1(L_h^2)}\lesssim
 2^{\frac{1}{2}j'}\|\dD_{j'}f\|_{L^2}
 \sum_{k''}2^{\frac{3k''}{2}}\|\dD_{j''}\dD_{k''}^hg\|_{L^2}.
\end{aligned}
\]
For $R^h$, we sum over $k$ using
$\sum_{k\leq k'+N_0}2^{\frac{k-k'}{2}}\lesssim1$.
We then apply Bernstein's inequality in the vertical variable
to $\dD_j$ and multiply by $2^{\frac{j}{2}}$ to obtain
\[
 2^{j-j'}2^{\frac{j''-j'}{2}}
 \Bigl(2^{2j'}\|\dD_{j'}f\|_{L^2}\Bigr)
 \sum_{k''}2^{-\frac{1}{2}j''+\frac{3k''}{2}}
        \|\dD_{j''}\dD_{k''}^hg\|_{L^2}.
\]
Summing the estimates for both terms in
\eqref{app-radial-remainders}, using
$\sum_{j<j'+N_0}2^{j-j'}\lesssim1$ and
$\sum_{j''\leq j'+1}2^{\frac{j''-j'}{2}}\lesssim1$,
we obtain \eqref{app-low-remainders}.

The terms involving $R^2_{j,k}$ and $R^4_{j,k}$ can be addressed by
\begin{equation}\label{app-low-remainders-complementary}
 \sum_{j,k}2^{\frac j2+\frac k2}
 (\|R^2_{j,k}\|_{L^2}+\|R^4_{j,k}\|_{L^2})
 \lesssim\|f\|_{\cB^{\frac 12,\frac 12}}\|g\|_{\cB^{\frac 32,\frac 12}}.
\end{equation}
Indeed, for the unprojected products put the blocks of $f$ in
$L^\infty$. For the projected products use the $RT^h$ and $RR^h$
bounds{ \eqref{app-positive-RR} and \eqref{app-positive-RT}}
with $s_1=s_2=\frac{1}{2}$. In both cases,
$\|\partial_1g\|_{\cB^{\frac 12,\frac 12}}
\lesssim\|g\|_{\cB^{\frac 32,\frac 12}}$.

Finally, we decompose
\[
 f=-\partial_1(-\Delta)^{-1}\partial_1f
   -\partial_2(-\Delta)^{-1}\partial_2f.
\]
Apply \eqref{app-low-remainders-complementary} to the first term
and \eqref{app-low-remainders} to the second, using
\[
\begin{aligned}
 \|\partial_1(-\Delta)^{-1}\partial_1f\|_{\cB^{\frac 12,\frac 12}}
 &\lesssim\|\partial_1f\|_{\cB^{-\frac 12,\frac 12}},\quad  \|\partial_2(-\Delta)^{-1}\partial_2f\|_{\dot B^2_{2,1}}\lesssim\|\partial_2f\|_{\dot B^1_{2,1}}.
\end{aligned}
\]
The embeddings used for $R^1_{j,k}$, together with
\eqref{app-mixed-commutator-sum}, now lead to \eqref{hc2}.
\end{proof}

\begin{proof}[{Proof of Corollary~\ref{lem:decay-damped-commutator}}]
For $k>j+N_0$, the operator $\dot\Delta_j\dot\Delta_k^h$ and
its commutator with $f$ are zero. On the remaining indices,
$2^{2(k-j)}$ is bounded above, so
\begin{equation}\label{decay-damped-weight}
\begin{aligned}
 2^{-\frac{1}{2}j+\frac{5}{2}k}e^{-c2^{2(k-j)}\tau}
 &=2^{\frac{3j}{2}+\frac{k}{2}}\,2^{2(k-j)}e^{-c2^{2(k-j)}\tau}\lesssim\frac{2^{\frac{3j}{2}+\frac{k}{2}}}{1+\tau}.
\end{aligned}
\end{equation}
We multiply each commutator norm by this weight, sum over
$j,k$ and apply \eqref{hc1} to obtain \eqref{decay-damped-commutator}.
\end{proof}

\subsubsection{Endpoint mixed commutator estimate}
\label{app:endpoint-commutator}
\begin{proof}[Proof of Proposition~\ref{lem:endpoint-commutator}]
We retain the four terms $R^r_{j,k}$ and the frequency conventions
of the proof of Proposition~\ref{hcom}.
Only the output weight and horizontal summation change.

For $R^1_{j,k}$, \eqref{eq:ce2} and Bernstein give
\begin{equation}\label{endpoint-low-low-kernel}
\begin{aligned}
 \|R^1_{j,k}\|_{L^2}
 \lesssim\sum_{\substack{|j'-j|\leq N_0\\|k'-k|\leq N_0}}
 &\Bigl(2^{k'-j}\|f\|_{\cB^{\frac 32,\frac 12}}
       +\|\partial_1f\|_{\cB^{\frac 12,\frac 12}}\Bigr)\|\dD_{j'}\dD_{k'}^hg\|_{L^2}.
\end{aligned}
\end{equation}
After multiplication by $2^{{sj}-\frac{1}{2}k}$, the weights on $g$ become
$2^{{(s-1)j'}+\frac{k'}{2}}$ and $2^{{sj'}-\frac{1}{2}k'}$, respectively. These give
the first two terms on the right of \eqref{decay-mixed-commutator}.

For $R^3_{j,k}$, the same two-product decomposition uses
\[
\begin{aligned}
 \|\dS_{j'-1}\dD_{k'}^hf\|_{L_h^2(L_v^\infty)}
 &\lesssim2^{-k'}\sum_{j''\leq j'-2}
 2^{\frac{1}{2}j''}\|\dD_{j''}\dD_{k'}^h\partial_1f\|_{L^2},\\
 \|\dD_{j'}\dS_{k'+2}^h\partial_1g\|_{L^2}
 &\lesssim2^{\frac{3k'}{2}}\sup_{k''}
 2^{-\frac{1}{2}k''}\|\dD_{j'}\dD_{k''}^hg\|_{L^2}.
\end{aligned}
\]
When $\dD_j\dD_k^h$ acts on the whole product, we apply
horizontal Bernstein from $L_h^1(L_v^2)$ to $L^2$ and use
the two estimates above. When it acts only on $g$, we employ
\[
 \|\dD_j\dD_k^h\dD_{j'}\dS_{k'+2}^h\partial_1g
       \|_{L_h^\infty(L_v^2)}
 \lesssim
 2^{2k}\sup_{k''}2^{-\frac{1}{2}k''}
       \|\dD_{j'}\dD_{k''}^hg\|_{L^2}.
\]
Since $k'>k-N_0$, we have
$2^{-k'}2^{2k}\lesssim2^{\frac{1}{2}k+\frac{1}{2}k'}$.
Thus, in both cases, summing over $j''$ and $k'$ yields
\[
 \|R^3_{j,k}\|_{L^2}
 \lesssim
 2^{\frac{k}{2}}\|\partial_1f\|_{\cB^{\frac12,\frac12}}
 \sum_{|j'-j|\leq N_0}\sup_{k''}
 2^{-\frac{1}{2}k''}\|\dD_{j'}\dD_{k''}^hg\|_{L^2}.
\]
Multiplying by $2^{sj-\frac{1}{2}k}$, taking the supremum over $k$
and summing over $j$, we obtain
\[
 \sum_j2^{sj}\sup_k2^{-\frac{k}{2}}\|R^3_{j,k}\|_{L^2}
 \lesssim
 \|\partial_1f\|_{\cB^{\frac12,\frac12}}
 \|g\|_{\mathcal X^{s,-\frac12}},
\]
which is the second term on the right-hand side of
\eqref{decay-mixed-commutator}.

For $R^2_{j,k}$ and $R^4_{j,k}$, we reuse
\eqref{app-radial-remainders}. The unprojected products are
estimated by \eqref{app-unprojected-remainders}. After multiplication
by $2^{sj-\frac{1}{2}k}$, their sum is bounded by
\[
 \Bigl(\sum_{j'>j-N_0}2^{\frac{3}{2}(j-j')}
       \sum_{k'}2^{\frac{3}{2}j'+\frac{1}{2}k'}
       \|\dD_{j'}\dD_{k'}^hf\|_{L^2}\Bigr)
 2^{{(s-1)j}+\frac{1}{2}k}\|\dD_j\dD_k^hg\|_{L^2}.
\]
Taking the horizontal supremum and summing in $j$ gives
$\|f\|_{\cB^{\frac 32,\frac 12}}\|g\|_{\mathcal X^{{s-1},\frac 12}}$.

For the projected products, the estimates for $RT^h$ and $RR^h$
in the proof of Proposition~\ref{lem:endpoint-product}{
(see {\eqref{app-mixed-product-RT-main}}, \eqref{app-mixed-product-RR-low}
and \eqref{app-mixed-product-RR-high})}, applied to
$(f,\partial_1g)$, give
\[
\begin{aligned}
 \|\dD_j\dD_k^hRT^h(f,\partial_1g)\|_{L^2}
  +\|\dD_j\dD_k^hRR^h(f,\partial_1g)\|_{L^2}\lesssim2^{\frac{k}{2}}\|g\|_{\cB^{\frac 12,\frac 32}}
 \sum_{j'>j-N_0}\sup_{k'}
 2^{-\frac{1}{2}k'}\|\dD_{j'}\dD_{k'}^hf\|_{L^2}.
\end{aligned}
\]

Here the summation uses
$\sum_{j<j'+N_0}2^{{s(j-j')}}\lesssim1$, yielding
$\|f\|_{\mathcal X^{{s},-\frac 12}}\|g\|_{\cB^{\frac 12,\frac 32}}$.
Combining the four terms proves \eqref{decay-mixed-commutator}.
\end{proof}

\section{Technical estimates in isopycnal coordinates}
\label{app:dc-tech}

This appendix contains the technical estimates used in
Section~\ref{dcsec}.
The isopycnal coordinate maps $X$ and $Y$ are defined in
Lemma~\ref{lem:density-reduction}.

\subsection{Composition by the isopycnal coordinates}
\label{app:dc-composition}
We control composition by the coordinate map and its inverse in
the isotropic spaces used in the main estimates.
\begin{proof}[Proof of Lemma~\ref{dcmp}]
 Let $\Phi$
be a bi-Lipschitz diffeomorphism with Jacobian bounded above and
below, and assume that
\begin{equation}\label{aD}
 D\Phi-I\in\dot B^1_{2,1} \quad \text{and} \quad D\Phi^{-1}-I\in\dot B^1_{2,1}.
\end{equation}
The change-of-variables formula and the chain rule show that
composition with $\Phi$ is bounded on $L^2$ and $\dot H^1$.
By real interpolation, using
$(L^2,\dot H^1)_{\sigma,q}=\dot B^\sigma_{2,q}$, we obtain
\[
 \|f\circ\Phi\|_{\dot B^\sigma_{2,q}}
 \lesssim \|f\|_{\dot B^\sigma_{2,q}},
 \qquad 0<\sigma<1,\quad q\in\{1,\infty\}.
\]
To extend this estimate to negative regularity, we use duality.
The adjoint of the composition operator is given by
\[
g\mapsto (g\circ\Phi^{-1})|\det D\Phi^{-1}|.
\]
Since $|\det D\Phi^{-1}|-1\in\dot B^1_{2,1}$,
the product law
$\dot B^1_{2,1}\cdot\dot B^\sigma_{2,\infty}
\hookrightarrow\dot B^\sigma_{2,\infty}$,
combined with the composition estimate for $\Phi^{-1}$,
shows that the adjoint is bounded on $\dot B^\sigma_{2,\infty}$.
By duality, composition with $\Phi$ is therefore bounded on
$\dot B^{-\sigma}_{2,1}$.
Finally, interpolating between $\dot B^{-\sigma}_{2,1}$
and $\dot B^\sigma_{2,1}$ yields the desired bound
in $\dot B^0_{2,1}$.
The chain rules
\[
\begin{aligned}
 \nabla(f\circ\Phi)
 &=(D\Phi)^T(\nabla f\circ\Phi),\\
 D^2(f\circ\Phi)
 &=(D\Phi)^T(D^2f\circ\Phi)D\Phi
   +\sum_{\ell=1}^2(\partial_\ell f\circ\Phi)D^2\Phi_\ell
\end{aligned}
\]
provide the bounds at indices one and two, using
$\dot B^1_{2,1}\cdot\dot B^0_{2,1}\hookrightarrow\dot B^0_{2,1}$
and $D^2\Phi\in\dot B^0_{2,1}$. Thus
\begin{equation}\label{dc-composition-standard}
 \|f\circ\Phi\|_{\dot B^s_{2,1}}
 \lesssim\|f\|_{\dot B^s_{2,1}},\qquad s=0,1,2.
\end{equation}
The constants are uniform when the bi-Lipschitz constants and the
two norms in \eqref{aD} are bounded.

We now check \eqref{aD} for $X$ and $Y$. Since
$Y(x)=(x_1,x_2-b(x))$, we have
\[
 DY-I=-\begin{pmatrix}0&0\\\partial_1b&\partial_2b\end{pmatrix},
 \qquad
 \|DY-I\|_{\dot B^1_{2,1}}\lesssim\|b\|_{\dot B^2_{2,1}}.
\]
For the inverse, we set
\begin{equation}\label{Ytau}
 Y_\tau(x){:=}(x_1,x_2-\tau b(x)),\qquad
 X_\tau{:=}Y_\tau^{-1},\qquad 0\leq\tau\leq1.
\end{equation}
The bound $\|\nabla b\|_{L^\infty}\lesssim\|b\|_{\dot B^2_{2,1}}$
ensures uniform bi-Lipschitz control. Differentiating
$Y_\tau\circ X_\tau=\operatorname{Id}$ gives
\[
 DX_\tau-I
 =\tau\left(\begin{pmatrix}0&0\\\partial_1b&\partial_2b\end{pmatrix}
       \circ X_\tau\right)DX_\tau.
\]
As long as $\|DX_\tau-I\|_{\dot B^1_{2,1}}\leq1$,
\eqref{dc-composition-standard} and the algebra property imply
\[
 \|DX_\tau-I\|_{\dot B^1_{2,1}}
 \lesssim\|b\|_{\dot B^2_{2,1}}
       \bigl(1+\|DX_\tau-I\|_{\dot B^1_{2,1}}\bigr).
\]
Smallness and continuity in $\tau$, starting from $DX_0=I$, yield
\[
 \sup_{0\leq\tau\leq1}\|DX_\tau-I\|_{\dot B^1_{2,1}}
 \lesssim\|b\|_{\dot B^2_{2,1}}.
\]
Thus both $X$ and $Y$ satisfy \eqref{aD}.
Applying \eqref{dc-composition-standard} to $X$ and $Y$
proves \eqref{dcn}. 
Due to $h=b\circ X$ and $b=h\circ Y$, it also gives rise to
$\|h\|_{\dot B^s_{2,1}}\simeq\|b\|_{\dot B^s_{2,1}}$ for $s=1,2$.
In particular, one has
\begin{equation}\label{dcm}
 \|h(0)\|_{\dot B^s_{2,1}}
 \lesssim\|b_0\|_{\dot B^s_{2,1}},\qquad s=1,2,
\end{equation}
which, together with \eqref{dca}, proves \eqref{dc3}.
For the converse statement, start from
$X(y)=(y_1,y_2+h(y))$ and apply the same argument to $X$ and its
inverse $Y$, interchanging $b$ and $h$.
\end{proof}

\begin{remark}
The argument above does not extend directly to the anisotropic space
$\cB^{\frac 32,\frac 12}$, because the change of variables mixes horizontal and
vertical derivatives. Indeed, we have
\[
\partial_1(f\circ Y)
=
(\partial_1f)\circ Y
-
\partial_1b\,(\partial_2f)\circ Y.
\]
Then under $f\in\cB^{\frac 32,\frac 12}$, it holds that
\[
\partial_1f\in\cB^{\frac{3}{2},-\frac{1}{2}},
\qquad
\partial_2f\in\cB^{\frac 12,\frac 12}.
\]
Thus the horizontal chain rule couples derivatives with different
anisotropic regularities. The preceding isotropic argument does not
establish boundedness of composition on $\cB^{\frac 32,\frac 12}$.
\end{remark}

\subsection{The Riesz transform in isopycnal coordinates}
\label{app:dc-riesz}

\begin{proof}[Proof of Lemma~\ref{dctechnical}]
\noindent\emph{(i) The transformed div--curl system.} 
First, we show that the coordinate error $G(h)$ satisfies
\begin{equation}\label{dc-G-equation}
 \Delta G(h)=\partial_2F+\partial_1Q,
\end{equation}
where $F$ and $Q$ are defined by
\begin{equation}\label{dc-errorsc8}
\begin{aligned}
 &F:=-\partial_2h\,\partial_1U_1
   +\partial_1h\,\partial_2U_1,\qquad Q:=-\partial_2h\,\partial_1U_2
   +\partial_1h\,\partial_2U_2.
\end{aligned}
\end{equation}
By \eqref{dc4} and $h=b\circ X$,
\[
 (\partial_{x_1}b)\circ X
 =\partial_1h-\frac{\partial_1h}{J}\partial_2h
 =\frac{\partial_1h}{J},
\]
where we used $J=1+\partial_2h$.
Applying \eqref{dc4} to $\div_xu=0$ and
$\curl_xu=-\partial_{x_1}b$, we obtain
\begin{equation}\label{obtain}
\begin{aligned}
 J\partial_1U_1-\partial_1h\,\partial_2U_1+\partial_2U_2
 &=J\bigl((\div_xu)\circ X\bigr)=0,\\
 J\partial_1U_2-\partial_1h\,\partial_2U_2-\partial_2U_1
 &=J\bigl((\curl_xu)\circ X\bigr)=-J\bigl((\partial_{x_1}b)\circ X\bigr)
 =-\partial_1h.
\end{aligned}
\end{equation}
Substituting $J=1+\partial_2h$ into \eqref{obtain}, we obtain
\begin{equation}\label{dc-errors}
\begin{aligned}
 \partial_1U_1+\partial_2U_2
 &=-\partial_2h\,\partial_1U_1
   +\partial_1h\,\partial_2U_1=F,\\
 \partial_1U_2-\partial_2U_1+\partial_1h
 &=-\partial_2h\,\partial_1U_2
   +\partial_1h\,\partial_2U_2=Q.
\end{aligned}
\end{equation}
Thus, with $F$ and $Q$ defined by \eqref{dc-errors},
\begin{equation}\label{dc-divcurl}
 \div U=F,\qquad \curl U=-\partial_1h+Q.
\end{equation}
Using \eqref{dc-divcurl}, we obtain
\begin{equation}\label{dc-hodge}
\begin{aligned}
 \Delta U_2
 &=\partial_2\div U+\partial_1\curl U
 =-\partial_1^2h+\partial_2F+\partial_1Q,\qquad \partial_1U_1=F-\partial_2U_2.
\end{aligned}
\end{equation}
Finally, $G(h)=U_2-\cR_1^2h$ by \eqref{eq:defG}.
Combining \eqref{dc-hodge} with
$\Delta(\cR_1^2h)=-\partial_1^2h$, we obtain
\[
 \Delta G(h)
 =\Delta U_2-\Delta(\cR_1^2h)
 =\partial_2F+\partial_1Q,
\]
which proves \eqref{dc-G-equation}.
We use \eqref{dc-hodge} to estimate $U_2$ and $\partial_1U_1$.

\medskip\noindent\emph{(ii) Control of the full gradient.}
We reconstruct $\nabla U$ from its divergence and curl with
order-zero multipliers and apply Proposition~\ref{AnisoPL} (i),
with $s_1=s_2=\frac{1}{2}$, to obtain
\[
 \|\nabla U\|_{\cB^{\frac12,\frac12}}
 \lesssim \|h\|_{\cB^{\frac12,\frac32}}
 +\|F\|_{\cB^{\frac12,\frac12}}
 +\|Q\|_{\cB^{\frac12,\frac12}}
 \lesssim \|h\|_{\cB^{\frac12,\frac32}}
 +\|h\|_{\cB^{\frac32,\frac12}}
   \|\nabla U\|_{\cB^{\frac12,\frac12}}.
\]
If $\|h\|_{\cB^{\frac32,\frac12}}$ is sufficiently small, we absorb the
last term and obtain
\begin{equation}\label{dc-full-gradient}
 \|\nabla U\|_{\cB^{\frac12,\frac12}}
 \lesssim \|h\|_{\cB^{\frac12,\frac32}},
\end{equation}
which proves \eqref{dch}.

\medskip\noindent\emph{(iii) Recovering the extra horizontal derivative.}
In \eqref{dc-errors}, the only term involving $\partial_2U_1$
is $\partial_1h\,\partial_2U_1$. By
\eqref{app-critical-algebra} and \eqref{dc-full-gradient},
\[
 \|\partial_1h\,\partial_2U_1\|_{\cB^{\frac12,\frac12}}
 \lesssim
 \|\partial_1h\|_{\cB^{\frac12,\frac12}}
 \|\partial_2U_1\|_{\cB^{\frac12,\frac12}}
 \lesssim \|h\|_{\cB^{\frac12,\frac32}}^2.
\]
Applying \eqref{app-critical-algebra} to the remaining terms
in \eqref{dc-errors}, we obtain
\[
 \begin{aligned}
 \|F\|_{\cB^{\frac12,\frac12}}
 +\|Q\|_{\cB^{\frac12,\frac12}}
 \lesssim{}&\|h\|_{\cB^{\frac32,\frac12}}
 \bigl(\|\partial_1U_1\|_{\cB^{\frac12,\frac12}}
       +\|\nabla U_2\|_{\cB^{\frac12,\frac12}}\bigr)
 +\|h\|_{\cB^{\frac12,\frac32}}^2.
 \end{aligned}
\]
Combining this estimate with \eqref{dc-hodge} and
$\|\nabla\cR_1^2h\|_{\cB^{\frac12,\frac12}}
\lesssim\|h\|_{\cB^{-\frac12,\frac52}}$ gives
\begin{equation}\label{eq:u1u2}
 \begin{aligned}
 \|\partial_1U_1\|_{\cB^{\frac12,\frac12}}
 +\|\nabla U_2\|_{\cB^{\frac12,\frac12}}
 \lesssim{}&\|h\|_{\cB^{-\frac12,\frac52}}\\
 &+\|h\|_{\cB^{\frac32,\frac12}}
 \bigl(\|\partial_1U_1\|_{\cB^{\frac12,\frac12}}
       +\|\nabla U_2\|_{\cB^{\frac12,\frac12}}\bigr)
 +\|h\|_{\cB^{\frac12,\frac32}}^2.
 \end{aligned}
\end{equation}
We absorb the second term on the right-hand side of \eqref{eq:u1u2}
and use the interpolation inequality
\[
 \|h\|_{\cB^{\frac12,\frac32}}^2
 \lesssim\|h\|_{\cB^{\frac32,\frac12}}
 \|h\|_{\cB^{-\frac12,\frac52}}
\]
to obtain
\begin{equation}\label{dc-good-gradient}
 \|\partial_1U_1\|_{\cB^{\frac12,\frac12}}
 +\|\nabla U_2\|_{\cB^{\frac12,\frac12}}
 \lesssim \|h\|_{\cB^{-\frac12,\frac52}}.
\end{equation}
This proves \eqref{dci}. Finally, using \eqref{dc-G-equation}
and the product bounds above, we obtain
\[
 \begin{aligned}
 \|G(h)\|_{\cB^{\frac 32,\frac 12}}
 &\simeq\|\nabla G(h)\|_{\cB^{\frac 12,\frac 12}}\\
 &\lesssim\|F\|_{\cB^{\frac 12,\frac 12}}+\|Q\|_{\cB^{\frac 12,\frac 12}}\\
 &\lesssim \|h\|_{\cB^{\frac 12,\frac 32}}^2
 +\|h\|_{\cB^{\frac 32,\frac 12}}\|h\|_{\cB^{-\frac 12,\frac 52}}.
 \end{aligned}
\]
This is \eqref{dcb}. The embedding
$\cB^{\frac 12,\frac 12}\hookrightarrow L^\infty$, together with
$U_2=\cR_1^2h+G(h)$, leads to \eqref{dcu2lip}.
\end{proof}

\subsection{Isotropic estimates for $G(h)$} 
\label{app:dc-isotropic}
In this section, we prove the product estimates used in the isotropic energy inequalities of Section~\ref{dcsec}.
\begin{proof}[{Proof of Lemma~\ref{diso}}]
It follows from the embedding \eqref{eq:emlh} and
Proposition~\ref{prop:inter} with $\theta=\frac{1}{2}$ that
\begin{equation}\label{dc-isotropic-interpolation}
 \|h\|_{\cB^{-\frac 12,\frac 52}}
 \lesssim\|h\|_{\cB^{\frac 32,\frac 12}},\qquad
 \|h\|_{\cB^{\frac 12,\frac 32}}^2
 \lesssim\|h\|_{\cB^{\frac 32,\frac 12}}\|h\|_{\cB^{-\frac 12,\frac 52}}.
\end{equation}
By \eqref{dc-full-gradient} and
\eqref{dc-good-gradient}, we have
\[
 \|\nabla U\|_{L^\infty}\lesssim \|h\|_{\cB^{\frac 12,\frac 32}},\qquad
 \|\partial_1U_1\|_{L^\infty}+\|\nabla U_2\|_{L^\infty}\lesssim \|h\|_{\cB^{-\frac 12,\frac 52}}.
\]

We use the classical Bony decomposition $fg=T(f,g)+R(f,g)$ from
Section~\ref{app:product-proof}. For $i=1,2$, set
\[
 \begin{aligned}
 P_i&{:=}R(U_i,\partial_2h),\qquad
 N_i{:=}R(U_i,\partial_1h),\qquad E_i{:=}-T(\partial_1U_i,\partial_2h)
       +T(\partial_2U_i,\partial_1h).
 \end{aligned}
\]
Using the divergence form of \eqref{dc-errors}, we obtain
\begin{align}\label{FQ}
 F=-\partial_1P_1+\partial_2N_1+E_1,\qquad
 Q=-\partial_1P_2+\partial_2N_2+E_2.
\end{align}
Indeed, $R(U_i,\partial_\ell h)=U_i\partial_\ell h
-T(U_i,\partial_\ell h)$, and the second derivatives of $h$ cancel
when differentiating the two paraproducts.

We first establish the estimate for $U$ in \eqref{dgt}.
By \eqref{dc-divcurl}, we have
\begin{equation}\label{dc-U-reconstruction}
 U=\begin{pmatrix}-\cR_1\cR_2h\\\cR_1^2h\end{pmatrix}
 +\Delta^{-1}
  \begin{pmatrix}\partial_1F-\partial_2Q\\\partial_2F+\partial_1Q\end{pmatrix}.
\end{equation}
After substituting \eqref{FQ} into \eqref{dc-U-reconstruction},
the operators acting on $P_i,N_i$ have order $0$, and those
acting on $E_i$ have order $-1$.

The standard paraproduct estimates
\cite[Theorems~2.47 and~2.52]{BCD}, together with
$L^\infty\hookrightarrow\dot B^0_{\infty,\infty}$ and
$\cB^{\frac 12,\frac 12}\hookrightarrow L^\infty$, yield
\begin{equation}\label{dc-isotropic-PNE}
 \begin{aligned}
 \sum_i(\|P_i\|_{\dot B^s_{2,1}}+\|N_i\|_{\dot B^s_{2,1}})
 &\lesssim \|h\|_{\cB^{\frac 32,\frac 12}}\|U\|_{\dot B^s_{2,1}},\\
 \sum_i\|E_i\|_{\dot B^{s-1}_{2,1}}
 &\lesssim\|\nabla U\|_{L^\infty}\|h\|_{\dot B^s_{2,1}}
 \lesssim \|h\|_{\cB^{\frac 12,\frac 32}}\|h\|_{\dot B^s_{2,1}}.
 \end{aligned}
\end{equation}
Consequently, for $s=1,2$, we have
\begin{equation}\label{dc-isotropic-velocity}
 \begin{aligned}
 \|U\|_{\dot B^s_{2,1}}\lesssim{}&
 \|\cR_1h\|_{\dot B^s_{2,1}}
 +\|h\|_{\cB^{\frac 32,\frac 12}}\|U\|_{\dot B^s_{2,1}}
 +\|h\|_{\cB^{\frac 12,\frac 32}}\|h\|_{\dot B^s_{2,1}}.
 \end{aligned}
\end{equation}
After absorbing the second term on the right-hand side of
\eqref{dc-isotropic-velocity}, we obtain the estimate for $U$
in \eqref{dgt}.

Define
\begin{equation}
\label{def:g1g2}
    \begin{aligned}
 G_1&{:=}\cR_2P_1+\cR_1P_2,\qquad G_2{:=}-\cR_2^2N_1-\cR_1\cR_2N_2
       +\Delta^{-1}(\partial_2E_1+\partial_1E_2).
 \end{aligned}
\end{equation}
Using \eqref{dc-G-equation} and \eqref{FQ}, we obtain
\begin{equation}\label{dc-G-decomposition}
 G=\cR_1G_1+G_2.
\end{equation}
By \eqref{def:g1g2}, the estimate for $G_1$ follows from the
bound for $P_i$ in \eqref{dc-isotropic-PNE}.

For $G_2$, we retain the horizontal derivative on $h$ whenever the
coefficient is $\partial_2U_1$. Since
$\|\partial_1h\|_{\dot B^{s-1}_{2,1}}\simeq
\|\cR_1h\|_{\dot B^s_{2,1}}$, one has
\[
 \begin{aligned}
 \|E_1\|_{\dot B^{s-1}_{2,1}}
 &\lesssim \|h\|_{\cB^{-\frac 12,\frac 52}}\|h\|_{\dot B^s_{2,1}}
       +\|h\|_{\cB^{\frac 12,\frac 32}}\|\cR_1h\|_{\dot B^s_{2,1}},\\
 \|E_2\|_{\dot B^{s-1}_{2,1}}
 &\lesssim \|h\|_{\cB^{-\frac 12,\frac 52}}\|h\|_{\dot B^s_{2,1}},\\
 \sum_i\|N_i\|_{\dot B^s_{2,1}}&\lesssim \|h\|_{\cB^{\frac 12,\frac 32}}\|U\|_{\dot B^s_{2,1}}.
 \end{aligned}
\]
The estimate for $U$ in \eqref{dgt} and the interpolation
inequality in \eqref{dc-isotropic-interpolation} yield
\begin{equation}\label{dc-decomposition-estimates}
 \begin{aligned}
 \|G_1\|_{\dot B^s_{2,1}}
 &\lesssim \|h\|_{\cB^{\frac 32,\frac 12}}\|\cR_1h\|_{\dot B^s_{2,1}}
       +\|h\|_{\cB^{\frac 32,\frac 12}} \|h\|_{\cB^{\frac 12,\frac 32}}\|h\|_{\dot B^s_{2,1}},\\
 \|G_2\|_{\dot B^s_{2,1}}
 &\lesssim \|h\|_{\cB^{\frac 12,\frac 32}}\|\cR_1h\|_{\dot B^s_{2,1}}
       +\|h\|_{\cB^{-\frac 12,\frac 52}}\|h\|_{\dot B^s_{2,1}}.
 \end{aligned}
\end{equation}
These prove the first two bounds in \eqref{dgt}.
Using \eqref{dc-isotropic-interpolation}, \eqref{dc-G-decomposition}
and \eqref{dc-decomposition-estimates}, together with the smallness of
$\|h\|_{\cB^{\frac 32,\frac 12}}$, we also obtain \eqref{dm13}.

For the time-integrated estimates, we apply H\"older's inequality to each dyadic
block in the bounds used for \eqref{dc-isotropic-PNE},
as in \cite[Section~2.6.3]{BCD}.
\end{proof}

\subsection{A low-order estimate for $G(h)$}
\label{app:dc-loworder}
\begin{proof}[Proof of Lemma~\ref{lem:dc-loworder}]
Applying \eqref{dm12} with $s=2$ and using
$\|h\|_{\cB^{\frac 12,\frac 32}}
\lesssim\|\cR_1h\|_{\dot B^2_{2,1}}$, we obtain
\[
 \|\partial_2U_1\|_{\dot B^1_{2,1}}
 \lesssim\|\cR_1h\|_{\dot B^2_{2,1}},
\]
provided that $\|h\|_{\dot B^2_{2,1}}$ is sufficiently small.
The formulas \eqref{dc-errors}, the transformed div--curl identities
\eqref{dc-hodge} and Proposition~\ref{AnisoPL}\textup{(ii)} lead to
\begin{align}
 &\|(\partial_1U_1,\nabla U_2)\|_{\cB^{-\frac 12,\frac 12}}
  \lesssim\|h\|_{\cB^{-\frac{3}{2},\frac{5}{2}}}+\|F\|_{\cB^{-\frac 12,\frac 12}}+\|Q\|_{\cB^{-\frac 12,\frac 12}},\\
  &\|(F,Q)\|_{\cB^{-\frac 12,\frac 12}}\lesssim\|h\|_{\dot B^2_{2,1}}
   \bigl(\|\partial_1U_1\|_{\cB^{-\frac 12,\frac 12}}
        +\|\nabla U_2\|_{\cB^{-\frac 12,\frac 12}}\bigr)+\|\partial_1h\|_{\cB^{-\frac 12,\frac 12}}
               \|\partial_2U_1\|_{\dot B^1_{2,1}}.\label{dc-low-FQ-energy}
\end{align}
Since $\|\partial_1h\|_{\cB^{-\frac 12,\frac 12}}
\lesssim\|\cR_1h\|_{\dot B^1_{2,1}}$, smallness allows us to absorb
the velocity-gradient terms and yields
\begin{equation}\label{dc-low-gradient}
\begin{aligned}
 \|\partial_1U_1\|_{\cB^{-\frac 12,\frac 12}}
  +\|\nabla U_2\|_{\cB^{-\frac 12,\frac 12}}\lesssim\|h\|_{\cB^{-\frac{3}{2},\frac{5}{2}}}
        +\|\cR_1h\|_{\dot B^1_{2,1}}\|\cR_1h\|_{\dot B^2_{2,1}}.
\end{aligned}
\end{equation}
Substituting \eqref{dc-low-gradient} into
\eqref{dc-low-FQ-energy}, we have
\[
\begin{aligned}
 \|F\|_{\cB^{-\frac 12,\frac 12}}+\|Q\|_{\cB^{-\frac 12,\frac 12}}\lesssim\|h\|_{\dot B^2_{2,1}}\|h\|_{\cB^{-\frac{3}{2},\frac{5}{2}}}
        +\|\cR_1h\|_{\dot B^1_{2,1}}\|\cR_1h\|_{\dot B^2_{2,1}}.
\end{aligned}
\]
Together with \eqref{dc-G-equation}, this proves \eqref{dglo}.
Applying \eqref{hc2} with $f=U_1$ and $g=h$ then controls the
low-order commutator in \eqref{dllo}. The divergence term is bounded
by $C\|\partial_1U_1\|_{L^\infty}\|h\|_{\cB^{\frac 12,\frac 12}}$.
\end{proof}

\subsection{Mixed estimates for $G(h)$}
\label{app:mixed-source}
\begin{proof}[Proof of Lemma~\ref{lem:mixed-source-estimates}]
We first prove \eqref{zero-mixed-source}.
The time-dependent versions of \eqref{dc-full-gradient} and
\eqref{dc-good-gradient}, together with
\eqref{decay-critical-coefficients}, imply
\begin{equation}\label{zero-mixed-gradient-time}
 \begin{aligned}
 &\|U_1\|_{\widetilde L_T^2(\cB^{\frac32,\frac12})}
 +\|\nabla U\|_{\widetilde L_T^2(\cB^{\frac12,\frac12})}+\|\partial_1U_1\|_{L_T^1(\cB^{\frac12,\frac12})}
 +\|\nabla U_2\|_{L_T^1(\cB^{\frac12,\frac12})}
 \lesssim\delta_0.
 \end{aligned}
\end{equation}
The estimates for $P_i,N_i$ in \eqref{FQ} follow from
\eqref{app-mixed-product-RT-main}, \eqref{app-mixed-product-RR-low} and \eqref{app-mixed-product-RR-high} with $s=\frac{1}{2}$,
while $E_i$ is estimated as in Appendix~\ref{app:endpoint-product}. Indeed, one has
\begin{equation}\label{zero-mixed-paraproduct-time}
 \begin{aligned}
 &\sum_{i=1}^2\left(
 \|P_i\|_{\widetilde L_T^2(\mathcal X^{\frac12,-\frac12})}
 +\|N_i\|_{\widetilde L_T^2(\mathcal X^{\frac12,-\frac12})}
 \right)\lesssim
 \|h\|_{\widetilde L_T^\infty(\cB^{\frac32,\frac12})}
 \|U\|_{\widetilde L_T^2(\mathcal X^{\frac12,-\frac12})},\\
 &\sum_{i=1}^2
 \|E_i\|_{\widetilde L_T^2(\mathcal X^{-\frac12,-\frac12})}\lesssim
 \|\nabla U\|_{\widetilde L_T^2(\cB^{\frac12,\frac12})}
 \|h\|_{\widetilde L_T^\infty(\mathcal X^{\frac12,-\frac12})}.
 \end{aligned}
\end{equation}
Substituting \eqref{FQ} into \eqref{dc-U-reconstruction} and using
\eqref{decay-critical-coefficients}, \eqref{zero-mixed-gradient-time}
and \eqref{zero-mixed-paraproduct-time}, we obtain
\begin{equation}\label{zero-mixed-velocity-time}
 \begin{aligned}
 \|U\|_{\widetilde L_T^2(\mathcal X^{\frac12,-\frac12})}
 &\lesssim
 \|h\|_{\widetilde L_T^2(\mathcal X^{-\frac12,\frac12})}
 +\delta_0
 \|h\|_{\widetilde L_T^\infty(\mathcal X^{\frac12,-\frac12})}.
 \end{aligned}
\end{equation}
For $G_2$ in \eqref{def:g1g2}, we retain the horizontal derivative
on $h$ in $\partial_1h\,\partial_2U_1$.
Since Bernstein's inequalities imply
$\|\partial_2h\|_{\mathcal X^{-\frac12,-\frac12}}
 \lesssim\|h\|_{\mathcal X^{\frac12,-\frac12}}$,
$\|\partial_1h\|_{\mathcal X^{-\frac12,-\frac12}}
 \lesssim\|h\|_{\mathcal X^{-\frac12,\frac12}}$ and
$\|h\|_{\mathcal X^{-\frac12,\frac12}}\lesssim\|h\|_{\mathcal X^{\frac12,-\frac12}}$, applying the estimates used in \eqref{zero-mixed-paraproduct-time}
to $N_i$, $E_1$ and $E_2$, we arrive at
\begin{align}
 \sum_{i=1}^2\|N_i\|_{\widetilde L_T^1(\mathcal X^{\frac12,-\frac12})}
 &\lesssim
 \|h\|_{\widetilde L_T^2(\cB^{\frac{1}{2},\frac{3}{2}})}
 \|U\|_{\widetilde L_T^2(\mathcal X^{\frac12,-\frac12})},
 \label{zero-mixed-N-time}\\
 \|E_1\|_{\widetilde L_T^1(\mathcal X^{-\frac12,-\frac12})}
 &\lesssim
 \|\partial_1U_1\|_{L_T^1(\cB^{\frac12,\frac12})}
 \|h\|_{\widetilde L_T^\infty(\mathcal X^{\frac12,-\frac12})}+
 \|\partial_2U_1\|_{\widetilde L_T^2(\cB^{\frac12,\frac12})}
 \|h\|_{\widetilde L_T^2(\mathcal X^{-\frac12,\frac12})},
 \label{zero-mixed-E1-time}\\
 \|E_2\|_{\widetilde L_T^1(\mathcal X^{-\frac12,-\frac12})}
 &\lesssim
 \|\nabla U_2\|_{L_T^1(\cB^{\frac12,\frac12})}
 \|h\|_{\widetilde L_T^\infty(\mathcal X^{\frac12,-\frac12})}.
 \label{zero-mixed-E2-time}
\end{align}
Combining \eqref{def:g1g2}, \eqref{zero-mixed-paraproduct-time}
and \eqref{zero-mixed-N-time}--\eqref{zero-mixed-E2-time}
with \eqref{decay-critical-coefficients},
\eqref{zero-mixed-gradient-time} and \eqref{zero-mixed-velocity-time} leads to
\begin{align}
 \|G_1\|_{\widetilde L_T^2(\mathcal X^{\frac12,-\frac12})}
 &\lesssim
 \delta_0\|h\|_{\widetilde L_T^2(\mathcal X^{-\frac12,\frac12})}
 +\delta_0^2\|h\|_{\widetilde L_T^\infty(\mathcal X^{\frac12,-\frac12})},
 \label{zero-mixed-G1-time}\\
 \|G_2\|_{\widetilde L_T^1(\mathcal X^{\frac12,-\frac12})}
 &\lesssim\delta_0\left(
 \|h\|_{\widetilde L_T^\infty(\mathcal X^{\frac12,-\frac12})}
 +\|h\|_{\widetilde L_T^2(\mathcal X^{-\frac12,\frac12})}
 \right).
 \label{zero-mixed-G2-time}
\end{align}
Since $\delta_0\leq1$, \eqref{zero-mixed-G1-time} and
\eqref{zero-mixed-G2-time} imply \eqref{zero-mixed-source}.

To prove \eqref{decay-mixed-source}, we estimate $F,Q$ in
\eqref{dc-errors} and use \eqref{dc-G-equation}.
Arguing as in \eqref{zero-mixed-velocity-time}, with the radial exponent $5/2$ and using \eqref{zero-mixed-gradient-time}, we obtain
\begin{equation}\label{decay-mixed-velocity-time}
 \|U\|_{\widetilde L_T^2(\mathcal X^{\frac 52,-\frac 12})}
 \lesssim\|h\|_{\widetilde L_T^2(\mathcal X^{\frac32,\frac12})}
 +\delta_0\|h\|_{\widetilde L_T^\infty(\mathcal X^{\frac 52,-\frac 12})}.
\end{equation}
Applying \eqref{decay-mixed-product} with time exponents $(2,2)$
and using \eqref{decay-mixed-velocity-time}, we derive
\[
 \begin{aligned}
 \|\partial_1h\,\partial_2U_1
           \|_{\widetilde L_T^1(\mathcal X^{\frac{3}{2},-\frac{1}{2}})}&\lesssim
 \|h\|_{\widetilde L_T^2(\cB^{\frac 12,\frac 32})}
 \|U_1\|_{\widetilde L_T^2(\mathcal X^{\frac 52,-\frac 12})}+
 \|\partial_2U_1\|_{\widetilde L_T^2(\cB^{\frac 12,\frac 12})}
 \|h\|_{\widetilde L_T^2(\mathcal X^{\frac32,\frac12})}\\
 &\lesssim\delta_0\bigl(
 \|h\|_{\widetilde L_T^\infty(\mathcal X^{\frac 52,-\frac 12})}
 +\|h\|_{\widetilde L_T^2(\mathcal X^{\frac32,\frac12})}\bigr).
 \end{aligned}
\]

The other terms in \eqref{dc-errors} are estimated by
\eqref{decay-mixed-product} with time exponents $(\infty,1)$.
As in the proof of \eqref{dc-good-gradient}, we can employ \eqref{dc-hodge}
and the smallness in \eqref{decay-critical-coefficients}
to absorb the terms involving $\partial_1U_1$ and $\nabla U_2$.
Then, \eqref{dc-G-equation} implies \eqref{decay-mixed-source}, which finishes the proof of Lemma~\ref{lem:mixed-source-estimates}.
\end{proof}

\bigbreak

\bigbreak
\noindent \textbf{Acknowledgments.} 
The authors would like to thank Valentin Lemarié and Marius Paicu for stimulating discussions and valuable insights related to this work.

RB is supported by the GNAMPA group of INdAM. 
TCB is supported by the project ANR-24-CE40-3260 -- Hyperbolic Equations, Approximations \& Dynamics (HEAD) and the project ANR-25-CE40-5565 (Cookie). LYS is supported by the National Natural Science Foundation of China (12671258). SW is supported in part by the Spanish Ministry of Science and Innovation through the ``Severo Ochoa Programme for Centres of Excellence in R\&D'' (CEX2023-001347-S).

\bigbreak
\noindent\textbf{Manuscript history and declaration of AI use.}
This work improves on a previous result by the first and second authors, together with Marius Paicu, published in ARMA in 2024, where stability was proved in $H^{3+}(\mathbb{R}^2)$. The present project, aimed at lowering the regularity threshold to $H^{2+}$, started in 2024 and formed part of Shan Wang's postdoctoral research. While the anisotropic functional framework is the same as in the 2024 paper, the present work introduces two key new ingredients: the use of isopycnal coordinates, previously employed by the first author in joint work with V.~Duch\^ene published in CPDE in 2024, and a new commutator estimate developed by Shan Wang in 2024. AI tools played no role in the mathematical research and results. They were used only for limited editorial purposes, namely to shorten some lengthy sentences, mainly in the introduction and appendix, and to polish the figure.

\smallbreak
\noindent \textbf{Data availability statement.}
Data sharing is not applicable to this article, as no datasets were generated or analyzed during the current study.

\smallbreak
\noindent \textbf{Conflicts of interest.}
The authors have no competing interests to declare that are relevant to the content of this article.

\begin{small}

\end{small}

\vspace{10mm}

\end{document}